\documentclass[a4paper, 10pt, american]{amsart}
\usepackage{times,latexsym,amssymb}
\usepackage{amsmath,amsthm,bm}
\usepackage{color}
\usepackage[colorlinks,pdfpagelabels,pdfstartview = FitH,bookmarksopen= true,bookmarksnumbered = true,linkcolor = blue,plainpages = false,hypertexnames = false,citecolor = red,pagebackref=false]{hyperref}

\usepackage{amsbsy}
\usepackage{amstext}
\usepackage{amssymb}
\usepackage{esint}
\usepackage{stmaryrd}
\usepackage[pdftex]{graphicx}
\usepackage{floatflt}
\usepackage{appendix}
\usepackage{mathrsfs}
\usepackage{cancel}

\allowdisplaybreaks
\newtheorem{theorem}{Theorem}
\newtheorem{lemma}[theorem]{Lemma}
\newtheorem{definition}[theorem]{Definition}
\newtheorem{proposition}[theorem]{Proposition}
\newtheorem{corollary}[theorem]{Corollary}
\newtheorem{remark}[theorem]{Remark}
\numberwithin{theorem}{section}
\numberwithin{equation}{section}

\newcommand{\mint}{- \mskip-19,5mu \int}

\def\N{\mathbb{N}}
\def\R{\mathbb{R}}

\renewcommand{\d}{\mathrm{d}}
\newcommand{\dx}{\mathrm{d}x}
\newcommand{\dy}{\mathrm{d}y}

\newcommand{\dt}{\mathrm{d}t}
\newcommand{\ds}{\mathrm{d}s}

\renewcommand{\epsilon}{\varepsilon}
\newcommand{\radius}{{\mathfrak r}}

\DeclareMathOperator{\Div}{div}

\DeclareMathOperator{\loc}{loc}
\renewcommand{\epsilon}{\varepsilon}
\newcommand{\eps}{\varepsilon}
\renewcommand{\rho}{\varrho}

\def\eqn#1$$#2$${\begin{equation}\label#1#2\end{equation}}

\newcommand{\pl}{\partial}

\newcommand{\dsty}{\displaystyle}
\newcommand{\al}{\alpha}

\newcommand{\dl}{\delta}

\newcommand{\lm}{\lambda}

\newcommand{\varep}{\varepsilon}
\newcommand{\vp}{\varphi}

\newcommand{\Om}{\Omega}

\newcommand{\z}{\zeta}

\newcommand{\bom}{\boldsymbol{\mathfrak m}}

\let\TeXchi\chi
\newbox\chibox
\setbox0 \hbox{\mathsurround0pt $\TeXchi$}
\setbox\chibox \hbox{\raise\dp0 \box 0 }
\def\chi{\copy\chibox}

\newcommand{\power}[2]{\bm{#1^{\mbox{\unboldmath{\scriptsize$#2$}}}}}
\newcommand{\powerexp}[2]{\bm{#1^{\mbox{\unboldmath{\tiny$#2$}}}}}
\newcommand{\abs}[1]{|#1|}
\newcommand{\babs}[1]{\big|#1\big|}

\def\Xint#1{\mathchoice
    {\XXint\displaystyle\textstyle{#1}}%
    {\XXint\textstyle\scriptstyle{#1}}%
    {\XXint\scriptstyle\scriptscriptstyle{#1}}%
    {\XXint\scriptscriptstyle\scriptscriptstyle{#1}}%
    \!\int}
\def\XXint#1#2#3{\setbox0=\hbox{$#1{#2#3}{\int}$}
    \vcenter{\hbox{$#2#3$}}\kern-0.5\wd0}
\def\bint{\Xint-}
\def\dashint{\Xint{\raise4pt\hbox to7pt{\hrulefill}}}

\def\Xiint#1{\mathchoice
    {\XXiint\displaystyle\textstyle{#1}}%
    {\XXiint\textstyle\scriptstyle{#1}}%
    {\XXiint\scriptstyle\scriptscriptstyle{#1}}%
    {\XXiint\scriptscriptstyle\scriptscriptstyle{#1}}%
    \!\iint}
\def\XXiint#1#2#3{\setbox0=\hbox{$#1{#2#3}{\iint}$}
    \vcenter{\hbox{$#2#3$}}\kern-0.5\wd0}
\def\biint{\Xiint{-\!-}}

\subjclass[2020]{35B65, 35K67, 35K40, 35K55}

\begin{document}
\title[Local boundedness and higher integrability]{Local boundedness and higher integrability for the sub-critical singular porous medium system}
\date{\today}

\author[V. B\"ogelein]{Verena B\"{o}gelein}
\address{Verena B\"ogelein\\
Fachbereich Mathematik, Universit\"at Salzburg\\
Hellbrunner Str. 34, 5020 Salzburg, Austria}
\email{verena.boegelein@plus.ac.at}

\author[F. Duzaar]{Frank Duzaar}
\address{Frank Duzaar\\
Fachbereich Mathematik, Universit\"at Salzburg\\
Hellbrunner Str. 34, 5020 Salzburg, Austria}
\email{frankjohannes.duzaar@plus.ac.at}

\author[U. Gianazza]{Ugo Gianazza}
\address{Ugo Gianazza\\
Dipartimento di Matematica ``F. Casorati",
Universit\`a di Pavia\\ 
via Ferrata 5, 27100 Pavia, Italy}
\email{ugogia04@unipv.it}

\author[N. Liao]{Naian Liao}
\address{Naian Liao\\
Fachbereich Mathematik, Universit\"at Salzburg\\
Hellbrunner Str. 34, 5020 Salzburg, Austria}
\email{naian.liao@plus.ac.at}

\begin{abstract}
The gradient of weak solutions to porous medium-type equations or systems possesses a higher integrability than the one assumed in the pure notion of a solution. We settle the critical and sub-critical, singular case and complete the program.
\end{abstract}

\maketitle

%****************************************************************

\setcounter{tocdepth}{1}
\tableofcontents

\section{Introduction and main results}
In this manuscript we finally bring to a conclusion the long-standing quest for higher integrability of the gradient of weak solutions to porous medium-type equations or systems:
\begin{equation}\label{por-med-eq}
  \partial_t u -\Div \boldsymbol A(x,t,u,D\power{u}{m}) =\Div F
  \quad\mbox{in $\Omega_T$,}
\end{equation}
where $m>0$, $\Omega\subset \R^N$ is an open bounded set, $\Omega_T:=\Omega\times(0,T]$, and $\boldsymbol A\colon \Omega_T\times\R^k\times\R^{kN}\to \R^{kN}$ is a Carath\'eodory vector field satisfying  
\begin{equation}\label{growth}
\left\{
\begin{array}{c}
	\boldsymbol A(x,t,u,\xi)\cdot\xi\ge \nu|\xi|^2\, ,\\[6pt]
	| \boldsymbol A(x,t,u,\xi)|\le L|\xi|\, ,\\[6pt]
    \dsty\sum_{\al=1}^N \boldsymbol A_\alpha(x,t,u,\xi)\cdot\power{u}{m}\, \xi_\al\cdot\power{u}{m}\ge 0,
\end{array}
\right.
\end{equation}
for a.e.~$(x,t)\in \Omega_T$ and any $(u,\xi)\in \R^k\times\R^{kN}$.
Note that~\eqref{growth}$_3$ asserts some sort of diagonal structure and is superfluous in the case of a single equation; see Remark~\ref{rem:diagonal} below.

Recent advances accessed this difficult subject while left open the critical and subcritical singular case $m\in(0,\frac{N-2}{N+2}]$. Here, our task is to settle this case and complete the picture. Indeed, our main result reads as follows:

\begin{theorem}\label{thm:higherint}
Assume that $N\ge3$,
$$
	0<m\le m_c:=\frac{N-2}{N+2},\quad p>\frac{N+2}2,
$$
and $r>0$ that satisfies
$$
\boldsymbol \lambda_r := N(m-1)+2r>0.$$
Then, there exists $\eps_o{=}\eps_o(N,m,\nu,L, p,r){\in} (0,1]$
such that whenever $F\in L^{2p}_{\rm loc}(\Omega_T,\R^{kN})$, $u$
is a weak solution to the porous medium-type system~\eqref{por-med-eq} in the sense of Definition~\ref{def:weak_solution}
under the assumptions~\eqref{growth}, and $u\in L^r_{\loc}(\Omega_T,\R^k)$, then  we have
$$
  	D\power{u}{m} 
	\in 
	L^{2+2\eps_o}_{\rm loc}\big(\Omega_T,\R^{kN}\big).
$$
Moreover, for every $\eps\in(0,\epsilon_o]$ and every cylinder 
$Q_R(z_o):=B_R(x_o)\times\big(t_o-R^{\frac{1+m}{m}},t_o+R^{\frac{1+m}{m}}\big)$ with 
$
  Q_{2R}(z_o)\subseteq\Omega_T
$,
we have the quantitative local higher integrability estimate
\begin{align}\label{eq:higher-int}
	\biint_{Q_{R}} &
	|D\power{u}{m}|^{2+2\epsilon} \,\dx\dt \nonumber\\
    &\qquad\le
	C\, \mathcal U(2R)^{\epsilon d}
	\biint_{Q_{2R}}
	|D\power{u}{m}|^{2} \,\dx\dt +
	C\bigg[ \biint_{Q_{2R}}  |F|^{2p} \,\dx\dt \bigg]^{\frac{1+\epsilon}{p}},
\end{align}
where 
\begin{equation*}
    \mathcal U(2R)
    :=
    1+\bigg[ \biint_{Q_{2R}} 
 	\frac{\abs{u}^{r}}{R^{\frac{r}m}}  \,\dx\dt\bigg]^{\frac{1+m}{r}} +
 	\bigg[ \biint_{Q_{2R}} |F|^{2p}\,\dx\dt \bigg]^{\frac1p}
\end{equation*}
and
\begin{equation}\label{def:d}
  	d
	:=
  	\frac{2r}{\boldsymbol\lambda_r}
\end{equation}
denotes the scaling deficit, and $C=C(N,m,\nu,L, p,r)$.
\end{theorem}

\begin{remark}
\textup{The scaling deficit $d$ that appears in the higher integrability
  estimate reflects the inhomogeneous scaling behavior of the porous
  medium system. 
In particular, we would have  $d=1$ only if
  $m=1$, which corresponds to the case of the classical heat
  equation. On the other hand, our assumption on $r$ ensures that $d$ is finite and strictly larger than 1.}
\end{remark}
\begin{remark} 
\textup{ The above higher integrability result can  be extended to the case of vector-fields $\boldsymbol A$ satisfying  the more general growth and coercivity conditions 
\begin{equation*}
\left\{
\begin{array}{c}
	\boldsymbol A(x,t,u,\xi)\cdot\xi\ge \nu|\xi|^2-h_1\, ,\\[6pt]
	| \boldsymbol A(x,t,u,\xi)|\le L|\xi| +h_2,
\end{array}
\right.
\end{equation*}
with non-negative measurable functions $h_1,h_2\colon \Omega_T\to [0,\infty]$, so that
$h_1+ h_2^2\in L^{p}(\Omega_T)$ for the exponent $p>\frac{N+2}2$ from Theorem~\ref{thm:higherint}.
}\end{remark}

\begin{remark}\label{rem:diagonal}\upshape
    Conditions \eqref{growth}$_1$--\eqref{growth}$_2$ are quite standard in the field. Concerning \eqref{growth}$_3$, we discuss its significance in Remark~\ref{rem:C_o} in the Appendix.
In the scalar case, i.e. when $k=1$,  condition \eqref{growth}$_3$ is redundant. 
However, for systems it is well-known that the first two  conditions in \eqref{growth}
%$_3$ 
are insufficient to guarantee local boundedness of solutions, which is what we need in the following (see Section~\ref{SS:sign} below).
\end{remark}

The  quantitative local estimate \eqref{eq:higher-int} can easily be converted into an estimate on standard
parabolic cylinders $C_R(z_o):=B_{R}(x_o)\times(t_o-R^2,t_o+R^2)$. The precise statement is:

\begin{corollary}\label{cor:higher-int}
  Under the assumptions of Theorem \ref{thm:higherint},    on any cylinder $C_{2R}(z_o)\subseteq\Omega_T$  and for every $\eps\in(0,\eps_o]$ we have
   \begin{align*}
	\biint_{C_{ R}(z_o)} &
    |D\power{u}{m}|^{2+2\epsilon} \,\dx\dt \\
	&\le
	\frac{C}{R^{2\epsilon}}\,\mathcal U(2R)^{\epsilon d}
	\biint_{C_{2R}(z_o)}
	|D\power{u}{m}|^{2} \,\dx\dt +
	C
    \bigg[
	\biint_{C_{2R}(z_o)} |F|^{2p} \,\dx\dt\bigg]^\frac{1+\epsilon}{p} .
   \end{align*}
where 
$$
    \mathcal U(2R)
    :=
    1+\bigg[ \biint_{C_{2R}(z_o)} \abs{u}^{r} \,\dx\dt\bigg]^\frac{1+m}{r}
    +
    R^2\bigg[\biint_{C_{2R}(z_o)}
    |F|^{2p} \,\dx\dt \bigg]^\frac1p
$$
and $C=C(N,m,\nu,L, p,r)$. 
\end{corollary}
%%%%%%%%%%%%%%

%%%%%%%%%%%%%%%%

System \eqref{por-med-eq} has a different behavior when $m>1$ or $0<m<1$. The first case is called
\emph{slow diffusion} regime, since disturbances propagate with finite speed and free boundaries may occur, while
in the second case disturbances propagate with infinite speed and 
extinction in finite time may occur; this is usually referred to as \emph{fast diffusion} range. 

For more information on the general theory of the porous medium equation and related regularity results, we refer to \cite{CaVaWo, DiBenedetto_Holder} 
and the references therein.
Here, we focus on the fast diffusion range; besides their intrinsic mathematical interest, equations and systems with $m\in(0,1)$ arise in a number of different applications, some of which we discuss next.

\subsection{Models}
The diffusion of a hydrogen plasma across a purely poloidal octupole magnetic field, as discussed in \cite{berryman-77,berryman-78}, may be written in normalized units as
\begin{equation}\label{Eq:1bis:2:1}
\frac{\partial}{\partial x}\left[D(n)\frac{\partial n}{\partial x}\right]=F(x)\frac{\partial n}{\partial t},\qquad
0\le x\le 1,
\end{equation}
where $n$ is the density of the plasma, $x$ is the spatial variable corresponding to magnetic field potential and $t$ is the time.
The geometric factor $F(x)>0$ is a function determined by the octupole geometry, and, in general, $D(n)\approx n^\dl$, where $\dl\ge-1$.

In \cite{Car} Carleman proposed an interesting singular porous medium model to study the diffusive limit of kinetic equations. He considered two types of particles in a one-dimensional setting, moving with opposite speeds, that is $c$ and $-c$. If we denote the densities respectively with $u$ and $v$, we can write their simple dynamics as
\begin{equation*}
\left\{
\begin{aligned}
&\partial_t u+c\partial_x u=k(u,v)(v-u),\\
&\partial_t v-c\partial_x v=k(u,v)(u-v),
\end{aligned}
\right.
\end{equation*}
where $k(u,v)\ge0$ is a proper interaction kernel. In a typical case, one assumes $\dsty k=(u+v)^\al c^2$
with $\al\in(0,1)$. If we now write the equations for $\rho=u+v$ and $j=c(u-v)$ and pass to the limit as 
$c=\frac1\eps\to\infty$, we will obtain to the first order in powers of $\eps$ (see \cite{lions-toscani})
\begin{equation*}
\partial_t\rho=\frac12\partial_x\left(\frac1{\rho^\al}\partial_x\rho\right);
\end{equation*}
once we rewrite the right-hand side as $\frac1{2(1-\alpha)}(\rho^{1-\alpha})_{xx}$, this is exactly the model case of \eqref{por-med-eq} with $m=1-\al\in(0,1)$. In Carleman's model $\alpha=1$, which yields $m=0$, is also allowed, but we do not take this possibility into account here.

Suppose we are given a Riemannian manifold $(M,g_o)$ in space dimension $N\ge3$. We look for a one-parameter family $g=g(x,t)$ of metrics solution to the evolution problem
\begin{equation*}
\left\{
\begin{aligned}
&\partial_t g=-Rg\\
&g(0)=g_o
\end{aligned}
\right.
\end{equation*}
on $M$, where $R=R_g$ denotes the scalar curvature of the metric $g$. If we write $\Delta_o$
for the Laplace-Beltrami operator of $g_o$, it is not hard to show that, after
rescaling the time variable, the problem is equivalent to
 \begin{equation*}
\left\{
\begin{aligned}
&\partial_t (u^{\frac{N+2}{N-2}})=Lu\\
&u(0)=1
\end{aligned}
\right.
\end{equation*}
on $M$, where 
\[
Lu:=\kappa\Delta_o u-R_{g_o}u,\qquad\kappa=\frac{4(N-1)}{N-2}.
\]
Let now $(M,g_o)$ be $\R^N$ with the standard flat metric, so that $R_{g_o}=0$. Putting
$\dsty u^{\frac{N+2}{N-2}}=v$, and properly rescaling the time, yields the so-called \emph{Yamabe flow} in $\R^N$:
\begin{equation*}%\label{Eq:1bis:5:1}
\partial_t  v=\Delta v^{m_c},\qquad m_c=\frac{N-2}{N+2}.
\end{equation*}

Rosenau \cite{rosenau} has pointed out that the conduction of heat in metallic and ceramic materials gives rise to a singular porous medium equation, very similar to \eqref{Eq:1bis:2:1}.

Porous medium \emph{systems} can be derived from a time-dependent $p$-curl system which describes Bean’s critical-state model in the superconductivity theory \cite{Bean}. In particular, one ends up with the system
\[
\partial_t U - \Delta(|U|^{m-1} U)=0,\qquad U=(u_1,u_2,\dots,u_k),
\]
where $U=\Delta\psi$, and $\psi$ is the magnetic vectorial potential, which is assumed to have compact support and satisfy the gauge calibration $\Div\psi=0$. It must be said that, in this case, for physically significant models, one generally assumes $m>1$.
For mathematical results on this type of problem, see \cite{Yin}, and the references therein.

Just for the sake of completeness, let us mention that systems of porous medium-type having a different structure with respect to the one we deal with here, are studied in \cite{ros-hym} (see also the references therein).

\subsection{Novelty and Significance}\label{SS:sign}
With regard to the higher integrability of the gradient of a weak solution, we are referring to the property of the gradient to enjoy a slightly larger degree of summability than the one that comes from the sheer definition. Such a fact is interesting in itself, and at the same time has a number of important consequences. The first step in this field is represented by Gehring's paper \cite{Gehring}, which was then extended to elliptic equations in \cite{Meyers-Elcrat}. An inaccuracy in Gehring's work about the support of the function was independently fixed in \cite{Stred} and \cite{Giaq-Mod} (see also \cite[Chapter~V, Theorem~2.1]{Giaquinta:book} and \cite[Section~6.4]{Giusti:book}). Since these first works, the higher integrability property for solutions to elliptic problems has been widely extended in many different ways.

In the parabolic setting, the first result is due to Giaquinta \& Struwe \cite{Giaquinta-Struwe} for quasi-linear systems with bounded and measurable coefficients. For quite some time there were no substantial progresses as far as the growth condition is concerned, since all generalizations dealt just with a growth of order $2$. A breakthrough was realized by Kinnunen \& Lewis \cite{Kinnunen-Lewis:1} in 2000, when they proved the higher integrability of the gradient of weak solutions to a wide class of non-linear parabolic systems of $p$-Laplacian type, provided that $p>\frac{2N}{N+2}$.

For porous medium equations and systems the higher integrability of the gradient $Du^m$, which is the natural quantity in this context rather than $Du$, remained elusive for another twenty years, until the works of Gianazza \& Schwarzacher, and of B\"ogelein, Duzaar, Korte \& Scheven.

Gianazza \& Schwarzacher \cite{Gianazza-Schwarzacher,Gianazza-Schwarzacher:m<1} dealt only with non-negative solutions to equations, and this is largely due to their use of technical tools which are available only for this kind of solutions. On the other hand, B\"ogelein, Duzaar, Korte \& Scheven \cite{BDKS-higher-int,BDS-sing} simplified the proof and extended it to signed solutions of equations, and also to systems.

All these results have in common the assumption $m>\frac{(N-2)_+}{N+2}$, and this is not surprising, since such a lower bound on $m$ somehow corresponds to the lower bound on $p$ found by Kinnunen \& Lewis. It is well known that when $p\in(1,\frac{2N}{N+2}]$ for the parabolic $p$-Laplacian, and $m\in(0,\frac{(N-2)_+}{N+2}]$ for the porous medium equation, the behavior of the solutions changes dramatically, in such a way that a complete and satisfactory theory is still missing (from this point of view, see also \cite{DiBe}, and \cite[Appendices~A and B]{DBGV-book}).

This caution notwithstanding, here we give a positive answer
to the higher integrability issue precisely in the so-called \emph{critical and sub-critical} fast diffusion range for $m$, that is, when $m$ satisfies 
\begin{equation}\label{lower-bound}
0<m\le m_c:=\frac{(N-2)_+}{N+2}.
\end{equation}
In this way, the description is complete, as far as porous medium equations and systems are concerned. The corresponding theory for the parabolic $p$-Laplacian will be considered in a separate work.

A careful check of the proofs, both in \cite{BDS-sing} and in \cite{Gianazza-Schwarzacher:m<1} shows that the lower bound $\frac{(N-2)_+}{N+2}$ for $m$ is a consequence of the subtle and delicate interplay between energy estimates and Sobolev-Poincar\'e inequalities, which is needed to prove a reverse H\"older inequality, the key step in any higher integrability result.
On the other hand, for the construction of sub-intrinsic cylinders (see~\S\,\ref{sec:cylinders}) one needs $u\in L^r_{\loc}(\Omega_T,\R^k)$ with $r$ as large as possible, but to realize that the sheer notion of a solution does not allow $m$ to go below $m_c$. 
Hence, the key idea in the subcritical and critical range is to assume proper summabilities on $u$ and on $F$, higher than the ones required by Definition~{\upshape\ref{def:weak_solution}}, so that $u$ eventually turns out to be bounded. That is why we require $u\in L^r_{\loc}(\Omega_T,\R^k)$ with $r$ such that $N(m-1)+2r>0$, and $F\in L^{2p}_{\loc}(\Omega_T,\R^{kN})$ with $p>\frac{N+2}2$. We show that these assumptions yield the boundedness of $u$; first, we give a qualitative estimate using Moser's techniques, and then, we refine it in a quantitative fashion using De Giorgi's approach, according to what we need in the sequel of the proof of Theorem~\ref{thm:higherint}. This is definitely one of the novelties of this work. Indeed, boundedness estimates for solutions to \emph{homogeneous} equations when $m$ satisfies \eqref{lower-bound} are known (see \cite[Appendix~B, Section~B4]{DBGV-book} and \cite[Theorems~2.1 and 2.4]{bonforte}), but to our knowledge estimates for systems and/or with a right-hand side different from zero are new.

It might look as if these hypotheses on $u$ and $F$ were ad-hoc, specifically linked to the methods we employ, and that a different approach would yield higher integrability results under less restrictive assumptions. However, in Section~\ref{SS:sharp} we show that they are somehow optimal, in the sense explained below.

\subsection{Sharpness of the result}\label{SS:sharp}
Unbounded solutions do exist for $m$ as in \eqref{lower-bound}, and the integrability of their gradient in general does not possess a self-improving property.
Concerning boundedness, the function
\begin{equation*}%\label{Eq:4:18:10}
\begin{array}{c}
{\dsty u(x,t)=\Big(\frac{2m|\lm|}{1-m}\Big)^{\frac1{1-m}}
\frac{(T-t)_+^{\frac{1}{1-m}}}{|x|^{\frac{2}{1-m}}}},\\
{}\\
{\dsty 0<m<\frac{(N-2)_+}{N+2},\qquad
\lm=N(m-1)+2.}
\end{array}
\end{equation*}
is a non-negative, local, weak solution 
to the prototype scalar porous medium equation 
\begin{equation}\label{PME} 
  \partial_t u - \Delta\, u^{m} = \Div F,
\end{equation}
in $\R^N\times\R$ for $F=0$. 
Such a solution is unbounded near $x=0$ for all $t<T$ and 
finite otherwise. It is known that when $m>\frac{(N-2)_+}{N+2}$, local boundedness is inherent in the notion of weak solutions. Hence, this example shows such a result is almost sharp.

A direct computation shows that this solution satisfies
\[
\nabla u^m \in L^{2+\varepsilon}_{\loc}(\R^N\times\R),\quad \forall \, \varep\in(0,\varepsilon_o),\quad \varepsilon_o=-\frac{N(m-1)+2(1+m)}{1+m}.
\]
However, we note that $\varep_o\downarrow0$ as $m\uparrow m_c$. In general, unbounded solutions cannot have the above improvement of the integrability with an $\varep_o$ that remains away from $0$ as $m\uparrow m_c$.
To wit, let us consider the porous medium equation \eqref{PME} with $F=0$. Suppose $u$ is a non-negative, local, weak solution satisfying $\nabla u^m\in L^{2+\varepsilon}_{\loc}(\Om_T)$ for some $\varepsilon$ that can be determined a priori and satisfies $\varepsilon\in[\varepsilon_*,1]$ for some $\varepsilon_*>0$ and all $m\in(0,\frac{N-2}{N+2}]$. Then, by the embedding theorem, $u\in L^{r}_{\loc}(\Om_T)$ for $r=(2+\varep)\frac{mN+1+m}{N}$. Now, by some elementary algebra one can determine a positive $\dl=\dl(N,\varepsilon_*)$ such that whenever $m\in(\frac{N-2}{N+2}-\dl, \frac{N-2}{N+2}]$, we have $N(m-1)+2r>0$. In fact, denoting $m_c=\frac{N-2}{N+2}$ and using $m_c-\dl<m\le m_c$ we compute
\begin{align*}
    N(m-1)+2r &> N(m_c -1 -\dl) + 2(2+\varepsilon)\frac{(N+1)(m_c-\dl)+1}{N}\\
    &=\underbrace{N (m_c-1)+4\frac{(N+1)m_c+1}{N}}_{=0}\\
    &\quad  + 2\varepsilon\frac{(N+1)m_c+1}{N} -\dl\bigg[N+2(2+\varepsilon)\frac{N+1}{N}\bigg].
\end{align*}
To make the above quantity positive, we can, for instance, choose
\[
\dl=\frac{\varepsilon_*\frac{(N+1)m_c+1}{N}}{N+6\frac{N+1}{N}}.
\]
As a result of  Proposition~\ref{prop:sup-est-qual}, such a solution should be locally bounded in $\Om_T$. 

In the same range for $m$, again for $F=0$, another non-negative, local, weak solution 
to the prototype scalar porous medium equation \eqref{PME} 
in $\R^N\times\R$ 
is given by the function
\[
{\dsty u(x,t)=(T-t)_+^{\frac{N-2-2m}{2m^2}}\left[A|x|^{\frac{(N-2)(1-m)}{m}}+
\frac{1-m}{2m|\lambda|}|x|^2(T-t)_+^{\frac{|\lambda|}{2m^2}}\right]^{-\frac1{1-m}}},
\]
where $A$ is an arbitrary, non-negative real parameter (see \cite[Section~4]{King} and \cite[Proposition~7]{kosov}). We remark that, as before, $u$ is unbounded near $x=0$ for all $t<T$ and finite otherwise. Moreover, if $A=0$, $u$ reduces to the previous example.

As for the self-improving property, consider the following example, built relying on the results of Kosov \& Semenov \cite[Sections 1 and 3]{kosov}. Given $T>0$, in the set $B_R\times(0,\frac T2]$ with $R=\left(\frac T2\right)^{\frac{1}{2m}}$, %and $\lambda=-\frac4{N-2}$, 
the function
\begin{equation*}%\label{Eq:III-counter}
u(x,t)=\left(\frac{(N-2)(T-t)_+}2\right)^{\frac{N+2}4}\frac1{|x|^{\frac{N+2}2}\left(-\ln\left[|x|(T-t)_+^{-\frac{1}{2m}}\right]\right)^{\frac{N+2}4}},
\end{equation*}
is an explicit, non-negative, local, weak solution 
to the equation
\[
\partial_t u-\tfrac1m\Delta\,  u^m=\Div F,
\]
where $m=\frac{N-2}{N+2}$, $N\ge3$, and $F$ is the time-dependent, radial vector field, given by 
\[
F(r,t)=-\frac{N+2}4\left(\frac{N-2}2\right)^{\frac{N+2}4}\left(T-t
\right)_+^{\frac{N-2}4}\frac1{r^{N-1}}\, G(r,t)\,{\mathbf e}_r,
\]
where $r=|x|$, ${\mathbf e}_r=\frac x{|x|}$, and $G(r,t)$ is the %antiderivative 
primitive with respect to $r$ of 
\[
    H(r,t)
    =
    r^{\frac{N-4}2}\left[-2\left(-\ln\left[r(T-t)_+^{-\frac1{2m}}\right]\right)^{-\frac{N-2}4}+\left(-\ln\left[r(T-t)_+^{-\frac1{2m}}\right]\right)^{-\frac{N+2}4}\right].
\]
If we look at the behavior of $F(r,t)$ as $r\downarrow 0$, it is apparent that 
\[
    H(r,t)
    \sim 
    -2\frac{r^{\frac{N-4}2}}{\left(-\ln\left[r(T-t)_+^{-\frac1{2m}}\right]\right)^{\frac{N-2}4}},
\]
whence, by a straightforward application of L'H\^{o}pital's rule, we conclude that 
\[
    G(r,t)
    \sim -\frac 4{N-2} 
    \frac{r^{\frac{N-2}2}}{\left(-\ln\left[r(T-t)_+^{-\frac1{2m}}\right]\right)^{\frac{N-2}4}},
\]
and also
\[
F(r,t)\sim-\frac{N+2}{N-2}\left(\frac{N-2}2\right)^{\frac{N+2}4}\frac{\left(T-t\right)_+^{\frac{N-2}4}}{r^{\frac N2}\left(-\ln\left[r(T-t)_+^{-\frac1{2m}}\right]\right)^{\frac{N-2}4}}\,{\mathbf e}_r.
\]
Provided that $N$ is sufficiently large (namely, $N>4$), we have $F\in L^2_{\loc}(B_R\times(0,\frac T2])$.
Moreover, it is not hard to check that there exists no $\varepsilon>0$ such that $\nabla u^{m}\in L^{2+\varepsilon}_{\loc}$, that is, the self-improving property does not hold in this case.

\color{black}
\subsection{Main technical tools}
We have already discussed the issue of the boundedness of $u$. From a purely technical point of view, the main idea in the proof of Theorem~\ref{thm:higherint} remains the same as in \cite{BDS-sing}, that is, the construction of  suitable intrinsic cylinders $Q_{r,s}(z_o):=B_r(x_o)\times(t_o-s,t_o+s)$ with $z_o=(x_o,t_o)$; however, there is an important difference: now we use cylinders whose space-time scaling depends on the mean values of  $|u|^{r}$. This choice is directly linked to the assumption on the extra summability of $u$.  
The construction of a system of such intrinsic cylinders is far from trivial, and the sub-critical and critical range poses extra difficulties.

Moreover, assuming $u\in L^r_{\loc}(\Omega_T,\R^k)$ allows us to use a technical trick, that was originally employed by  Saari \& Schwarzacher in \cite{SaSch} to overcome the upper bound in the admissible range of parameter found by B\"ogelein, Duzaar, Kinnunen \& Scheven \cite{BDKS-doubly}, when they studied the higher integrability of the gradient of solutions to the so-called \emph{Trudinger equation}.
The new scaling for cylinders and the use of the idea first employed in \cite{SaSch} represent another novelty in the work.

\subsection{Structure of the paper}
In Section~\ref{S:2} we collect the basic notation, and the definition of solution, while Section~\ref{S:3} is devoted to some auxiliary results which will be used throughout the whole paper. In Section~\ref{S:4} we deal with qualitative and quantitative boundeness of solutions, and Section~\ref{sec:energy} contains the basic energy estimates. In Section~\ref{sec:poin} we study the Sobolev-Poincar\'e inequalities that are the fundamental tool needed in Section~\ref{sec:revholder} to prove proper reverse H\"older inequalities, and Section~\ref{sec:hi} is devoted to the proofs of the main results. Finally, in the Appendix we give some estimates which are the starting point for the boundedness results of Section~\ref{S:4}.
\vskip.2truecm
\noindent
{\bf Acknowledgments.} V.~B\"ogelein has been supported by the FWF-Project P31956-N32
``Doubly nonlinear evolution equations." 
U. Gianazza is a  member of the Gruppo Nazionale per l'Analisi Matematica, la Probabilità e le loro Applicazioni (GNAMPA) of the Istituto Nazionale di Alta Matematica (INdAM), and was partly supported by GNAMPA-INdAM  Project 2023 ``Regolarità per problemi ellittici e parabolici con crescite non standard" (CUP\_E53C22001930001). N.~Liao has been supported by the FWF-project P36272-N ``On the Stefan type problems." 

\section{Notation and definition}\label{S:2}

\subsection{Notation}
The characteristic function of a set is denoted as $\mathbf 1_{\{...\}}$.
To keep formulations as simple as possible, we define the {\em
  power of a vector} or of a possibly {\it negative number}  by
$$
	\power{u}{\alpha}
	:=
	|u|^{\alpha-1}u,
	\quad\mbox{for $u\in\R^k$ and $\alpha>0$,}
$$
which in the case $u=0$  and $\alpha\in (0,1)$ we interpret as $\power{u}{\alpha}=0$. 
Throughout the paper we write $z_o=(x_o,t_o)\in \R^N\times\R$ for points in space-time. We use space-time cylinders 
\begin{equation}\label{def-Q}
	Q_\rho^{(\theta)}(z_o)
	:=
	B_\rho^{(\theta)}(x_o)\times\Lambda_\rho(t_o),
\end{equation}
where 
\begin{equation*}
  	B_\rho^{(\theta)}(x_o)
  	:=
  	\Big\{x\in\R^N: |x-x_o|<\theta^{\frac{m(m-1)}{1+m}}\rho\Big\}
\end{equation*}
and
\begin{equation*}
  	\Lambda_\rho(t_o)
  	:=
  	\big(t_o-\rho^{\frac{1+m}{m}},t_o+\rho^{\frac{1+m}{m}}\big)
\end{equation*}
with some scaling  parameter $\theta >0$. 
In the case $\theta =1$, we 
omit the parameter $\theta$ in the notation and write
$$
	Q_\rho(z_o)
	:=
	B_\rho(x_o)\times\big(t_o-\rho^{\frac{1+m}{m}},t_o+\rho^{\frac{1+m}{m}}\big)
$$
instead of $Q_\rho^{(1)}(z_o)$. If the center $z_o$ is clear from the context we omit it in the notation.

For the $L^\infty$-estimates in \S\,\ref{S:4} and Appendix\,\ref{sec:app-energy} we use another kind of space-time cylinders, namely
\begin{equation}\label{cy-infty}
Q_{R,S}(z_o)=B_R(x_o)\times(t_o-S,t_o].
\end{equation}
Its parabolic boundary is given by 
$$
    \partial_P Q_{R,S}(z_o)
    :=
    \big(\partial B_R(x_o)\times(t_o-S,t_o]\big)
    \cup \overline{B_R(x_o)}\times\{t_o-S\}.
$$
Once more, if the center $z_o$ is clear from the context, we omit it in the notation. If $S=R^2$, we simply write $Q_R(z_o)$.

For a map $u\in L^1(0,T;L^1(\Omega,\R^k))$ and given measurable sets $A\subset\Omega$  and $E\subset\Omega_T$ with positive Lebesgue measure the slicewise mean $\langle u\rangle_{A}\colon (0,T)\to \R^k$ of $u$ on  $A$ is defined
by
\begin{equation*}
	\langle u\rangle_{A}(t)
	:=
	\mint_{A} u(t)\,\dx,
	\quad\mbox{for a.e.~$t\in(0,T)$,}
\end{equation*}
whereas  the mean value $(u)_{E}\in \R^k$ of $u$ on  $E$ is defined by
\begin{equation*}
	(u)_{E}
	:=
	\biint_{E} u\,\dx\dt.
\end{equation*}
Note that if $u\in C^0((0,T);L^2(\Omega,\R^k))$, the slicewise means are defined for any $t\in (0,T)$. 
If  $A$ is a ball $B_\rho^{(\theta)}(x_o)$,  we write $\langle
u\rangle_{x_o;\rho}^{(\theta)}(t):=\langle
u\rangle_{B_\rho^{(\theta)}(x_o)}(t)$. Similarly,   if $E$ is a
cylinder of the form $Q_\rho^{(\theta)}(z_o)$, we use the shorthand
notation $(u)^{(\theta)}_{z_o;\rho}:=(u)_{Q_\rho^{(\theta)}(z_o)}$.

Finally, given a \emph{scalar} function $v\colon\Omega_T\to\R$, by $\nabla v$ we mean its spatial gradient, whereas for a \emph{vector valued} function $u\colon \Omega_T\to\R^k$, its gradient will be denoted by $D u$.

\subsection{Notion of weak solution}
In the statement of Theorem~\ref{thm:higherint}, we relied on the notion of {\it weak  solution}, which we now formulate. The definition holds for any $m>0$, and for the definition to make sense, it suffices to assume that $F$ has a low degree of integrability.

\begin{definition}\label{def:weak_solution}
\upshape
Let $m>0$ and $\boldsymbol A\colon \Omega_T\times \R^k\times\R^{kN}\to\R^{kN}$ be a vector field satisfying \eqref{growth} and $F\in L^{2}(\Omega_T,\R^{kN})$. 
A function 
\begin{equation*}%\label{spaces}
	u\in C^0 \big((0,T); L^{1+m}(\Omega,\R^k)\big)
	\quad\mbox{with}\quad 
	\power{u}{m}\in L^2\big(0,T;W^{1,2}(\Omega,\R^k)\big)
\end{equation*} 
is a \textit{weak solution} to the porous medium-type system \eqref{por-med-eq} if and only if the identity
\begin{align}\label{weak-solution}
	\iint_{\Omega_T}\big[u\cdot\partial_t\varphi - \boldsymbol A(x,t,u,D\power{u}{m})\cdot D\varphi\big]\dx\dt
    	=
    	\iint_{\Omega_T} F\cdot D\varphi \,\dx\dt
\end{align}
holds true, for any testing function $\varphi\in
C_0^\infty(\Omega_T,\R^k)$. 
\qed
\end{definition}
%%%%%%%%%%%%%%
%\color{black}
\section{Auxiliary material}\label{S:3}
In this section we provide the necessary tools which will be used later.
To ``re-absorb'' certain terms, we shall frequently use the following iteration lemma, cf. \cite[Lemma 6.1]{Giusti:book}.

\begin{lemma}\label{lem:tech-classical}
Let $0<\eps<1$, $A,B\ge 0$ and $\alpha\ge 0 $. Then there exists a universal constant  $C = C(\alpha,\eps)$
such that whenever $\phi\colon [R_o,R_1]\to \R$  is non-negative bounded function satisfying
\begin{equation*}
	\phi(t)
	\le
	\eps \phi(s) + \frac{A}{(s-t)^\alpha}  + B
	\qquad \text{for all $R_o\le t<s\le R_1$,}
\end{equation*}
then we have
\begin{equation*}
	\phi(t)
	\le
	C\,  \bigg[\frac{A}{(s-t)^\alpha}  + B\bigg]
\end{equation*}
for any $R_o\le t<s\le R_1$.
\end{lemma}

The following lemma can be deduced as in \cite[Lemma~8.3]{Giusti:book}.

\begin{lemma}\label{lem:Acerbi-Fusco}
For any $\alpha>0$, there exists a constant $c=c(\alpha)$ such that,
for all $a,b\in\R^k$, $k\in\N$, the following inequality holds true:
\begin{align*}
	\tfrac1c\big|\power{b}{\alpha} - \power{a}{\alpha}\big|
	\le
	\big(|a| + |b|\big)^{\alpha-1}|b-a|
	\le
	c \big|\power{b}{\alpha} - \power{a}{\alpha}\big|.
\end{align*}
\end{lemma}

The next lemma is an immediate consequence of Lemma~\ref{lem:Acerbi-Fusco}.

\begin{lemma}\label{lem:a-b}
For any $\alpha\ge 1$, there exists a constant $c=c(\alpha)$ such that,
for all $a,b\in\R^k$, $k\in\N$, the following inequality holds true:
\begin{align*}
	|b-a|^\alpha
    \le
    c\big|\power{b}{\alpha} - \power{a}{\alpha}\big|.
\end{align*}
\end{lemma}

It is well known that mean values over subsets $A\subset B$ are quasi-minimizers of the mapping $
\R^N\ni a\mapsto \int_B |u-a|^p \,\dx$. The following lemma shows that this also applies to powers 
 $\power{u}{\alpha}$  of $u$, provided  $\alpha\ge\frac1p$.  For $p=2$ and $A=B$, the lemma has been proved in \cite[Lemma 6.2]{Diening-Kaplicky-Schwarzacher}.  The general version is established
in \cite[Lemma~3.5]{BDKS-doubly}.

\begin{lemma}\label{lem:alphalemma}
For any $p\ge 1$ and $\alpha\ge\frac1p$, there exists a universal constant $c=c(\alpha,p)$ such that whenever $A\subset B\subset \R^N$, $N\in\N$, are two bounded domains with positive measure, then for any  $u \in L^{\alpha p}(B,\R^k)$ and any  $a\in\R^k$, we have 
$$
	\mint_B \big|\power{u}{\alpha}-\power{(u)_A}{\alpha}\big|^p \dx 
	\le 
	\frac{c\,|B|}{|A|} 
	\mint_B \big|\power{u}{\alpha}-\power{a}{\alpha}\big|^p \dx.
$$
\end{lemma}

\section{Local boundedness of weak solutions}\label{S:4}
%, $L^r_{\rm loc}$-$L^\infty_{\rm loc}$-estimates}

This section consists of two propositions. The first one is based on Moser's iteration and shows weak solutions are locally bounded if they are sufficiently integrable.  The second one uses De Giorgi's iteration and refines the quantitative $L^{\infty}$-estimate. Throughout this section we use the one-sided cylinders introduced in~\eqref{cy-infty} and write $Q_\rho(z_o)=Q_{\rho,\rho^2}(z_o)$.

\begin{proposition}[qualitative $L^\infty$-estimate]\label{prop:sup-est-qual}
Let $m\in(0,m_c]$, $p>\frac{N+2}2$, and $r\ge 1$ satisfying
\begin{equation}\label{Eq:B:4:1}
    \boldsymbol\lambda_r:=N(m-1)+2r>0.
\end{equation}
Then any local, weak solution to \eqref{por-med-eq} in $\Omega_T$ in the sense of Definition~{\upshape\ref{def:weak_solution}}  with right-hand side $F\in L^{2p}_{\loc}(\Omega_T,\R^{kN})$ that satisfies
$
    u\in L^r_{\rm loc}(\Omega_T,\R^k),
$
is locally bounded in $\Omega_T$. Moreover, there exists a constant $C=C(N,m,\nu,L ,p,r)$ such that on any cylinder $Q_\rho(z_o)$ we have
\begin{equation*}
    \sup_{Q_{\frac12\rho}(z_o)}|u|
    \le
    C
    \bigg[
    \Big(
    1+[F]_{2p,Q_{\rho}(z_o)}^{\frac{2p(N+2)}{2p-(N+2)}}
    \Big)
    \biint_{Q_{\rho}(z_o)}\big[|u|^r+1\big]\,\dx\dt
    \bigg]^{\frac{2}{\boldsymbol\lambda_r}},
\end{equation*}
where
\begin{equation*}
    [F]_{2p,Q_{\rho}(z_o)}
    :=
    \bigg[
    \rho^{2p}\biint_{Q_{\rho}(z_o)}|F|^{2p}\,\dx\dt
    \bigg]^{\frac1{2p}}.
\end{equation*}
\end{proposition}
%%%%%%%%%%%%%%%%%
\begin{proof} 
The proof is divided into two steps. The first step towards the local $L^\infty$-bound of $u$ is an improvement of the energy-type inequality from Lemma~\ref{lem:energy-est}, assuming a certain higher integrability of $u$ in form of $|u|\in L_{\rm loc}^{m+1+2m\alpha}(\Omega_T)$
for some $\alpha \ge 0$. This is done by establishing a quantitative reverse H\"older inequality for $|u|$. In the second step, we iterate this reverse Hölder inequality, starting from the $L^r_{\rm loc}(\Omega_T)$-assumption
and successfully improve the integrability of $|u|$. During the iteration process, we keep track of the constants that occur. In this way, we can pass to the limit and derive the local $L^\infty$-bound.

\textit{Step 1.}
We consider $\alpha \ge 0$ and a cylinder $Q_{R}(z_o)\Subset\Omega_T$.
Moreover, we assume that
\begin{equation*}%\label{u-Lalpha}
    |u|\in L^{m+1+2m\alpha}(Q_{R}(z_o)).
\end{equation*}
In the following we omit the reference point $z_o$ in the notation. In the energy estimate~\eqref{En-Est-Phi} from Lemma~\ref{lem:energy-est} we choose 
\begin{equation}\label{def:phi-alpha-l-k}
    \phi_{\alpha, k,\ell}(s)=
    \left\{
    \begin{array}{cl}
        k^{2m\alpha}, &\mbox{if $0\le s\le k^{2m}$,}\\[5pt]
        s^\alpha,&\mbox{if $k^{2m}<s<\ell^{2m}$,}\\[5pt]
        \ell^{2m\alpha},&\mbox{if $s\ge\ell^{2m}$,}
    \end{array}
    \right.
\end{equation}
with parameters $k\in (0,1]$ and $\ell\in(1,\infty)$. Note that $\phi_{\alpha, k,\ell}$ is a bounded non-decreasing Lipschitz continuous function and that the constant $C_{\phi_{\alpha,k,\ell}}$ from \eqref{condition-Phi} is bounded by $\alpha$. Therefore, $\phi_{\alpha,k,\ell}$ is admissible in the energy inequality \eqref{En-Est-Phi}.
To proceed further we need to compute $v$ in \eqref{En-Est-Phi}. 
Notice that 
\begin{equation*}
    \phi_{\alpha, k,\ell}\big(s^{\frac{2m}{m+1}}\big)
    =
    \left\{
    \begin{array}{cl}
        k^{2m\alpha}, &\mbox{if $0\le s\le k^{m+1}$,}\\[5pt]
        s^{\frac{2m\alpha}{m+1}},&\mbox{if $k^{m+1}<s<\ell^{m+1}$,}\\[5pt]
        \ell^{2m\alpha},&\mbox{if $s\ge\ell^{m+1}$.}
    \end{array}
    \right.
\end{equation*}
%\color{black}
If $|u|\le k$ we have 
\begin{align*}
    v
    &=
    \int_0^{|u|^{m+1}} \phi_{\alpha,k,\ell }\big(s^{\frac{2m}{m+1}}\big)\,\ds
    =
    k^{2m\alpha} |u|^{m+1}.
\end{align*}
Next, if $k<|u|<\ell$, we have
\begin{align*}
     v
    &=
    \int_0^{k^{m+1}} \phi_{\alpha,k,\ell }\big(s^{\frac{2m}{m+1}}\big)\,\ds
    +
    \int_{k^{m+1}}^{|u|^{m+1}} \phi_{\alpha,k,\ell }\big(s^{\frac{2m}{m+1}}\big)\,\ds\\
    &=
    k^{m+1+2m\alpha}+ \tfrac{m+1}{m+1+2m\alpha}\big[|u|^{m+1+2m\alpha}- k^{m+1+2m\alpha}\big].
\end{align*}
Finally, if $|u|\ge\ell$ we have
\begin{align*}
     v
    &=
    \int_{0}^{\ell^{m+1}} \phi_{\alpha,k,\ell }\big(s^{\frac{2m}{m+1}}\big)\,\ds +
    \int_{\ell^{m+1}}^{|u|^{m+1}} \phi_{\alpha,k,\ell }\big(s^{\frac{2m}{m+1}}\big)\,\ds\\
    &=
     k^{m+1+2m\alpha}+ \tfrac{m+1}{m+1+2m\alpha}\big[\ell^{m+1+2m\alpha}- k^{m+1+2m\alpha}\big]\\
     &\phantom{=\,}
     +
     \ell^{2m\alpha}\big[|u|^{m+1}-\ell^{m+1}\big].
\end{align*}
Joining the different cases it is straightforward to check that
\begin{align*}
    v&\ge
    \tfrac{m+1}{m+1+2 m\alpha}\phi_{\alpha ,k,\ell}
    \big(|\power{u}{m}|^2\big)^{\frac{m+1+2 m \alpha}{2m\alpha}}
    -k^{m+1+2m\alpha}\\
    &\ge 
     \tfrac{m+1}{m+1+2 m\alpha}\phi_{\alpha ,k,\ell}
    \big(|\power{u}{m}|^2\big)^{\frac{m+1+2 m \alpha}{2m\alpha}}-1,
\end{align*}
and
\begin{align*}
    v&\le 
    |\power{u}{m}|^{\frac{m+1+2 m \alpha}{m}}+k^{m+1+2m\alpha}\\
    &\le 
    |\power{u}{m}|^{\frac{m+1+2 m \alpha}{m}}+1.
\end{align*}
By differentiation, i.e.~the chain rule for Sobolev functions, and Kato's inequality we obtain
\begin{align*}
    &\Big|\nabla \big[\phi_{\alpha,k,\ell}\big(|\power{u}{m}|^2\big)\big]^\frac{1+\alpha}{2\alpha}\Big|^2 \\
    &\qquad=
    \Big|\tfrac{1+\alpha}{2\alpha } \big[\phi_{\alpha,\ell ,k}
    \big(|\power{u}{m}|^2\big)\big]^{\frac{1+\alpha}{2\alpha }-1} \phi_{\alpha,k,\ell}'\big(|\power{u}{m}|^2\big)\nabla |\power{u}{m}|^2\Big|^2\\
    &\qquad=
    \big(\tfrac{1+\alpha}{\alpha } \big)^2
    \phi_{\alpha,k,\ell }\big(|\power{u}{m}|^2\big) 
    \Big[\big[\phi_{\alpha,k,\ell}\big(|\power{u}{m}|^2\big)\big]^{\frac{1}{2\alpha}-1} 
    \phi_{\alpha,k,\ell}'\big(|\power{u}{m}|^2\big)\Big]^2 
    |\power{u}{m}|^2\big|\nabla |\power{u}{m}|\big|^2\\
    &\qquad\le
     \big(\tfrac{1+\alpha}{\alpha } \big)^2
     \phi_{\alpha,k,\ell }
     \big(|\power{u}{m}|^2\big) |D\power{u}{m}|^2 
  \Big[\big[\phi_{\alpha,k,\ell}\big(|\power{u}{m}|^2\big)\big]^{\frac{1}{2\alpha}-1} 
    \phi_{\alpha,k,\ell }'\big(|\power{u}{m}|^2\big)|\power{u}{m}|\Big]^2 .
\end{align*}
If $k^{2m}<|\power{u}{m}|^2<\ell^{2m}$ we have 
\begin{align*}
     \big[\phi_{\alpha,k,\ell}\big(|\power{u}{m}|^2\big)\big]^{\frac{1}{2\alpha}-1} 
     \phi_{\alpha,k,\ell}'\big(|\power{u}{m}|^2\big)|\power{u}{m}| 
     &=
    \alpha
     |\power{u}{m}|^{2\alpha(\frac{1}{2\alpha}-1)} 
     |\power{u}{m}|^{2(\alpha -1)}|\power{u}{m}|
     \le
     \alpha,
\end{align*}
while for $|\power{u}{m}|^2> \ell^{2m}$, respectively $|\power{u}{m}|^2<k^{2m}$ we have 
$\phi_{\alpha,k,\ell }'(|\power{u}{m}|^2)=0$. 
Hence, in any case we conclude that 
\begin{align*}
    \Big|\nabla \big[\phi_{\alpha,k,\ell}\big(|\power{u}{m}|^2\big)\big]^\frac{1+\alpha}{2\alpha }\Big|^2
     &\le
     (1+\alpha)^2
    |D\power{u}{m}|^2\phi_{\alpha,k,\ell }\big(|\power{u}{m}|^2\big).
\end{align*}
Finally, we have
\begin{align*}
    \phi_{\alpha,k,\ell}\big(|\power{u}{m}|^2\big)
    +
    |\power{u}{m}|^2 \phi_{\alpha,k,\ell}^{\prime}\big(|\power{u}{m}|^2\big)
    &\le
    (1+\alpha) \phi_{\alpha,k,\ell}\big(|\power{u}{m}|^2\big).
\end{align*}
If we insert these preparatory inequalities into \eqref{En-Est-Phi}, we obtain
\begin{align*}
    \underbrace{\tfrac{m+1}{m+1+2 m\alpha}}_{\ge\frac{m}{1+\alpha}}
    \sup_{\tau\in [t_o-R^2,t_o]}&\int_{B_R\times\{\tau\}} \phi_{\alpha ,k,\ell}
    \big(|\power{u}{m}|^2\big)^{\frac{m+1+2 m \alpha}{2m\alpha}}\,\zeta^2 \,\dx\\
    &\phantom{\le\,}
    + \tfrac{\nu }{(1+\alpha)^2}
    \iint_{Q_{R}}\Big|\nabla \phi_{\alpha,k,\ell}\big(|\power{u}{m}|^2\big)^\frac{1+\alpha}{2\alpha }\Big|^2 \zeta^2 \,\dx\dt \\
    &\le 
    C  \iint_{Q_{R}}|\power{u}{m}|^{2 } 
    \phi_{\alpha,k,\ell}(|\power{u}{m}|^{2})|\nabla \zeta|^2\,\dx\dt\\
    &\phantom{\le\,}+
    2\iint_{Q_{R}}
    \big[|\power{u}{m}|^{\frac{m+1+2 m \alpha}{m}}+1\big] |\partial_{t}\zeta^2|\,\dx\dt\\
    &\phantom{\le\,}
    +C (1+\alpha)^2 
    \iint_{Q_{R}} |F|^2  \phi_{\alpha,k,\ell}\big(|\power{u}{m}|^2\big)\zeta^2 \,\dx\dt\\
    &\phantom{\le\,} +
     \sup_{\tau\in [t_o-R^2,t_o]}\int_{B_R\times\{\tau\}} \zeta^2 \,\dx
\end{align*}
where $C=C(\nu,L)$.
Using the assumption $\zeta^2(\cdot,t_o-R^2)=0$, the last term on the right-hand side can be absorbed into the second integral, since
$$
    |\zeta^2(\cdot,\tau)|
    =
    \big|\zeta^2(\cdot,\tau)-\zeta^2(\cdot,t_o-R^2)\big|
    \le
    \int_{t_o-R^2}^\tau |\partial_t \zeta^2(\cdot, t)|\,\dt
    \le
    \|\partial_t \zeta^2\|_{L^\infty}R^2,
$$
immediately implies that
\begin{align*}
    \sup_{\tau\in [t_o-R^2,t_o]}\int_{B_R\times\{\tau\}} \zeta^2\,\dx
    &\le
    \|\partial_t \zeta^2\|_{L^\infty} |B_R|R^2=
    \|\partial_t \zeta^2\|_{L^\infty}
    \iint_{Q_{R}}1\, \dx\dt.
\end{align*}
Since
\begin{align*}
    \Big|\nabla \big[\phi_{\alpha,k,\ell} \big(|\power{u}{m}|^2\big)^\frac{1+\alpha}{2\alpha }\zeta\big]\Big|^2 
    \le
    2 \Big|\nabla \phi_{\alpha,k,\ell}\big(|\power{u}{m}|^2\big)^\frac{1+\alpha}{2\alpha }\Big|^2 \zeta^2
    +
    2 \phi_{\alpha,k,\ell}\big(|\power{u}{m}|^2\big)^\frac{1+\alpha}{\alpha }|\nabla\zeta |^2,
\end{align*}
the energy estimate becomes
\begin{align*}
    \tfrac{m}{1+\alpha}
    \sup_{\tau\in [t_o-R^2,t_o]}&\int_{B_R\times\{\tau\}} \phi_{\alpha ,k,\ell}
    \big(|\power{u}{m}|^2\big)^{\frac{m+1+2 m \alpha}{2m\alpha}}\zeta^2 \,\dx\\
    &\phantom{\le\,}
    + \tfrac{\nu }{(1+\alpha)^2}
    \iint_{Q_{R}}\Big|\nabla \big[\phi_{\alpha,k,\ell}\big(|\power{u}{m}|^2\big)^\frac{1+\alpha}{2\alpha }\zeta\big]\Big|^2  \,\dx\dt \\
    &\le 
    C  \iint_{Q_{R}}\phi_{\alpha,k,\ell}\big(
    |\power{u}{m}|^{2}\big)^\frac{1+\alpha}{\alpha }|\nabla \zeta|^2\,\dx\dt\\
    &\phantom{\le\,}+
    C\, \|\partial_{t}\zeta^2\|_{L^\infty}\iint_{Q_{R}}\big[ |\power{u}{m}|^{\frac{m+1+2 m \alpha}{m}}
    +1\big]\,\dx\dt\\
    &\phantom{\le\,}
    +C (1+\alpha)^2
    \iint_{Q_{R}} |F|^2  \phi_{\alpha,k,\ell}\big(|\power{u}{m}|^2\big)\zeta^2\,\dx\dt.
\end{align*}
As $\frac{m+1+2 m \alpha}{m} \geq 2+2 \alpha$ and $k\le 1$, we have 
%\color{red}
\[
    \phi_{\alpha,k,\ell}\big(|\power{u}{m}|^{2}\big)^\frac{1+\alpha}{\alpha }
     \le
     |\power{u}{m}|^{\frac{m+1+2 m \alpha}{m}}+1.
\]
%\color{black}
Therefore, in averaged form we get
\begin{align}\label{Eq:energy-est-av}\nonumber
    \boldsymbol E
    &:=
    \sup_{\tau\in [t_o-R^2,t_o]}\frac1{R^2}\bint_{B_R\times\{\tau\}} \phi_{\alpha ,k,\ell}
    \big(|\power{u}{m}|^2\big)^{\frac{m+1+2 m \alpha}{2m\alpha}}\zeta^2 \,\dx\\\nonumber
     &\qquad\qquad\qquad\qquad\quad + 
    \biint_{Q_{R}}\Big|\nabla \big[\phi_{\alpha,k,\ell}\big(|\power{u}{m}|^2\big)^\frac{1+\alpha}{2\alpha }\zeta\big]\Big|^2  \,\dx\dt \\\nonumber
    &\,\le 
    C (1+\alpha)^2 
    \left(\|\nabla \zeta\|_{L^\infty}^2+\|\partial_t \zeta^2\|_{L^{\infty}}\right) \biint_{Q_{R}}\big[|\power{u}{m}|^{\frac{m+1+2 m \alpha}{m}}+1\big]\,\dx\dt\\
    &\phantom{\le\,}
    +C (1+\alpha)^4
    \biint_{Q_{R}} |F|^2  \phi_{\alpha,k,\ell}\big(|\power{u}{m}|^2\big)\zeta^2\,\dx\dt.
\end{align}
By H\"older's inequality we get
\begin{align*}
    \biint_{Q_{R}}& |F|^2
    \phi_{\alpha,k,\ell}\big(|\power{u}{m}|^2\big)\zeta^2
    \,\dx\dt\\
    &\le
    \bigg[\biint_{Q_{R}} |F|^{2p}\,\dx\dt\bigg]^{\frac1p}
    \bigg[\biint_{Q_{R}}
    \big( \phi_{\alpha,k,\ell}\big(|\power{u}{m}|^2\big)\zeta^2\big)^{\frac{p}{p-1}}\,\dx\dt\Big]^{\frac{p-1}p}.
\end{align*}
Since $p>\frac{N+2}2$, we have $\frac{p}{p-1}<\frac{N+2}{N}$. Therefore, by interpolation
for $1<\frac p{p-1}<\frac{N+2}N$ we obtain for $\varep>0$ that
\begin{align}\label{after-inter}
    \bigg[\biint_{Q_{R}}
    \big( \phi_{\alpha,k,\ell}\big(|\power{u}{m}|^2\big)\zeta^2\big)^{\frac{p}{p-1}}\,\dx\dt
    \bigg]^{\frac{p-1}p} 
    &
    \le 
    \varepsilon
    \mathbf I^{\frac{N}{N+2}}
    +\varepsilon^{-\frac{N+2}{2p-(N+2)}}\mathbf{II},
\end{align}
where
\begin{equation*}
    \mathbf I
    :=
    \biint_{Q_{R}}\big( \phi_{\alpha,k,\ell}\big(|\power{u}{m}|^2\big)\zeta^2\big)^{\frac{N+2}{N}}\,\dx\dt,
\end{equation*}
and
\begin{equation*}
    \mathbf{II}:=  \biint_{Q_{R}} \phi_{\alpha,k,\ell}\big(|\power{u}{m}|^2\big)\zeta^2\,\dx\dt.
\end{equation*}
Let 
\begin{equation}\label{def:q-first}
    q
    :=
    2 \frac{N+\frac{m+1+2 m \alpha}{m(1+\alpha)}}{N}
    =
    2\Big(1+\tfrac{m+1+2 m \alpha}{Nm(1+\alpha)}\Big).
\end{equation}
Then, it is a matter of straightforward computations to verify that
$$
    2\alpha \frac{N+2}{N}
    \le
    2 \frac{N(1+\alpha)+\frac{m+1+2 m \alpha}{m}}{N}
    =
    (1+\alpha) q. 
$$
Therefore, taking into account that $0\le\zeta\le1$, we have
\begin{align*}
    \textbf{I}
    &=
    \biint_{Q_{R}}\big( \phi_{\alpha,k,\ell}\big(|\power{u}{m}|^2\big)^\frac{1}{2\alpha}\zeta^\frac1{\alpha}\big)^{2\alpha\frac{N+2}{N}}\,\dx\dt \\
    &\le
    \biint_{Q_{R}}\big[\big(\phi_{\alpha,k,\ell}\big(|\power{u}{m}|^2\big)^\frac{1}{2\alpha}\zeta^{\frac1\alpha}\big)^{(1+\alpha)q}+1\big]\,\dx\dt \\
    &\le
    \biint_{Q_{R}}
    \big(\phi_{\alpha,k,\ell}\big(|\power{u}{m}|^2\big)^\frac{1+\alpha}{2\alpha}\zeta\big)^q
    \,\dx\dt+1.
\end{align*}
Next, we apply Gagliardo-Nirenberg's inequality, cf. \cite[Chapter~I, Proposition~3.1]{DiBe}, with $(v,p,m,q)$ replaced by $\big( 
\phi_{\alpha,k,\ell}(|\power{u}{m}|^2)^\frac{1+\alpha}{2\alpha}\zeta,\,2,\, \frac{m+1+2 m\alpha}{m(1+\alpha)},q\big)$ and obtain
\begin{align*}
    \textbf{I}
    &\le
    C\biint_{Q_{R}}
    \big|\nabla\big[
    \phi_{\alpha,k,\ell}\big(|\power{u}{m}|^2\big)^\frac{1+\alpha}{2\alpha}\zeta\big]\big|^2\,\dx\dt\\
    &\qquad\qquad 
    \cdot\bigg[\sup_{\tau\in[ t_o-R^2,t_o]} 
    \int_{B_R\times\{\tau\}} 
    \big[
    \phi_{\alpha,k,\ell}\big(|\power{u}{m}|^2\big)^\frac{1+\alpha}{2\alpha}
    \zeta\big]^{\frac{m+1+2 m\alpha}{m(1+\alpha)}}\,\dx\bigg]^{\frac2N} + 1\\
    &\le 
    C\,R^{\frac{2(N+2)}{N}} \boldsymbol{E}^{\frac{N+2}{N}} + 1.
\end{align*}
If we substitute it back into \eqref{after-inter}, we obtain
\begin{align*}
    \bigg[\biint_{Q_{R}}
    \big( \phi_{\alpha,k,\ell}\big(|\power{u}{m}|^2\big)\zeta^2\big)^{\frac{p}{p-1}}\,\dx\dt\bigg]^{\frac{p-1}p}
    \le 
    C\,\varepsilon R^2 \boldsymbol{E}+\varepsilon+\varepsilon^{-\frac{N+2}{2p-(N+2)}} \mathbf{II}.
\end{align*}
This allows us to estimate the last term of \eqref{Eq:energy-est-av} by
\begin{align*}
    C & (1+\alpha)^{4}
    \biint_{Q_{R}}|F|^2  \phi_{\alpha,k,\ell}\big(|\power{u}{m}|^2\big)\zeta^2 \,\dx\dt\\
    &\le C\, \varepsilon
    (1+\alpha)^4
    [F]_{2p,Q_R}^2 
    \boldsymbol E +
    \frac{C}{R^2} (1+\alpha)^4 [F]_{2p,Q_R}^2 \Big[\varepsilon+\varepsilon^{-\frac{N+2}{2p-(N+2)}} \mathbf{II}\Big].
\end{align*}
If we insert it in \eqref{Eq:energy-est-av} and shift the term containing $\boldsymbol E$ to the left-hand side, we obtain
\begin{align*}
    &\Big[1-C\,\varepsilon (1+\alpha)^4 [F]_{2p,Q_R}^2 \Big]\boldsymbol{E}\\
    &\qquad\le 
    C\,(1+\alpha)^2\left(\|\nabla \zeta\|_{L^\infty}^2+\|\partial_t \zeta^2\|_{L^\infty}\right) 
    \biint_{Q_{R}}
    \big[|\power{u}{m}|^{\frac{m+1+2 m \alpha}{m}}+1\big]\,\dx\dt\\
    &\qquad\phantom{\le\,}
    + \frac{C}{R^2}\,\varepsilon (1+\alpha)^4[F]_{2p,Q_R}^2 
    + \frac{C}{R^2\varepsilon^{\frac{N+2}{2p-(N+2)}}}\, (1+\alpha)^4[F]_{2p,Q_R}^2 \mathbf{II}.
\end{align*}
In the above estimate, choose $\varepsilon$ such that 
\begin{equation*}
    C\,\varepsilon (1+\alpha)^4
    [F]_{2p,Q_R}^2
    =
    \tfrac12.
\end{equation*}
Then, the left-hand side becomes $\frac12\boldsymbol{E}$, the second term on the right-hand side is less than or equal to $\frac12R^{-2}$, and the third term is
\begin{align*}
    &\frac{C}{R^2\varepsilon^{\frac{N+2}{2p-(N+2)}}}\, 
    (1+\alpha)^4 [F]_{2p,Q_R}^2 \mathbf{II} \\
    &\qquad =
    \frac{C}{R^2}\,(1+\alpha)^4 \Big[2C\,(1+\alpha)^4 [F]_{2p,Q_R}^2 \Big]^{\frac{N+2}{2p-(N+2)}}
    [F]_{2p,Q_R}^2
    \mathbf{II}\\
    &\qquad\le 
    \frac{C}{R^2}
    (1+\alpha)^{\frac{8p}{2p-(N+2)}}
    [F]_{2p,Q_R}^{\frac{4p}{2p-(N+2)}}\mathbf{II} .
\end{align*}
Collecting all the terms and recalling the definition of $\boldsymbol{E}$ and $\mathbf{II}$, we get
\begin{align*}
    \sup_{\tau\in [t_o-R^2,t_o]}
    &
    \bint_{B_R\times\{\tau\}} \phi_{\alpha ,k,\ell}
    \big(|\power{u}{m}|^2\big)^{\frac{m+1+2 m \alpha}{2m\alpha}}\zeta^2 \,\dx\\
    &\phantom{\le\,}
    + 
    R^2
    \biint_{Q_{R}}\Big|\nabla \big[\phi_{\alpha,k,\ell}\big(|\power{u}{m}|^2\big)^\frac{1+\alpha}{2\alpha }\zeta\big]\Big|^2  \,\dx\dt\\ 
    &
    \le C\,R^2
    (1+\alpha)^2\big(\|\nabla \zeta\|_{L^\infty}^2+\|\partial_t \zeta^2\|_{L^\infty}\big) \biint_{Q_{R}}
    \big[|\power{u}{m}|^{\frac{m+1+2 m \alpha}{m}}+1\big]\,\dx\dt\\
    &\phantom{\le\,}
    +C
    + 
    C\,(1+\alpha)^{\frac{8p}{2p-(N+2)}}
    [F]_{2p,Q_R}^\frac{4p}{2p-(N+2)}
    \biint_{Q_{R}} \phi_{\alpha,k,\ell}
    \big(|\power{u}{m}|^2\big)\zeta^2
    \,\dx\dt,
\end{align*}
that is,
\begin{align}
     \sup_{\tau\in [t_o-R^2,t_o]}
    &
    \bint_{B_R\times\{\tau\}} \phi_{\alpha ,k,\ell}
    \big(|\power{u}{m}|^2\big)^{\frac{m+1+2 m \alpha}{2m\alpha}}\zeta^2 \,\dx\nonumber\\
    &\phantom{\le\,}
    + 
    R^2
    \biint_{Q_{R}}\Big|\nabla \big[\phi_{\alpha,k,\ell}\big(|\power{u}{m}|^2\big)^\frac{1+\alpha}{2\alpha }\zeta\big]\Big|^2  \,\dx\dt\label{En-Est-useful}\\
    &\le
    C \,(1+\alpha)^{\frac{8p}{2p-(N+2)}}
    \boldsymbol\mu(R)
     \biint_{Q_{R}}
    \big[|\power{u}{m}|^{\frac{m+1+2 m \alpha}{m}}+1\big]\,\dx\dt,\nonumber
\end{align}
where we let
\begin{align}\label{def:mu}
    \boldsymbol\mu(R)
    &:=
    1+ R^2 \big(\|\nabla \zeta\|_{L^\infty}^2
    +\|\partial_t \zeta^2\|_{L^\infty}\big)
    +
    [F]_{2p,Q_R}^\frac{4p}{2p-(N+2)}.
\end{align}
To obtain the last line we used the fact that
$$
    \phi_{\alpha,k,\ell}\big(|\power{u}{m}|^2\big)
    \le
    |\power{u}{m}|^{\frac{m+1+2 m \alpha}{m}}+1.
$$
By the Gagliardo Nirenberg's inequality \cite[Chapter~I, Proposition~3.1]{DiBe} applied with $(v,p,m,q)$ replaced by $\big(\phi_{\alpha,k,\ell}\big(|\power{u}{m}|^2\zeta\big)^\frac{1+\alpha}{2\alpha },2, \frac{m+1+2 m\alpha}{m(1+\alpha)},q\big)$, with $q$ defined in \eqref{def:q-first}, we conclude 
\begin{align*}
    \biint_{Q_{R}} &
    \phi_{\alpha,k,\ell}\big(|\power{u}{m}|^2\big)^{q\frac{1+\alpha}{2\alpha}}
    \zeta^q\,\dx\dt \\
    &\le C\,
    R^2\biint_{Q_{R}}\big|\nabla\big[
    \phi_{\alpha,k,\ell}\big(|\power{u}{m}|^2\big)^\frac{1+\alpha}{2\alpha }\zeta\big]\big|^2\,\dx\dt\\
    &\phantom{\le\,}\cdot
    \Bigg[\sup _{\tau\in[ t_o-R^2,t_o]} 
    \bint_{B_R\times\{\tau\}} 
    \big(\phi_{\alpha,k,\ell}\big(|\power{u}{m}|^2\big)^\frac{1+\alpha}{2\alpha }
    \zeta\big)^\frac{m+1+2 m\alpha}{m(1+\alpha)}
    \,\dx\Bigg]^{\frac{2}{N}}\\
    &\le 
    C\,
    R^2\biint_{Q_{R}}\big|\nabla\big[
    \phi_{\alpha,k,\ell}\big(|\power{u}{m}|^2\big)^\frac{1+\alpha}{2\alpha }\zeta\big]\big|^2\,\dx\dt\\
     &\phantom{\le\,}\cdot
    \Bigg[\sup _{\tau\in[ t_o-R^2,t_o]} 
    \bint_{B_R\times\{\tau\}} 
   \phi_{\alpha,k,\ell}\big(|\power{u}{m}|^2\big)^\frac{m+1+2 m\alpha}{2m\alpha}
    \zeta^2\,\dx\Bigg]^{\frac{2}{N}}.
\end{align*}
To obtain the last line we used the fact that $\frac{m+1+2m\alpha}{m(1+\alpha)}\ge 2$ and $0\le\zeta\le 1$. Both factors on the right-hand side can be estimated in terms of the left-hand side of \eqref{En-Est-useful}. This yields
\begin{align*}
    \bigg[ 
    \biint_{Q_{R}} &
    \phi_{\alpha,k,\ell}\big(|\power{u}{m}|^2\big)^{q\frac{1+\alpha}{2\alpha}}
    \zeta^q\,\dx\dt
    \bigg]^{\frac N{N+2}}\\
    &\le C\,(1+\alpha)^{\frac{8p}{2p-(N+2)}}
    \boldsymbol\mu(R)
     \biint_{Q_{R}}\big[|\power{u}{m}|^{\frac{m+1+2 m \alpha}{m}}+1\big]\,\dx\dt.
\end{align*}
Note that the right-hand side is finite and independent of $k$ and $\ell$. Moreover, we have
\begin{equation*}
    \lim_{\ell\to \infty}\lim_{k\downarrow 0}
    \phi_{\alpha,k,\ell}\big(|\power{u}{m}|^2\big)^{q\frac{1+\alpha}{2\alpha}}= |\power{u}{m}|^{(1+\alpha)q}.
\end{equation*}
Therefore, by Fatou's  lemma we can pass first to the limit $k\downarrow 0$ and then by the theorem of monotone convergence to
the limit $\ell\to\infty$ to obtain the {\it reverse  H\"older inequality}
\begin{align}\label{Moser-energy:0}\nonumber
    \bigg[ 
    \biint_{Q_{R}} &\zeta^q
    |\power{u}{m}|^{(1+\alpha)q}
    \,\dx\dt
    \bigg]^{\frac N{N+2}}\\
    &\le C\, (1+\alpha)^{\frac{8p}{2p-(N+2)}}
    \boldsymbol\mu(R)
     \biint_{Q_{R}}\big[|\power{u}{m}|^{\frac{m+1+2 m \alpha}{m}}+1\big]\,\dx\dt,
\end{align}
for a constant $C=C(N,\nu,L)$, and where $\boldsymbol \mu(R)$ is defined in \eqref{def:mu}. Note that \eqref{Moser-energy:0} is valid
for any cut-off function $\zeta\in W^{1,\infty}(Q_R,[0,1])$ that vanishes on 
the parabolic boundary $\partial_P Q_R$. Specifying $\zeta$ in such a way that $\zeta=1$ on $Q_{\frac12 R}$, $|\nabla\zeta|\le\frac2{R}$ and $|\partial_t\zeta|\le \frac{4}{R^2}$, we get
\begin{align*}%\label{reverse-Hoelder}\nonumber
    \bigg[ 
    \biint_{Q_{\frac12 R}} &
    |\power{u}{m}|^{(1+\alpha)q}
    \,\dx\dt
    \bigg]^{\frac N{N+2}}\\
    &\le C\, (1+\alpha)^{\frac{8p}{2p-(N+2)}}
     \Big( 1+
    [F]_{2p,Q_R}^\frac{4p}{2p-(N+2)}
   \Big)
     \biint_{Q_{R}}\big[|\power{u}{m}|^{\frac{m+1+2 m \alpha}{m}}+1\big]\,\dx\dt.
\end{align*}
By a covering argument we easily obtain the conclusion
\begin{equation*}
     |\power{u}{m}|\in L^\frac{m+1+2m\alpha}{m}_{\rm loc} (\Omega_T)
     \quad\Longrightarrow\quad
     |\power{u}{m}|\in L^{(1+\alpha)q}_{\rm loc} (\Omega_T).
\end{equation*}
The reverse H\"older inequality \eqref{Moser-energy:0} is the starting point for the Moser iteration, which eventually leads to the boundedness estimates. 

\textit{Step 2.} To set up the Moser iteration scheme, we let $Q_o:=Q_{\rho}(z_o)\subset\Omega_T$ and define a shrinking family of cylinders  $Q_i:=
Q_{\rho_i}(z_o)=B_{\rho_i}(x_o)\times(t_o-\rho_i^2,t_o]$ for $i\in\N_0$ with the same base point $z_o=(x_o,t_o)$, where
$\rho_i$  is defined according to
$$
    \rho_i=\frac{\rho}2 +\frac{\rho}{2^{i+1}},\quad \mbox{for $i\in\N_0$.}
$$ 
Observe that $\rho_i$ is decreasing, $\rho_o=\rho$ and $\rho_\infty=\frac12 \rho$.
Moreover, choose the cutoff function $\zeta\in W^{1,\infty}(Q_i,[0,1])$ such that $\zeta=1$ in 
$Q_{i+1}$, $\zeta =0$ on the parabolic boundary of $Q_i$, and 
\begin{equation*}
    |\nabla\zeta|\le\frac{2^{i+2}}{\rho},
    \quad\mbox{and}\quad
    |\partial_t\zeta|\le\frac{2^{2(i+2)}}{\rho^2}.
\end{equation*}
In this setting, we obtain from \eqref{Moser-energy:0} that
\begin{align*}
    & \bigg[\frac{|Q_{i+1}|}{|Q_i|} \biint_{Q_{i+1}} |\power{u}{m}|^{(1+\alpha) q}\,\dx\dt\bigg]^{\frac N{N+2}}\\
    &\ \ \le 
    C\, (1+\alpha)^{\frac{8p}{2p-(N+2)}}
    \bigg[1+\rho_i^2 \frac{2^{2(i+2)}}{\rho^2} +[F]_{2p,Q_{i}}^{\frac{4p}{2p-(N+2)}}\bigg]
    \biint_{Q_{i}}\big[|\power{u}{m}|^{\frac{m+1+2 m \alpha}{m}}+1\big]\,\dx\dt\\
    &\ \ \le
    C\,4^i (1+\alpha)^{\frac{8p}{2p-(N+2)}}
    \Big(1+[F]_{2p,Q_{o}}^{\frac{4p}{2p-(N+2)}}\Big)\biint_{Q_{i}}
    \big[|\power{u}{m}|^{\frac{m+1+2 m \alpha}{m}}+1\big]\,\dx\dt.
\end{align*}
Notice that the quantity
$[F]_{2p,Q_{o}}$ remains unchanged during the iteration. 
Observing that $\frac{\left|Q_{i+1}\right|}{\left|Q_i\right|} \geq \frac{1}{2^{N+2}}$, the last displayed inequality yields 
\begin{align}\label{Iter-start} \nonumber
    \bigg[
    \biint_{Q_{i+1}} & |\power{u}{m}|^{(1+\alpha) q}\,\dx\dt
    \bigg]^{\frac N{N+2}}\\
    &\le 
    C\, 4^i (1+\alpha)^{\frac{8p}{2p-(N+2)}}\Big(1+[F]_{2p,Q_{o}}^{\frac{4p}{2p-(N+2)}}\Big)
    \biint_{Q_{i}}\big[|\power{u}{m}|^{\frac{m+1+2 m \alpha}{m}}+1\big]\,\dx\dt.
\end{align}
Letting $\kappa:=1+\frac2N$, and considering 
\begin{equation*}
    (1+\alpha)q
    =
    2\alpha\kappa +2 \frac{N m+m+1}{N m},
    \quad
    \mbox{and}\quad
    \frac{m+1+2 \alpha m}{m}=2 \alpha+\frac{m+1}{m},
\end{equation*}
it is reasonable to define $\alpha_o$ by
\begin{equation*}
    \alpha_o=\frac{r-(m+1)}{2 m},
\end{equation*}
and  the sequence $(\alpha_i)_{i\in\N_0}$  recursively by
\begin{equation*}
    2 \alpha_{i+1}:=2 \alpha_i \kappa+2 \frac{N m+m+1}{N m}-\frac{m+1}{m}=2 \alpha_i \kappa+\frac{N(m-1)+2(m+1)}{N m}.
\end{equation*}
Then,
\begin{equation*}
    m \frac{m+1+2  m\alpha_o}{m}=r,
\end{equation*}
and therefore, the iteration process starts with the $L^r$-integral of $|u|$ on $Q_o$ on the right-hand side, which is finite by assumption.
Recalling the definition of $\boldsymbol \lambda_r$ in \eqref{Eq:B:4:1}, by simple calculations, we obtain that
$$
    \alpha_i=\frac{\boldsymbol\lambda_r}{4 m} \kappa^i-\frac{N(m-1)+2(m+1)}{4 m},
    \qquad\forall\, i\in \N_0.
$$
If we take the definition of the recursive sequence $\{\alpha_i\}$ into account, \eqref{Iter-start} leads to 
\begin{align*}
    \bigg[\biint_{Q_{i+1}}
    &
    |\power{u}{m}|^{2 \alpha_{i+1}+\frac{m+1}{m}}\,\dx\dt\bigg]^{\frac N{N+2}}\\
    &\le 
    C\,
    (1+\alpha_i)^{\frac{8p}{2p-(N+2)}}4^i
    \Big(1+[F]_{2p,Q_{o}}^{\frac{4p}{2p-(N+2)}}\Big)
    \biint_{Q_{i}}
    \big[|\power{u}{m}|^{2\al_i + \frac{m+1}{m}}+1\big]\,\dx\dt.
\end{align*}
To proceed further, we need to control $1+\alpha_i$ in terms of $\kappa^i$. Using in turn the definition of $\alpha_i$, and the fact that
$N(1-m)-2(m+1)\ge 0$, which holds since we are dealing with the sub-critical case $0<m\le\frac{N-2}{N+2}$,  we have
\begin{align*}
    1+\alpha_i 
    &=
    1+\frac{\boldsymbol\lambda_r}{4 m} \kappa^i-\frac{N(m-1)+2(m+1)}{4 m}\\
    & =
    1+\frac{N(m-1)+2 r}{4 m} \kappa^i+\frac{N(1-m)-2(m+1)}{4 m}\\
    &= 
    \frac{r}{2 m}\kappa^i -
    \frac{2(1-m)+N(1-m)(\kappa^i-1)}{4m}
    \le \frac{r}{2 m}\kappa^i.
\end{align*}
Inserting this in the second last displayed inequality and adding 1 to the left-hand side, we obtain the {\it recursive inequality}
\begin{align}\label{rec-Hoelder}\nonumber
    \bigg[\biint_{Q_{i+1}}&
    \big[|\power{u}{m}|^{2 \alpha_{i+1}+\frac{m+1}{m}}+1\big]\,\dx\dt\bigg]^{\frac N{N+2}}\\\nonumber
    &\le 
    C\, \Big(\frac{r}{2 m}\kappa^{i}\Big)^{\frac{8p}{2p-(N+2)}}
    \,4^i\Big(1+[F]_{2p,Q_{o}}^{\frac{4p}{2p-(N+2)}}\Big)
    \biint_{Q_{i}}\big[|\power{u}{m}|^{2\alpha_i+\frac{m+1}{m}}+1\big]\,\dx\dt \\\nonumber
    &\le 
    C\, \Big(\frac{r}{2 m}\Big)^{\frac{8p}{2p-(N+2)}}
    \Big(4\kappa^\frac{8p}{2p-(N+2)}\Big)^i
    \Big(1+[F]_{2p,Q_{o}}^{\frac{4p}{2p-(N+2)}}\Big)\\
    &\quad\qquad\qquad\qquad\qquad\qquad\cdot
    \biint_{Q_{i}}\big[|\power{u}{m}|^{2\alpha_i+\frac{m+1}{m}}+1\big]\,\dx\dt. 
\end{align}
The series of inequalities \eqref{rec-Hoelder}$_i$ must be understood as follows. If the right-hand side integral in \eqref{rec-Hoelder}$_i$ is finite, then the left-hand side integral is also finite. This ensures that the right-hand side in \eqref{rec-Hoelder}$_{i+1}$ is finite. For $i=0$ the right-hand side is finite, because in this case $2\alpha_o+\frac{m+1}{m}=\frac{r}{m}$. Let
\begin{align*}
    & A=C\, \Big(\frac{r}{2 m}\Big)^{\frac{8p}{2p-(N+2)}}
    \Big(1+[F]_{2p,Q_{o}}^{\frac{4p}{2p-(N+2)}}\Big),\\
    &\boldsymbol b=4\kappa^{\frac{8p}{2p-(N+2)}}, 
    \quad 
    p_i=2 \alpha_i+\frac{m+1}{m} = \frac{\boldsymbol \lambda_r}{2 m} \kappa^i+\frac{N(1-m)}{2 m},\\
    & \boldsymbol{Y}_i=\bigg[\biint_{Q_i}\big[|\power{u}{m}|^{p_i}+1\big]\,\dx\dt\bigg]^{\frac1{p_i}},
\end{align*}
and obtain
$$
    \boldsymbol Y_{i+1}^{p_{i+1}} 
    \le
    \big(
    A \boldsymbol b^i \boldsymbol Y_i^{p_i}
    \big)^{\kappa}, \quad \forall \, i \in \mathbb{N}_0.
$$
The iteration yields
\begin{equation*}
    \boldsymbol Y_i
    \le
    \prod_{j=1}^i A^{\frac{\kappa^{i-j+1}}{p_i}}
    \prod_{j=1}^i \boldsymbol b^{j\frac{\kappa^{i-j+1}}{p_i}} 
    \boldsymbol Y_o^{\frac{p_o \kappa^i}{p_i}},\qquad\forall\,i\in\mathbb{N}_0.
\end{equation*}
Note that $p_o=\frac rm$. In order to compute the products,
we observe that
\begin{align*}
    \frac{\kappa^{i-j+1}}{p_i}
    =
    \frac{\kappa^{i-j+1}}{\frac{\boldsymbol \lambda_r}{2 m} \kappa^i+\frac{N(1-m)}{2 m}} 
    \le \frac{2 m}{\boldsymbol \lambda_r} \kappa^{-j+1},
\end{align*}
and
\begin{align*}
    \frac{p_o \kappa^i}{p_i}
    =
    \frac{2 \alpha_o+\frac{m+1}{m}}{\frac{\boldsymbol \lambda_r}{2 m} \kappa^i+\frac{N(1-m)}{2 m}} \kappa^i 
    \le 
    \frac{2 \alpha_o+\frac{m+1}{m}}{\frac{\boldsymbol \lambda_r}{2 m}} 
    =
    \frac{2 r}{\boldsymbol \lambda_r}.
\end{align*}
This allows to estimate
\begin{align*}
     \sum_{j=1}^i \frac{\kappa^{i-j+1}}{p_i} 
     &\le 
     \frac{2 m}{\boldsymbol\lambda_r} 
     \sum_{j=1}^i \kappa^{-j+1} 
     \le 
     \frac{2 m}{\boldsymbol\lambda_r} \frac{1}{1-1 /\kappa}
     =
     \frac{2 m}{\boldsymbol\lambda_r} \frac{N+2}{2},
\end{align*}
and similarly
\begin{align*}
    \sum_{j=1}^i j \frac{\kappa^{i-j+1}}{p_i}
    &\le 
    \frac{2 m}{\boldsymbol \lambda_r} \sum_{j=1}^i j \kappa^{-j+1} 
    \le 
    \frac{2 m}{\boldsymbol \lambda_r}\Big(\frac{\kappa}{\kappa-1}\Big)^2
    =
    \frac{2 m}{\boldsymbol \lambda_r}\Big(\frac{N+2}{2}\Big)^2.
\end{align*}
Hence, we obtain
\begin{equation*}
  \limsup_{i \rightarrow \infty}  \boldsymbol{Y}_i 
  \le 
  A^{\frac{2 m}{\boldsymbol\lambda_r} \frac{N+2}{2}} 
  \boldsymbol b^{\frac{2 m}{\boldsymbol\lambda_r}(\frac{N+2}{2})^2} 
    \boldsymbol{Y}_o^{\frac{2r}{\boldsymbol \lambda_r}}.  
\end{equation*}
Straightforward computations show on the one hand that
\begin{align*}
    \boldsymbol{Y}_o^{\frac{2r}{\boldsymbol\lambda_r}}
    &=
    \bigg[\biint_{Q_o}
    \big[|\power{u}{m}|^{2\alpha_o+\frac{m+1}m}+1\big]\,\dx\dt
    \bigg]^{\frac{2r}{p_o\boldsymbol \lambda_r}}
    =
    \bigg[
    \biint_{Q_o}\big[|u|^r+1\big]\,\dx\dt
    \bigg]^{\frac{2m}{\boldsymbol\lambda_r}},
\end{align*}
and on the other hand that
\begin{align*}
    A^{\frac{N+2}{2}}
    &
    =
    \bigg[C\, \Big(\frac r{2m}\Big)^{\frac{8p}{2p-(N+2)}}
    \Big( 
    1+[F]_{2p,Q_{o}}^{\frac{4p}{2p-(N+2)}}\Big)
    \bigg]^{\frac{N+2}2}
    \equiv
    C\Big(1+[F]_{2p,Q_{o}}^{\frac{2p(N+2)}{2p-(N+2)}}\Big).
\end{align*}
Hence, we conclude that
\begin{equation*}
    \sup_{Q_{\frac12\rho}}|u|
    \le
    C\bigg[
    \Big(
    1+[F]_{2p,Q_{o}}^{\frac{2p(N+2)}{2p-(N+2)}}
    \Big)
    \biint_{Q_{\rho}}\big[|u|^r+1\big]\,\dx\dt
    \bigg]^{\frac{2}{\boldsymbol\lambda_r}},
\end{equation*}
where $C$ depends on $N$, $m$, $\nu$, $L$, $p$, and $r$.
\end{proof}
%%%%%%%%%%%%%%%%%%%%%%

Proposition~\ref{prop:sup-est-qual} was proved for the prototype porous medium equation, and without right-hand side in \cite{bonforte}. 
We now turn the previous qualitative estimate into a quantitative one.
Such a quantitative result was first proved in \cite{DiBe-Kwong} for the prototype porous medium equation without right-hand side.

\begin{proposition}[quantitative $L^\infty$-estimate]\label{prop:sup-est}
Let $m\in(0,m_c]$, $p>\frac{N+2}2$, and $r\ge 1$ satisfying \eqref{Eq:B:4:1}. Then there exists a constant
$C=C(N,m,\nu,L , p,r)$ such that whenever $u$ is a locally bounded, local, weak solution to \eqref{por-med-eq} in $\Omega_T$ in the sense of Definition~{\upshape\ref{def:weak_solution}}  with $F\in L^{2p}_{\loc}(\Omega_T,\R^{kN})$,  then for all cylinders
$Q_{\rho,\vartheta}\equiv Q_{\rho,\vartheta}(z_o)\Subset \Omega_T$ and all $\sigma\in [\frac12 ,1)$ we have 
\begin{align*}
    \sup_{Q_{\sigma\rho,\sigma\vartheta}}  |u|
    &\le
    \max\Bigg\{
    \bigg[
     \frac{C}{(1-\sigma)^{N+2}}
    \Big(\frac{\varrho^2}{\vartheta}\Big)^\frac{N}{2}
    \biint_{Q_{\varrho,\vartheta}}|u|^{r}\,\dx\dt
    \bigg]^\frac{2}{\boldsymbol\lm_{r}},\\
    &\qquad\qquad\;
    \bigg[
    \Big(\frac{\vartheta}{\varrho^2}\Big)^{p-\frac{N}2}
    \biint_{Q_{\varrho,\vartheta}}|\varrho F|^{2p}\,\dx\dt\bigg]^\frac{2}{\boldsymbol \lambda_{p(1+m)}},
    \Big(
    \frac{\vartheta}{\rho^2}
    \Big)^{\frac{1}{1-m}}
    \Bigg\}.
\end{align*}
\end{proposition}

\begin{proof}
Let $Q_o:= Q_{\varrho,\vartheta}=B_\rho \times(-\vartheta, 0]\Subset\Omega_T$.
We define a family of nested, shrinking cylinders with common vertex by $Q_n:=Q_{\rho_n,\vartheta_n}$, where $\rho_n$, $\vartheta_n$ 
are defined by
\[
    \rho_n
    =\sigma \rho+\frac{1-\sigma}{2^n} \rho, 
    \quad 
    \vartheta_n=\sigma \vartheta+\frac{1-\sigma}{2^n} \vartheta,
\]
and $\sigma \in(0,1)$ is fixed.
Moreover, we introduce  increasing sequences of levels by %\textcolor{blue}{($\tilde k_n$ is not used.)}
\[
    k_n:=k-\frac{k}{2^n}, 
    %\quad 
    %\tilde{k}_n=\frac{1}{2}\left(k_n+k_{n+1}\right),
\]
where $k>0$ is a number to be determined in a universal way later.
Taking into account the relationship between two consecutive levels
\begin{equation*}
    k_n=k_{n+1} \frac{2^{n+1}-2}{2^{n+1}-1},
\end{equation*}
we have on the set $\{|u|>k_{n+1}\}$ that
\begin{align*}
    |u|^m-k_n^m
    &=
    |u|^m-k_{n+1}^m \Big(\frac{2^{n+1}-2}{2^{n+1}-1}\Big)^m
    =
    |u|^m\bigg[1-\Big(\frac{k_{n+1}}{|u|}\Big)^m \Big(\frac{2^{n+1}-2}{2^{n+1}-1}\Big)^m \bigg] \\
    &>
    |u|^m\bigg[1- \Big(\frac{2^{n+1}-2}{2^{n+1}-1}\Big)^m \bigg]
    \ge
    m |u|^m \bigg(1-\frac{2^{n+1}-2}{2^{n+1}-1}\bigg) 
    \ge 
    \frac{m}{2^{n+1}} |u|^m,
\end{align*}
where we used the inequality $1-x^m\ge m(1-x)$ for $x\in [0,1]$. Hence, rewriting this inequality we have 
\begin{align}\label{est:DG-pre}
    |u|^m
    \le 
    \frac{2^{n+1}}{m} \big(|u|^m-k_n^m\big)
    \qquad \mbox{on $\{|u|>k_{n+1}\}$.}
\end{align}

The energy inequality will be applied on the cylinders $Q_{n+1}\subset Q_n$, together with a cut-off function $\zeta_n$ such that $\z_n=0$ on  $\pl_P Q_n$; $\z_n=1$ on  $Q_{n+1}$; and moreover,
\[
    \left|\nabla \zeta_n\right| 
    \le 
    \frac{2^{n+1}}{(1-\sigma) \rho}, \quad 
    |\partial_t \zeta_n|
    \le \frac{2^{n+1}}{(1-\sigma)\vartheta} .
\]
With respect to these choices, the  energy estimates from Lemma \ref{lem:new-energy-est} yield
\begin{align}\label{en:start}\nonumber
    \sup_{-\vartheta_n \le \tau \le 0} 
    \int_{B_n \times\{\tau\}}&
    \big(|u|^m-k_{n+1}^m\big)_+^{3+\frac{1}{m}} \zeta_n^2\,\dx 
    +
    \iint_{Q_n}\big|\nabla\big[\big(|u|^m-k_{n+1}^m\big)_+^2 \zeta^2_n
    \big]\big|^2\,\dx\dt\\\nonumber
    &\quad
    \le
    C\bigg[ \frac{2^{2(n+1)}}{(1-\sigma)^2\rho^2}\iint_{Q_n}|u|^{4 m} \mathbf 1_{\{|u|>k_{n+1}\}}
    \,\dx\dt\\\nonumber
    &\quad\phantom{\le C\bigg[} +
    \frac{2^{n+1}}{(1-\sigma)\vartheta}\iint_{Q_n}|u|\big(|u|^{m}-k_{n+1}^m\big)_+^3\,
   \,\dx\dt\bigg] \\
    &\quad\phantom{\le C\bigg[} +
    \iint_{Q_n}|F|^2|u|^{2 m} \mathbf 1_{\{|u|>k_{n+1}\}} \,\dx\dt
    .
\end{align}
Since $0<m<1$, it is apparent that $4 m<1+ 3 m$. Moreover, since $0<m \le \frac{(N-2)_+}{N+2}$, the required condition  $\boldsymbol{\lambda}_r=N(m-1)+2 r>0$ implies that $r>\frac{2 N}{N+2}$. The crucial point is the relationship between $r$ and $1+3 m$.

Let us first consider the {\bf case} $\boldsymbol{r<1+3 m}$. In the following, by $C(m)$ we denote constants which depend only on $m$, and which might change from line to line. In the previously indicated range, also taking into account inequality \eqref{est:DG-pre}, we obtain for the second integral on the right-hand side of \eqref{en:start} that
\begin{align*}%\label{est:|u|}
    \iint_{Q_n}
    |u|
    \big(|u|^m-k_{n+1}^m\big)_+^3\,\dx\dt
    &\le 
    \iint_{Q_n}
    \bigg[\frac{2^{n+1}}{m}\big(|u|^m-k_n^m\big)\bigg]^{\frac1m}  
    \big(|u|^m-k_{n+1}^m\big)_+^3\,\dx\dt \nonumber\\
    &\le C(m)\,
    2^\frac{n+1}{m}
     \iint_{Q_n}
    \big(|u|^m-k_n^m\big)_{+}^{3+\frac{1}{m}}\,\dx\dt.
\end{align*}
The estimate for the first integral on the right-hand side of \eqref{en:start} is similar. Namely,
\begin{align}\label{est:|u|-4m}
     \iint_{Q_n}|u|^{4 m} \mathbf 1_{\{|u|>k_{n+1}\}}\,\dx\dt
     &\le
     \frac{1}{k_{n+1}^{1-m}}
     \iint_{Q_n}|u|^{3 m+1} \mathbf 1_{\{|u|>k_{n+1}\}}\,\dx\dt \nonumber\\
     &\le
     \frac{2^{1-m}}{k^{1-m}}
     \iint_{Q_n}\bigg[\frac{2^{n+1}}{m}\big(|u|^m-k_n^m\big)\bigg]^{3+\frac{1}{m}} \,\dx\dt \nonumber\\
     &\le
     C(m)\frac{2^{(n+1)(3+\frac{1}{m})}}{k^{1-m}}
     \iint_{Q_n}
    \big(|u|^m-k_n^m\big)_{+}^{3+\frac{1}{m}}\, \dx\dt.
\end{align}
The same arguments can  be used to estimate the measure 
of the super-level set $Q_n\cap \{|u|>k_{n+1}\}$. Indeed, we have
\begin{align}\label{est:new-meas}
    \big|Q_n\cap \{|u|>k_{n+1}\}\big|
     &\le
     \frac{1}{k_{n+1}^{3m+1}}
     \iint_{Q_n}|u|^{3 m+1} \mathbf 1_{\{|u|>k_{n+1}\}}\,\dx\dt \nonumber\\
     &\le
     \frac{2^{3m+1}}{k^{3m+1}}
     \iint_{Q_n}\bigg[\frac{2^{n+1}}{m}\big(|u|^m-k_n^m\big)\bigg]^{3+\frac{1}{m}} \,\dx\dt \nonumber\\
    &\le
    C(m) \frac{2^{(n+1)(3+\frac{1}{m})}}{k^{m(3+\frac{1}{m})}}
    \iint_{Q_n}
    \big(|u|^m-k_n^m\big)_+^{3+\frac{1}{m}}\,\dx\dt.
\end{align}
Finally, we treat the integral containing the inhomogeneity $F$. Using H\"older's inequality with exponents
$p$, $2$, and $\frac{2p}{p-2}$ we get
\begin{align*}
    \iint_{Q_n}&|F|^2|u|^{2 m} \mathbf 1_{\{|u|>k_{n+1}\}}\,\dx\dt\\
    &\le
    \bigg[
    \iint_{Q_n}|F|^{2p}\,\dx\dt
    \bigg]^\frac{1}{p}
    \bigg[
    \iint_{Q_n}|u|^{4 m} \mathbf 1_{\{|u|>k_{n+1}\}}\,\dx\dt
    \bigg]^\frac12
    \big|Q_n\cap \{|u|>k_{n+1}\}\big|^{\frac12-\frac{1}{p}}\\
    &\le 
    C(m)\|F\|_{2p, Q_o}^{2}
    \frac{2^{(n+1)(3+\frac{1}{m})(1-\frac{1}{p})}}{k^{1+m-\frac{m}{p}(3+\frac{1}{m})}}
    \bigg[
    \iint_{Q_n}
    \big(|u|^m-k_n^m\big)_+^{3+\frac{1}{m}}\,\dx\dt    
    \bigg]^{1-\frac{1}{p}}.
\end{align*}
To obtain the last line we used \eqref{est:|u|-4m} and \eqref{est:new-meas}. Collecting all the terms in \eqref{en:start}, and taking into account that  $\zeta_n^2\ge \zeta_n^{3+\frac{1}{m}}$, we obtain
\begin{align*}
     \sup_{-\vartheta_n \le \tau \le 0} &
    \int_{B_n \times\{\tau\}}
    \big(|u|^m-k_{n+1}^m\big)_+^{3+\frac{1}{m}} \zeta_n^{3+\frac{1}{m}}\,\dx 
    +
    \iint_{Q_n}\big|\nabla\big[\big(|u|^m-k_{n+1}^m\big)_+^2 \zeta^2_n
    \big]\big|^2\,\dx\dt\\
    &
    \le
    C\Bigg[ \frac{2^{(n+1)(5+\frac{1}{m})}}{(1-\sigma)^2\rho^2} \frac{1}{k^{1-m}}
     \iint_{Q_n}
    \big(|u|^m-k_n^m\big)_{+}^{3+\frac{1}{m}}\, \dx\dt\\
    &\phantom{\le C\bigg[}
    +\frac{2^{(n+1)(1+\frac{1}{m})}}{(1-\sigma)\vartheta}
     \iint_{Q_n}
    \big(|u|^m-k_n^m\big)_{+}^{3+\frac{1}{m}}
   \,\dx\dt\\
    &\phantom{\le C\bigg[}+
    \|F\|_{2p, Q_o}^{2}\frac{2^{(n+1)(3+\frac{1}{m})(1-\frac{1}{p})}}{k^{1+m-\frac{m}{p}(3+\frac{1}{m})}}
    \bigg[
    \iint_{Q_n}
    \big(|u|^m-k_n^m\big)_+^{3+\frac{1}{m}}\,\dx\dt    
    \bigg]^{1-\frac{1}{p}}\Bigg]\\
    &\le
    C 2^{(n+1)(5+\frac{1}{m})}
    \Bigg[
    \bigg[
    \frac{\vartheta}{k^{1-m}\rho^2}+1
    \bigg] \frac{k^{m(3+\frac{1}{m})}}{(1-\sigma)^2\vartheta}\widetilde{\boldsymbol Y}_n + 
    k^{2m} \|F\|_{2p, Q_o}^{2}
    \widetilde{\boldsymbol Y}_n^{1-\frac{1}{p}} \Bigg],
\end{align*}
for a constant $C=C(m,\nu,L)$. In turn we used the abbreviation
\begin{equation*}
    \widetilde{\boldsymbol Y}_n
    :=
    \frac{1}{k^{m(3+\frac{1}{m})}}
    \iint_{Q_n}
    \big(|u|^m-k_n^m\big)_+^{3+\frac{1}{m}}\,\dx\dt.
\end{equation*}
Next, we choose $k$ large enough so that
\begin{equation}\label{k-1-1}
    \frac{1}{k^{1-m}} \frac{\vartheta}{\rho^2} \le 1 
    \quad 
    \Longleftrightarrow 
    \quad 
    k \ge\Big(\frac{\vartheta}{\rho^2}\Big)^{\frac{1}{1-m}}.
\end{equation}
With this choice we arrive at
\begin{align*}%\label{est:en-prel}
     \sup_{-\vartheta_n \le \tau \le 0} &
    \int_{B_n \times\{\tau\}}
    \big[\big(|u|^m-k_{n+1}^m\big)_+^2\zeta_n^2\big]^{\frac12(3+\frac{1}{m})} \,\dx \nonumber\\
    &\quad +
    \iint_{Q_n}\big|\nabla\big[\big(|u|^m-k_{n+1}^m\big)_+^2 \zeta_n^2
    \big]\big|^2\,\dx\dt \nonumber\\
    &\le
    C 2^{(n+1)(5+\frac{1}{m})}
    \bigg[\frac{k^{m(3+\frac{1}{m})}}{(1-\sigma)^2\vartheta}
    \widetilde{\boldsymbol Y}_n
    +
    k^{2m} \|F\|_{2p, Q_o}^{2}
    \widetilde{\boldsymbol Y}_n^{1-\frac{1}{p}}\bigg] 
    =:
    \boldsymbol{F} .
\end{align*}
Next, we apply Gagliardo-Nirenberg's inequality \cite[Chapter~I, Proposition~3.1]{DiBe} with $(v,q,p,m)$ replaced by
\[
    \bigg(\big(|u|^m-k_{n+1}^m\big)_{+}^2 \zeta_n^2,\, 2\frac{N+\frac12 (3 +\frac{1}{m})}{N}, \, 2, \, \tfrac12 (3 +\tfrac{1}{m})\bigg) ,
\]
and obtain
\begin{align}\label{est:q-int}\nonumber
    \iint_{Q_n}&\big[\big(|u|^m-k_{n+1}^m\big)_+^2 \zeta_n^2\big]^q\,\dx\dt\\\nonumber
    &\le
    C_{\rm Sob}^q 
    \bigg[\iint_{Q_n}\big|\nabla\big[
    \big(|u|^m-k^m_{n+1}\big)_+^2\zeta_n^2
    \big]\big|^2\,\dx\dt\bigg]\\\nonumber
    &\qquad\,\,\,\cdot 
    \bigg[\sup _{-\vartheta_n\le \tau \le 0}  
    \int_{B_n\times\{\tau\}}
    \big[\big(|u|^m-k_{n+1}^m\big)_+^2\zeta^2_n\big]^{\frac12 (3+\frac1m) }\,\dx\bigg]^{\frac2N} \\
    &\le
    C \boldsymbol{F}^{1+\frac2N},
\end{align}
where $q=2\frac{N+\frac12 (3 +\frac{1}{m})}{N}$. 
Using the assumption $\boldsymbol\lambda_r=N(m-1)+2r>0$ we deduce that
$q>\frac12 (3+\frac1m)$. In fact,
\begin{align*}
    q
    &=
    \frac{4Nm+2(3m+1)}{2Nm}
    =
    \frac12\Big( 3+\frac1m\Big) + \frac{N(m-1)+2(3m+1)}{2Nm}\\
    &\ge 
     \frac12\Big( 3+\frac1m\Big) + \frac{N(m-1)+2r}{2Nm}
     >
     \frac12\Big( 3+\frac1m\Big)=:\bom.
\end{align*}
The abbreviation $\bom$ proves to be useful, as it avoids excessive exponents in the calculations.
The fact that $q>\bom$ allows us to apply H\"older's inequality first, then the measure estimate \eqref{est:new-meas}, and finally \eqref{est:q-int} to conclude that
\begin{align*}
    \iint_{Q_{n+1}}&
    \big[\big(|u|^m-k_{n+1}^m\big)_{+}^2\big]^{\bom }\,\dx\dt 
    \le 
    \iint_{Q_n}
    \big[\big(|u|^m-k_{n+1}^m\big)_{+}^2 \zeta_n^2\big]^{\bom}\,\dx\dt \\
    &
    \le
    \bigg[\iint_{Q_n}
    \big[\big(|u|^m-k_{n+1}^m\big)_+^2\zeta_n^2\big]^q
    \,\dx\dt\bigg]^{\frac{\bom}{q}}
    \big|Q_n\cap \{|u|>k_{n+1}\} \big|^{1- \frac{\bom}{q}}\\
    &\le
    C \boldsymbol{F}^{\frac{\bom}{q}(1+\frac2N)}
    \Big[ 
    2^{(n+1)(3+\frac{1}{m})}\widetilde{\boldsymbol Y}_n
    \Big]^{1- \frac{\bom}{q}}\\
    &\le
    C 
    \Big[ 
    2^{2 (n+1)\bom}\widetilde{\boldsymbol Y}_n
    \Big]^{1- \frac{\bom}{q}}\\
    &\phantom{\le\,\,\, }\cdot
    2^{2(n+1)(1+\bom)\frac{\bom}{q}(1+\frac2N)}
    \bigg[\frac{k^{2m\bom}}{(1-\sigma)^2\vartheta}
    \widetilde{\boldsymbol Y}_n
    +
    k^{2m} \|F\|_{2p, Q_o}^{2}
    \widetilde{\boldsymbol Y}_n^{1-\frac{1}{p}}\bigg]
    ^{\frac{\bom}{q}(1+\frac2N)}.
\end{align*}
Since the power of $2$ in the last estimate can be estimated by
\begin{align*}
    2(n+1)\bom\Big(1- \frac{\bom}{q}\Big) +& 
    2(n+1)(1+\bom)\frac{\bom}{q}\Big(1+\frac2N\Big)\\
    &\le 2(n+1)\bom+ 2(n+1)(1+\bom)3\\
    &\le 8(n+1)\bom,
\end{align*}
the inequality becomes
\begin{align*}
    \iint_{Q_{n+1}}&
    \big[\big(|u|^m-k_{n+1}^m\big)_{+}^2\big]^{\bom }\,\dx\dt \\
    &\le
    C 2^{8(n+1)\bom}
    \widetilde{\boldsymbol Y}_n
    ^{1- \frac{\bom}{q}}
    \bigg[\frac{k^{2m\bom}}{(1-\sigma)^2\vartheta}
    \widetilde{\boldsymbol Y}_n
    +
    k^{2m}\|F\|_{2p, Q_o}^{2}
    \widetilde{\boldsymbol Y}_n^{1-\frac{1}{p}}\bigg]
    ^{\frac{\bom}{q}(1+\frac2N)} \\
    &=
    C 2^{8(n+1)\bom}
    k^{2m\bom \frac{\bom}{q}(1+\frac2N)}
    \widetilde{\boldsymbol Y}_n^{1+\frac{2\bom}{Nq}(1-\frac{N+2}{2p})}\\
    &\phantom{\le\,\,\, }\cdot
    \bigg[
    \frac{\widetilde{\boldsymbol Y}_n^\frac{1}{p}}{(1-\sigma)^2\vartheta}
    +
    \|F\|_{2p, Q_o}^{2}
    k^{2m(1-\bom)}
    \bigg]^{\frac{\bom}{q}(1+\frac{2}{N})}.
\end{align*}
We divide both sides by $k^{2m\bom}$. Taking into account 
\begin{align*}
    &2m\bom \frac{\bom}{q}\Big(1+\frac{2}{N}\Big)-2m\bom
    =
    \frac{2m\bom}{Nq}
    \big[\bom (N+2)-Nq\big]\\
    &\qquad=
    \frac{2m\bom}{Nq}
    \big[\bom (N+2)-2(N+\bom)\big]
    =
    \frac{2m\bom}{Nq}
    \big[\bom N -2N\big]\\
    &\qquad=
    \frac{m\bom}{q} (2\bom -4)
    =
    \frac{\bom}{q} (1-m),
\end{align*}
we obtain
\begin{align*}
    \widetilde{\boldsymbol Y}_{n+1}
    &\le
     C \boldsymbol b^{n+1}
     k^{\frac{\bom}{q}(1-m)}
    \widetilde{\boldsymbol Y}_n^{1+\frac{2\bom}{Nq}(1-\frac{N+2}{2p})}
    \bigg[
    \frac{\widetilde{\boldsymbol Y}_n}{(1-\sigma)^{2p}\vartheta^p}
    +
    \frac{\|F\|_{2p, Q_o}^{2p}}{k^{p(1+m)}}
    \bigg]^{\frac{2\bom}{Nq}\frac{N+2}{2p}},
\end{align*}
where
$\boldsymbol b:=2^{8\bom}$.
Taking integral means on both sides, and defining
\[
\boldsymbol Y_{n}:=\frac{1}{k^{m(3+\frac{1}{m})}}
    \biint_{Q_n}
    \big(|u|^m-k_n^m\big)_+^{3+\frac{1}{m}}\,\dx\dt,
\]
we arrive at
\begin{align*}
    \boldsymbol Y_{n+1}
    &\le
    C \boldsymbol b^{n}
    k^{\frac{\bom}{q}(1-m)}
    |Q_o|^\frac{2\bom}{Nq}
    \bigg[
    \frac{\boldsymbol Y_n}{(1-\sigma)^{2p}
    \vartheta^p}
    +\frac{\|F\|_{2p, Q_o}^{2p}}{|Q_o|k^{p(1+m)}}
    \bigg]^{\frac{2\bom}{Nq}\frac{N+2}{2p}}
    \boldsymbol Y_n^{1+\frac{2\bom}{Nq}(1-\frac{N+2}{2p})} \\
    &\le
    C \boldsymbol b^{n}
    k^{\frac{\bom}{q}(1-m)}
    |Q_o|^\frac{2\bom}{Nq}
    \bigg[
    \frac{2^{N+1}\boldsymbol Y_o}{(1-\sigma)^{2p}
    \vartheta^p}
    +\frac{\|F\|_{2p, Q_o}^{2p}}{|Q_o|k^{p(1+m)}}
    \bigg]^{\frac{2\bom}{Nq}\frac{N+2}{2p}}
    \boldsymbol Y_n^{1+\frac{2\bom}{Nq}(1-\frac{N+2}{2p})}.
\end{align*} 
At this stage we choose $k$ large enough such that we have
\begin{align}\label{k-1-2}
    \boldsymbol Y_{o}
    &\le
    \frac{\vartheta^\frac{N+2}2(1-\sigma)^{N+2}}{|Q_o|k^{\frac{N}{2}(1-m)}}
\end{align}
and 
\begin{equation}\label{k-1-3}
    k^{\frac12\boldsymbol \lambda_{p(1+m)}} 
    \ge
    \Big(\frac{\vartheta}{\varrho^2}\Big)^{p-\frac{N}2}
    \biint_{Q_o}|\varrho F|^{2p}\,\dx\dt.
\end{equation}
This is possible, since $\boldsymbol{\lm}_{3m+1}>\boldsymbol{\lm}_{r}>0$ and $\boldsymbol\lambda_{p(1+m)}>0$. Note that \eqref{k-1-3} implies
\begin{align*}
     \frac{\|F\|_{2p, Q_o}^{2p}}{k^{p(1+m)}}
     &\le
     \frac{|Q_o|}{\varrho^{2p}k^{\frac{N}{2}(1-m)}} 
     \Big(\frac{\varrho^2}{\vartheta}\Big)^{p-\frac{N}2} 
     = 
     \frac{\vartheta^{\frac{N+2}2-p}|B_1|}{k^{\frac{N}{2}(1-m)}} ,
\end{align*}
where $|B_1|$ is the measure of $B_1$. 
Such a choice of $k$ yields the recursive inequality
\begin{align*}
    \boldsymbol Y_{n+1}
    &\le
     C \boldsymbol b^{n}
    k^{\frac{\bom}{q}(1-m)}
    |Q_o|^\frac{2\bom}{Nq}\\
    &\phantom{\le\,\,\, }\cdot
    \bigg[
     \frac{\vartheta^{\frac{N+2}2-p}}{(1-\sigma)^{2p-(N+2)}|Q_o| k^{\frac{N}{2}(1-m)}}
    +\frac{\vartheta^{\frac{N+2}2-p}|B_1|}{|Q_o|k^{\frac{N}{2}(1-m)}}
    \bigg]^{\frac{2\bom}{Nq}\frac{N+2}{2p}}
    \boldsymbol Y_n^{1+\frac{2\bom}{Nq}(1-\frac{N+2}{2p})} \\
    &\le
     C \boldsymbol b^{n}
    k^{\frac{\bom}{q}(1-m)(1-\frac{N+2}{2p})}
    \bigg[\frac{|Q_o|}{(1-\sigma)^{N+2}\vartheta^\frac{N+2}2}\bigg]^{\frac{2\bom}{Nq}(1-\frac{N+2}{2p})}
    \boldsymbol Y_n^{1+\frac{2\bom}{Nq}(1-\frac{N+2}{2p})}.
\end{align*}
Apply the fast geometric convergence \cite[\S\,2, Lemma 5.1]{DBGV-book} with
$$
    \alpha:=\frac{2\bom}{Nq}\Big(1-\frac{N+2}{2p}\Big)
    \quad
    \mbox{and}
    \quad
    C\equiv 
    C  
     k^{\frac{N}{2}\alpha(1-m)}
     \bigg[\frac{|Q_o|}{(1-\sigma)^{N+2}\vartheta^{\frac{N+2}2}}
     \bigg]^{\alpha}
$$
to conclude that $\boldsymbol Y_n\to 0$ as $n\to\infty$ if we require that
\begin{align*}
    \boldsymbol Y_o=
    \frac{1}{k^{2m\bom}}
    \biint_{Q_o}|u|^{2m\bom}\,\dx\dt
    &\le
    \bigg[
    C  
     k^{\frac{N}{2}\alpha(1-m)}
     \bigg[\frac{|Q_o|}{(1-\sigma)^{N+2}\vartheta^{\frac{N+2}2}}\bigg]^{\alpha}
    \bigg]^{-\frac{1}{\alpha}}\boldsymbol b^{-\frac1{\alpha^2}}\\
    &\equiv
    C^{-1}(1-\sigma)^{N+2}
    \frac{\vartheta^{\frac{N+2}2}}{|Q_o|}
    k^{\frac{N}{2}(m-1)}.
\end{align*}
Solving this inequality for $k$ gives
\begin{align}\label{k-1-4}
    k^{\frac12\boldsymbol{\lambda}_{3m+1}}
    =
    k^{\frac{N}{2}(m-1)+2m\bom}
    \ge
    \frac{C}{(1-\sigma)^{N+2}}
    \Big(\frac{\varrho^2}{\vartheta}\Big)^\frac{N}{2}
    \biint_{Q_o}|u|^{3m+1}\,\dx\dt.
\end{align}
As a result of $\boldsymbol Y_\infty =0$, we find that
\begin{equation*}
    \sup_{Q_{\sigma\rho,\sigma\vartheta}}  |u|
    \le k.
\end{equation*}
Choosing $k$ according to \eqref{k-1-1}, \eqref{k-1-2}, \eqref{k-1-3}, \eqref{k-1-4} as
\begin{align*}
    k
    &:=
    \max\Bigg\{
    \bigg[
     \frac{C}{(1-\sigma)^{N+2}}
    \Big(\frac{\varrho^2}{\vartheta}\Big)^\frac{N}{2}
    \biint_{Q_o}|u|^{3m+1}\,\dx\dt
    \bigg]^\frac{2}{\boldsymbol\lm_{3m+1}},\\
    &\qquad\qquad\;\,
    \bigg[
    \Big(\frac{\vartheta}{\varrho^2}\Big)^{p-\frac{N}2}
    \biint_{Q_o}|\varrho F|^{2p}\,\dx\dt\bigg]^\frac{2}{\boldsymbol \lambda_{p(1+m)}},
    \Big(
    \frac{\vartheta}{\rho^2}
    \Big)^{\frac{1}{1-m}}
    \Bigg\}
\end{align*}
we get
\begin{align*}
    \sup_{Q_{\sigma\rho,
    \sigma\vartheta}}
    |u|
    &\le
    \max\Bigg\{
    \Big(
    \sup_{Q_{\rho,\vartheta}}|u|
    \Big)^{1-\frac{\boldsymbol\lm_r}{\boldsymbol\lm_{3m+1}}}
    \bigg[
     \frac{C}{(1-\sigma)^{N+2}}
    \Big(\frac{\varrho^2}{\vartheta}\Big)^\frac{N}{2}
    \biint_{Q_{\varrho,\vartheta}}|u|^{r}\,\dx\dt
    \bigg]^\frac{2}{\boldsymbol\lm_{3m+1}},\\
    &\qquad\qquad\;\,
    \bigg[
    \Big(\frac{\vartheta}{\varrho^2}\Big)^{p-\frac{N}2}
    \biint_{Q_{\varrho,\vartheta}}|\varrho F|^{2p}\,\dx\dt\bigg]^\frac{2}{\boldsymbol \lambda_{p(1+m)}},
    \Big(
    \frac{\vartheta}{\rho^2}
    \Big)^{\frac{1}{1-m}}
    \Bigg\}.
\end{align*}
The $sup$-term in the first entry of the maximum can be eliminated with the help of a standard interpolation argument which leads to the final quantitative $L^\infty$-estimate.  The precise argument is as follows.
For $\tau\in[\tfrac12,1)$, we define cylinders $\widehat Q_n=\widehat Q_{\rho_n,\theta_n}$ with
\begin{equation*}%\label{interpol-rho-n}
    \rho_n:=\tau \rho +(1-\tau)\rho \sum_{j=1}^n2^{-j}
    \quad\mbox{and}\quad
    \theta_n:=\tau \vartheta +(1-\tau)\vartheta \sum_{j=1}^n2^{-j},
\end{equation*}
and
\begin{equation}\label{interpol-M-n}
     \boldsymbol M_n
    :=
    \sup_{\widehat Q_n} |u|.
\end{equation}
Note that $\widehat  Q_o=Q_{\tau\rho,\tau\vartheta}(z_o)$ and $\widehat Q_\infty = Q_{\rho,\vartheta}(z_o)$. 
We apply the above estimate on $\widehat Q_n$, $\widehat Q_{n+1}$ instead of
$Q_{\sigma\rho, \sigma\vartheta}(z_o)$, $Q_{\rho,\vartheta}(z_o)$ so that  $1-\frac{\rho_n}{\rho_{n+1}}=
1-\frac{\vartheta_n}{\vartheta_{n+1}}$ plays the role of $1-\sigma$.  A simple computation yields the lower bound $1-\frac{\rho_n}{\rho_{n+1}}\ge \frac{1-\tau}{2^{n+1}}$. Moreover, we have
\begin{equation*}%\label{interpol-F-n}
    \Big(\frac{\vartheta_{n+1}}{\varrho_{n+1}^2}\Big)^{p-\frac{N}2}
    \biint_{\widehat Q_{n+1}}|\varrho_{n+1} F|^{2p}\,\dx\dt
    \le
    C(N,p) 
    \Big(\frac{\vartheta}{\varrho^2}\Big)^{p-\frac{N}2}
    \biint_{Q_{\varrho,\vartheta}}|\varrho F|^{2p}\,\dx\dt,
\end{equation*}
and
\begin{equation*}%\label{interpol-|u|^r-n}
    \Big(\frac{\varrho_{n+1}^2}{\vartheta_{n+1}}\Big)^\frac{N}{2}
    \biint_{\widehat Q_{n+1}}|u|^{r}\,\dx\dt
    \le
    2^{N+1}
    \Big(\frac{\varrho^2}{\vartheta}\Big)^\frac{N}{2}
    \biint_{Q_{\varrho,\vartheta}}|u|^{r}\,\dx\dt.
\end{equation*}
Inserting this above, we get
\begin{align*}
    \boldsymbol M_n
    &\le 
    \max\Bigg\{
    \boldsymbol M_{n+1}^{1-\frac{\boldsymbol\lm_r}{\boldsymbol\lm_{3m+1}}}
    \bigg[
     \frac{C2^{2Nn}}{(1-\tau)^{N+2}}
    \Big(\frac{\varrho^2}{\vartheta}\Big)^\frac{N}{2}
    \biint_{Q_{\varrho,\vartheta}}|u|^{r}\,\dx\dt
    \bigg]^\frac{2}{\boldsymbol\lm_{3m+1}},\\
    &\qquad\qquad\;\,
    \bigg[ C(N,p)
    \Big(\frac{\vartheta}{\varrho^2}\Big)^{p-\frac{N}2}
    \biint_{Q_{\varrho,\vartheta}}|\varrho F|^{2p}\,\dx\dt\bigg]^\frac{2}{\boldsymbol \lambda_{p(1+m)}},
    2^\frac{1}{1-m}\Big(
    \frac{\vartheta}{\rho^2}
    \Big)^{\frac{1}{1-m}}
    \Bigg\}.
\end{align*}
If for some $n\in \N_0$ the second or third term in the maximum dominates the first, there is nothing to prove, because we trivially have
\begin{align*}
    \sup_{Q_{\tau\rho,\tau\vartheta}} |u| 
    &=
   \boldsymbol  M_o
    \le 
   \boldsymbol  M_{n}\\
   &\le
   \max\Bigg\{
   \bigg[ C
    \Big(\frac{\vartheta}{\varrho^2}\Big)^{p-\frac{N}2}
    \biint_{Q_{\varrho,\vartheta}}|\varrho F|^{2p}\,\dx\dt\bigg]^\frac{2}{\boldsymbol \lambda_{p(1+m)}},
    2^\frac{1}{1-m}\Big(
    \frac{\vartheta}{\rho^2}
    \Big)^{\frac{1}{1-m}}
    \Bigg\}.
\end{align*}
Otherwise, the first term is dominant for all $n\in\N_0$, so that for any $n\in\N_0$ we have
\begin{equation*}
    \boldsymbol  M_{n}
    \le
    \underbrace{
    \bigg[
     \frac{C}{(1-\tau)^{N+2}}
    \Big(\frac{\varrho^2}{\vartheta}\Big)^\frac{N}{2}
    \biint_{Q_{\varrho,\vartheta}}|u|^{r}\,\dx\dt
    \bigg]^\frac{2}{\boldsymbol\lm_{3m+1}}}_{=:\mathfrak C}
    2^\frac{4Nn}{\boldsymbol\lm_{3m+1}}
    \boldsymbol M_{n+1}^{1-\frac{\boldsymbol\lm_r}{\boldsymbol\lm_{3m+1}}}.
\end{equation*}
Apply \cite[Chapter I, Lemma 4.3]{DiBe} with 
\begin{equation*}
    \alpha:= \frac{\boldsymbol\lm_r}{\boldsymbol\lm_{3m+1}}
    \quad
    \mbox{and}
    \quad
    \boldsymbol b
    := 2^\frac{4N}{\boldsymbol\lm_{3m+1}},
\end{equation*}
to conclude that
\begin{align*}
    \boldsymbol M_o
    &\le
    %\bigg(\frac{2\boldsymbol b\mathfrak C}{\boldsymbol b^{1-%\frac1{\alpha}}}
    %\bigg)^\frac1{\alpha}
    %\textcolor{blue}{
    \bigg(\frac{2\mathfrak C}{\boldsymbol b^{1-\frac1{\alpha}}}
    \bigg)^\frac1{\alpha}
    %}
    =
    %\big( 2 \boldsymbol b^{\frac1{\alpha}}
    %\mathfrak C\big)^\frac1{\alpha}
    %\equiv
    \bigg[
     \frac{C}{(1-\tau)^{N+2}}
    \Big(\frac{\varrho^2}{\vartheta}\Big)^\frac{N}{2}
    \biint_{Q_{\varrho,\vartheta}}|u|^{r}\,\dx\dt
    \bigg]^\frac{2}{\boldsymbol\lm_{r}}.
\end{align*}
Joining this with the first case,
we end up with
\begin{align*}
    \sup_{Q_{\tau\rho,
    \tau\vartheta}}
    |u|
    &\le
    \max\Bigg\{
    \bigg[
     \frac{C}{(1-\tau)^{N+2}}
    \Big(\frac{\varrho^2}{\vartheta}\Big)^\frac{N}{2}
    \biint_{Q_{\varrho,\vartheta}}|u|^{r}\,\dx\dt
    \bigg]^\frac{2}{\boldsymbol\lm_{r}},\\
    &\qquad\qquad\;\,
    \bigg[ C
    \Big(\frac{\vartheta}{\varrho^2}\Big)^{p-\frac{N}2}
    \biint_{Q_{\varrho,\vartheta}}|\varrho F|^{2p}\,\dx\dt\bigg]^\frac{2}{\boldsymbol \lambda_{p(1+m)}},2^\frac{1}{1-m}
    \Big(
    \frac{\vartheta}{\rho^2}
    \Big)^{\frac{1}{1-m}}
    \Bigg\}.
\end{align*}
The constants $C=C(N,p)$ in the second and $2^\frac{1}{1-m}$ in the third entry of the maximum can be replaced by 1 (by modifying \eqref{k-1-1} and \eqref{k-1-3} accordingly). This finishes the proof in the case $r<3m+1$.

Now we consider the {\bf case} $\boldsymbol{r\ge1+3 m}$. The starting point is again the energy inequality \eqref{en:start}, combined with inequality \eqref{est:DG-pre} on the super-level sets $\{|u|>k_{n+1}\}$. In fact, for the first term on the right-hand side of \eqref{en:start} we obtain
\begin{align*}
     \iint_{Q_n}|u|^{4 m} \mathbf 1_{\{|u|>k_{n+1}\}}\,\dx\dt
     &\le
     \frac{1}{k_{n+1}^{r-4m}}
     \iint_{Q_n}|u|^{r} \mathbf 1_{\{|u|>k_{n+1}\}}\,\dx\dt \nonumber\\
     &\le
     \frac{2^{r-4m}}{k^{r-4m}}
     \iint_{Q_n}\bigg[\frac{2^{n+1}}{m}\big(|u|^m-k_n^m\big)\bigg]^{\frac{r}{m}} \,\dx\dt \nonumber\\
     &\le
    C(m,r)
    \frac{2^{(n+1)\frac{r}{m}}}{k^{r-4m}}
     \iint_{Q_n}
    \big(|u|^m-k_n^m\big)_{+}^{\frac{r}{m}}
    \,\dx\dt.
\end{align*}
Similarly, we get for the second term on the right-hand side of \eqref{en:start} that
\begin{align*}
    &\iint_{Q_n}
    |u|
    \big(|u|^m-k_{n+1}^m\big)_+^3\,\dx\dt \\
    &\qquad\le 
    \frac{1}{k_{n+1}^{r-3m-1}} \iint_{Q_n}
    |u|^{r-3m}
    \big(|u|^m-k_{n+1}^m\big)_+^3\,\dx\dt \\
    &\qquad\le 
    \frac{2^{r-3m-1}}{k^{r-3m-1}} \iint_{Q_n}
    \bigg[\frac{2^{n+1}}{m}\big(|u|^m-k_n^m\big)\bigg]^{r-3m}
    \big(|u|^m-k_{n+1}^m\big)_+^3\,\dx\dt \\
    &\qquad\le 
    C(m,r) \frac{2^{(n+1)(\frac{r}{m}-3)}}{k^{r-3m-1}}
    \iint_{Q_n}
    \big(|u|^m-k_n^m\big)_{+}^{\frac{r}{m}}
    \,\dx\dt.
\end{align*}
We can use the same line of 
reasoning to estimate the measure 
of the super-level set $Q_n\cap \{|u|>k_{n+1}\}$. In fact, 
\begin{align}\label{est:new-meas-r}
    \big|Q_n\cap \{|u|>k_{n+1}\}\big|
     &\le
     \frac{1}{k_{n+1}^{r}}
     \iint_{Q_n}|u|^{r} \mathbf 1_{\{|u|>k_{n+1}\}}\,\dx\dt \nonumber\\
     &\le
     \frac{2^{r}}{k^{r}}
     \iint_{Q_n}\bigg[\frac{2^{n+1}}{m}\big(|u|^m-k_n^m\big)\bigg]^{\frac{r}{m}} \,\dx\dt \nonumber\\
    &\le
    C(m,r) \frac{2^{(n+1)\frac{r}{m}}}{k^{r}}
    \iint_{Q_n}
    \big(|u|^m-k_n^m\big)_+^{\frac{r}{m}}\,\dx\dt.
\end{align}
Finally, we have
\begin{align*}
    \iint_{Q_n}&|F|^2|u|^{2 m} \mathbf 1_{\{|u|>k_{n+1}\}}\,\dx\dt\\
    &\le
    \bigg[
    \iint_{Q_n}|F|^{2p}\,\dx\dt
    \bigg]^\frac{1}{p}
    \bigg[
    \iint_{Q_n}|u|^{4 m} \mathbf 1_{\{|u|>k_{n+1}\}}\,\dx\dt
    \bigg]^\frac12
    \big|Q_n\cap \{|u|>k_{n+1}\}\big|^{\frac12-\frac{1}{p}}\\
    &\le 
    C\|F\|_{2p, Q_o}^{2}
    \bigg[
    \frac{2^{(n+1)\frac{r}{m}}}{k^{r-4m}}
     \iint_{Q_n}
    \big(|u|^m-k_n^m\big)_{+}^{\frac{r}{m}}
    \,\dx\dt
    \bigg]^\frac12\\
    &\qquad\qquad\,\,\,\cdot
    \bigg[
     \frac{2^{(n+1)\frac{r}{m}}}{k^{r}}
    \iint_{Q_n}
    \big(|u|^m-k_n^m\big)_+^{\frac{r}{m}}\,\dx\dt
    \bigg]^{\frac12 -\frac1p}\\
    &\le
    C2^{(n+1)\frac{r}{m}(1-\frac{1}{p})}k^{2m}\|F\|_{2p, Q_o}^{2}
    \bigg[
    \frac{1}{k^{r}}
     \iint_{Q_n}
    \big(|u|^m-k_n^m\big)_{+}^{\frac{r}{m}}
    \,\dx\dt
    \bigg]^{1-\frac{1}{p}},
\end{align*}
where $C=C(m,r)$.
Joining the preceding estimates with the energy inequality \eqref{en:start} yields 
\begin{align*}
    &\sup_{-\vartheta_n \le \tau \le 0}
    \int_{B_n \times\{\tau\}}
    \big(|u|^m-k_{n+1}^m\big)_+^{3+\frac{1}{m}} \zeta_n^2\,\dx 
    +
    \iint_{Q_n}\big|\nabla\big[\big(|u|^m-k_{n+1}^m\big)_+^2 \zeta_n
    \big]\big|^2\,\dx\dt\\\nonumber
    &\quad
    \le
    C\bigg[\frac{2^{(n+1)(2+\frac{r}{m})}k^{4m}}{(1-\sigma)^2\rho^2}
    \frac{1}{k^{r}}
     \iint_{Q_n}
    \big(|u|^m-k_n^m\big)_{+}^{\frac{r}{m}}
    \,\dx\dt\\
    &\quad\phantom{\le C\bigg[}
    + \frac{2^{(n+1)(\frac{r}{m}-2)}k^{3m+1}}{(1-\sigma)\vartheta}
        \frac{1}{k^{r}}
       \iint_{Q_n}
    \big(|u|^m-k_n^m\big)_{+}^{\frac{r}{m}}
    \,\dx\dt
    \bigg]
    \\
    &\quad\phantom{\le C\bigg[}
    + 2^{(n+1)\frac{r}{m}(1-\frac{1}{p})}k^{2m}\|F\|_{2p, Q_o}^{2}
    \bigg[
    \frac{1}{k^{r}}
     \iint_{Q_n}
    \big(|u|^m-k_n^m\big)_{+}^{\frac{r}{m}}
    \,\dx\dt
    \bigg]^{1-\frac{1}{p}}\\
    &\quad\le 
    C 
    2^{(n+1)(2+\frac{r}{m})}
    \Bigg[
    \bigg[
    \frac{\vartheta}{k^{1-m}\rho^2}
    +
    1
    \bigg]
    \frac{k^{3m+1}}{(1-\sigma)^2\vartheta} \widetilde{\boldsymbol Y}_n 
    +k^{2m}\|F\|_{2p, Q_o}^{2}
    \widetilde{\boldsymbol Y}_n^{1-\frac{1}{p}}
    \Bigg],
\end{align*}
where we defined
$$
    \widetilde{\boldsymbol Y}_n
    :=
    \frac{1}{k^{r}}
       \iint_{Q_n}
    \big(|u|^m-k_n^m\big)_{+}^{\frac{r}{m}}
    \,\dx\dt.
$$
We stipulate that
\begin{equation}\label{k-2-1}
    \frac{1}{k^{1-m}}\frac{\vartheta}{\rho^2}\le 1
    \quad\Longleftrightarrow\quad
    k\ge \Big(\frac{\vartheta}{\varrho^2}\Big)^\frac{1}{1-m}.
\end{equation}
This simplifies the previous inequality to
\begin{align*}%\label{est:en-r}
    \sup_{-\vartheta_n \le \tau \le 0}&
    \int_{B_n \times\{\tau\}}
    \big[\big(|u|^m-k_{n+1}^m\big)_+^2\zeta_n^2\big]^{\frac12(3+\frac{1}{m})}\,\dx\nonumber\\ 
    +&
    \iint_{Q_n}\big|\nabla\big[\big(|u|^m-k_{n+1}^m\big)_+^2 \zeta_n^2
    \big]\big|^2\,\dx\dt\\ 
    &\qquad\le 
    C
    2^{(n+1)(2+\frac{r}{m})}
    \bigg[
    \frac{k^{3m+1}}{(1-\sigma)^2\vartheta}
    \widetilde{\boldsymbol Y}_n
    +
    k^{2m}\|F\|_{2p, Q_o}^{2}
    \widetilde{\boldsymbol Y}_n^{1-\frac{1}{p}}
    \bigg] 
    =:
    \boldsymbol{F}. \nonumber
\end{align*}
Next, we apply Gagliardo-Nirenberg's inequality \cite[Chapter~I, Proposition~3.1]{DiBe} with $(v,q,p,m)$ replaced  by
\[
    \bigg(\big(|u|^m-k_{n+1}^m\big)_{+}^2 \zeta_n^2,\, 2\frac{N+\frac12 (3 +\frac{1}{m})}{N}, \, 2, \, \tfrac12 (3 +\tfrac{1}{m})\bigg) ,
\]
and obtain
\begin{align}\label{est:q-int-*}\nonumber
    \iint_{Q_n}&\big[\big(|u|^m-k_{n+1}^m\big)_+^2 \zeta_n^2\big]^q\,\dx\dt\\\nonumber
    &\le
    C_{\rm Sob}^q 
    \bigg[\iint_{Q_n}\big|\nabla\big[
    \big(|u|^m-k^m_{n+1}\big)_+^2\zeta_n^2
    \big]\big|^2\,\dx\dt\bigg]\\\nonumber
    &\qquad\,\,\,\,\cdot 
    \bigg[\sup _{-\vartheta_n\le \tau \le 0}  
    \int_{B_n\times\{\tau\}}
    \big[\big(|u|^m-k_{n+1}^m\big)_+^2\zeta_n^2\big]^{\frac12 (3+\frac1m) }\,\dx\bigg]^{\frac2N} \\
    &\le
    C\boldsymbol{F}^{1+\frac2N},
\end{align}
where $q=2\frac{N+\frac12 (3 +\frac{1}{m})}{N}$. 
Observe that we do not have $2q<\frac{r}{m}$ for all admissible values of $r$ and $m$. However, we can ensure that $q<\frac{r}{m}$, which implies $r-mq>0$. Hence, to estimate the left-hand side of the previous inequality from below, we first apply H\"older's inequality to raise the exponent from $\frac{r}{m}$ to $\frac{2r}{m}$, and afterwards we reduce the integrability exponent of
$\big(\big(|u|^m-k_{n+1}^m\big) \zeta_n\big)^\frac{2r}{m}$ to $2q$ by extracting the $L^\infty$-norm of $|u|$. To the second term we apply the measure estimate \eqref{est:new-meas-r}. Finally, we apply \eqref{est:q-int-*}. In this way we get
\begin{align*}
    \iint_{Q_{n+1}}&\big(|u|^m-k_{n+1}^m\big)^\frac{r}
    {m}\,\dx\dt
    \le
    \iint_{Q_{n}}\big[\big(|u|^m-k_{n+1}^m\big)\zeta_n\big]^\frac{r}
    {m}\,\dx\dt\\
    &\le
    \bigg[
     \iint_{Q_{n}}\big[\big(|u|^m-k_{n+1}^m\big)^2\zeta_n^2\big]^\frac{r}
    {m}\,\dx\dt
    \bigg]^\frac12
    \big| Q_n\cap\{|u|>k_{n+1}\}\big|^\frac12\\
    &\le
    C\|u\|_{\infty ,Q_o}^{r-mq}
    \bigg[
     \iint_{Q_{n}}\big[\big(|u|^m-k_{n+1}^m\big)^2\zeta_n^2\big]^q\,\dx\dt
    \bigg]^\frac12
    \big[
     2^{(n+1)\frac{r}{m}}
     \widetilde{\boldsymbol Y}_n
    \big]^\frac12\\
    &\le
     C\|u\|_{\infty ,Q_o}^{r-mq}
     \boldsymbol{F}^{\frac12 (1+\frac2N)}\big[
     2^{(n+1)\frac{r}{m}}
     \widetilde{\boldsymbol Y}_n
    \big]^\frac12\\
    &=
    C\|u\|_{\infty ,Q_o}^{r-mq}\big[
     2^{(n+1)\frac{r}{m}}
     \widetilde{\boldsymbol Y}_n
    \big]^\frac12\\
    &\phantom{=\,}\cdot
    \Bigg[
    2^{(n+1)(2+\frac{r}{m})}
    \bigg[
    \frac{k^{3m+1}}{(1-\sigma)^2\vartheta}
    \widetilde{\boldsymbol Y}_n
    +
    k^{2m}\|F\|_{2p, Q_o}^{2}
    \widetilde{\boldsymbol Y}_n^{1-\frac{1}{p}}
    \bigg]
    \Bigg]^{\frac12 (1+\frac{2}{N})}\\
    &\le
    C\|u\|_{\infty ,Q_o}^{r-mq} 2^{(n+1)(2+\frac{r}{m})(1+\frac{1}{N})}
    k^{\frac12 (3m+1)(1+\frac{2}{N})}\\
    &\phantom{=\,}\cdot
    \widetilde{\boldsymbol Y}_n^{1+\frac{1}{N}(1-\frac{N+2}{2p})}
    \bigg[
    \frac{\widetilde{\boldsymbol Y}_n^\frac{1}{p}}{(1-\sigma)^2\vartheta}
    +
    \frac{\|F\|_{2p, Q_o}^{2}}{k^{1+m}}
    \bigg]^{\frac12 (1+\frac{2}{N})}.
\end{align*}
The exponent of $k$ can be re-written as 
\begin{align*}
    \tfrac12 (3m+1)\big(1+\tfrac{2}{N}\big)
    &=
    \tfrac12(3m+1) +\tfrac{3m+1}{N} \\
    &=
    \tfrac12(3m+1) + m(q-2) 
    =
    mq + \tfrac12(1-m).
\end{align*}
Let
\begin{equation*}
    \boldsymbol b:= 2^{(2+\frac{r}{m})(1+\frac{1}{N})},
\end{equation*}
and divide the last inequality by $k^r$ to obtain
\begin{align*}
    \widetilde{\boldsymbol Y}_{n+1}
    &\le
    C\|u\|_{\infty ,Q_o}^{r-mq}
    k^{mq-r + \frac12(1-m)}
    \boldsymbol b^{n+1}\\
     &\phantom{=\,}\cdot
     \widetilde{\boldsymbol Y}_n^{1+\frac{1}{N}(1-\frac{N+2}{2p})}
    \bigg[
    \frac{\widetilde{\boldsymbol Y}_n}{(1-\sigma)^{2p}\vartheta^p}
    +
    \frac{\|F\|_{2p, Q_o}^{2p}}{k^{p(1+m)}}
    \bigg]^{\frac{1}{N}\frac{N+2}{2p}}.
\end{align*}
Taking means on both sides, and defining $\boldsymbol Y_{n}:=|Q_n|^{-1} \widetilde{\boldsymbol Y}_n$ as in the case $r<3m+1$, we find that
\begin{align*}
    \boldsymbol Y_{n+1}
    &\le
    C\|u\|_{\infty ,Q_o}^{r-mq}
    k^{mq-r + \frac12(1-m)}|Q_o|^\frac{1}{N}
    \boldsymbol b^{n+1}\\
     &\phantom{=\,}\cdot
     \boldsymbol Y_n^{1+\frac{1}{N}(1-\frac{N+2}{2p})}
    \bigg[
    \frac{\boldsymbol Y_n}{(1-\sigma)^{2p}\vartheta^p}
    +
    \frac{\|F\|_{2p, Q_o}^{2p}}{|Q_o|k^{p(1+m)}}
    \bigg]^{\frac{1}{N}\frac{N+2}{2p}}\\
    &\le
    C\|u\|_{\infty ,Q_o}^{r-mq}
    k^{mq-r + \frac12(1-m)}|Q_o|^\frac{1}{N}
    \boldsymbol b^{n+1}\\
     &\phantom{=\,}\cdot
     \boldsymbol Y_n^{1+\frac{1}{N}(1-\frac{N+2}{2p})}
    \bigg[
    \frac{\boldsymbol Y_o}{(1-\sigma)^{2p}\vartheta^p}
    +
    \frac{\|F\|_{2p, Q_o}^{2p}}{|Q_o|k^{p(1+m)}}
    \bigg]^{\frac{1}{N}\frac{N+2}{2p}}.
\end{align*}
As in the case $r<3m+1$ we first choose $k$ large enough such that both
\begin{equation}\label{k-2-2}
    \boldsymbol Y_{o}
    \le
    \frac{\vartheta^\frac{N+2}2(1-\sigma)^{N+2}}{|Q_o|k^{\frac{N}{2}(1-m)}}\,\,\Leftrightarrow\,\,k\ge\bigg[
     \frac{|B_1|}{(1-\sigma)^{N+2}}
    \Big(\frac{\varrho^2}{\vartheta}\Big)^\frac{N}{2}
    \biint_{Q_o}|u|^{r}\,\dx\dt
    \bigg]^\frac{2}{\boldsymbol\lm_{r}}
\end{equation}
and 
\begin{equation}\label{k-2-3}
     \frac{\|F\|_{2p, Q_o}^{2p}}{k^{p(1+m)}}
     \le
    \frac{\vartheta^\frac{N+2}2|B_1|}{\vartheta ^p k^{\frac{N}{2}(1-m)}}\,\,\Leftrightarrow\,\,k\ge \bigg[
    \Big(\frac{\vartheta}{\varrho^2}\Big)^{p-\frac{N}2}
    \biint_{Q_o}|\varrho F|^{2p}\,\dx\dt\bigg]^\frac{2}{\boldsymbol \lambda_{p(1+m)}}
\end{equation}
hold true. 
This allows us to write the recursive inequality in a more compact form. Namely,
\begin{align*}
    \boldsymbol Y_{n+1}
    &\le
    \underbrace{C\|u\|_{\infty ,Q_o}^{r-mq}
    k^{mq-r+\frac{N}{2}(1-m)\alpha}
    \bigg[\frac{|Q_o|}{\vartheta^\frac{N+2}2(1-\sigma)^{N+2}}\bigg]^{\alpha}}_{=:\mathfrak C}
    \boldsymbol b^{n}\boldsymbol Y_n^{1+\alpha}.
\end{align*}
where
\begin{equation*}
    \alpha:= \frac1{N}\Big(1-\frac{N+2}{2p}\Big).
\end{equation*}
To apply \cite[\S\,2, Lemma 5.1]{DBGV-book} with $(\boldsymbol Y_n, \boldsymbol b, \mathfrak C,\alpha)$ we have to require that $\boldsymbol Y_o$ satisfies
\begin{align}\label{k-2-4}
    \boldsymbol Y_o&=\frac{1}{k^r}\biint_{Q_o}|u|^r\,\dx\dt
    \le
    \mathfrak C^{-\frac{1}{\alpha}}\boldsymbol b^{-\frac{1}{\alpha^2}} \nonumber\\
    &\equiv
    C^{-1} \|u\|_{\infty ,Q_o}^{-\frac{r-mq}{\alpha}}
    \frac{\vartheta^\frac{N+2}2(1-\sigma)^{N+2}}{|Q_o|}
     k^{\frac{r-mq}{\alpha}+\frac{N}{2}(m-1)}.
\end{align}
However, recalling that $r+\frac{N}{2}(m-1)=\frac12\boldsymbol \lm_r$, this is equivalent to a condition for $k$, i.e.
\begin{equation*}
     k^{\frac12\boldsymbol\lm_r+\frac{r-mq}{\alpha} }
    \ge 
    \frac{C\|u\|_{\infty ,Q_o}^{\frac{r-mq}{\alpha}}}{(1-\sigma)^{N+2}}
    \Big(\frac{\varrho^2}{\vartheta}\Big)^\frac{N}{2}\biint_{Q_o}|u|^r\,\dx\dt.
\end{equation*}
The implication of \cite[\S\,2, Lemma 5.1]{DBGV-book} on fast geometric convergence is that $\boldsymbol Y_n\to 0$ as $n\to \infty$, so that $\boldsymbol Y_\infty =0$, leading to
\begin{equation*}
      \sup_{Q_{\sigma\rho,\sigma\vartheta}}  |u|
    \le k.
\end{equation*}
Now we choose $k$ as
\begin{align*}
    k
    &:=
    \max\Bigg\{
    \|u\|_{\infty ,Q_o}^{1-\frac{\frac12\boldsymbol\lm_{r}}{\frac12\boldsymbol\lm_{r}
    +\frac{r-mq}{\alpha}}}
    \bigg[
    \frac{C}{(1-\sigma)^{N+2}}
    \Big(\frac{\varrho^2}{\vartheta}\Big)^\frac{N}{2}\biint_{Q_o}|u|^r\,\dx\dt
    \bigg]^\frac{1}{\frac12\boldsymbol\lm_{r}
    +\frac{r-mq}{\alpha}},\\
    &\qquad\qquad\;\,
    \bigg[
     \frac{C}{(1-\sigma)^{N+2}}
    \Big(\frac{\varrho^2}{\vartheta}\Big)^\frac{N}{2}
    \biint_{Q_o}|u|^{r}\,\dx\dt
    \bigg]^\frac{2}{\boldsymbol\lm_{r}},\\
    &\qquad\qquad\;\,
    \bigg[|B_1|
    \Big(\frac{\vartheta}{\varrho^2}\Big)^{p-\frac{N}2}
    \biint_{Q_o}|\varrho F|^{2p}\,\dx\dt\bigg]^\frac{2}{\boldsymbol \lambda_{p(1+m)}},
    \Big(
    \frac{\vartheta}{\rho^2}
    \Big)^{\frac{1}{1-m}}
    \Bigg\},
\end{align*}
so that all the conditions, i.e. \eqref{k-2-1}, \eqref{k-2-2}, \eqref{k-2-3}, \eqref{k-2-4}, imposed on $k$ so far are fulfilled.

The final quantitative $L^\infty$-estimate is obtained similarly as in the first case $r<1+3 m$. With the notation introduced there -- in particular \eqref{interpol-M-n} -- we obtain
\begin{align*}
    \boldsymbol M_n
    &\le 
    \max\Bigg\{
    \boldsymbol M_{n+1}^{1-\frac{\frac12\boldsymbol\lm_{r}}{\frac12\boldsymbol\lm_{r}
    +\frac{r-mq}{\alpha}}}
    \bigg[
     \frac{C2^{2Nn}}{(1-\tau)^{N+2}}
    \Big(\frac{\varrho^2}{\vartheta}\Big)^\frac{N}{2}
    \biint_{Q_{\varrho,\vartheta}}|u|^{r}\,\dx\dt
    \bigg]^\frac{1}{\frac12\boldsymbol\lm_{r}
    +\frac{r-mq}{\alpha}},\\
    &\qquad\qquad\;\,
    \bigg[
     \frac{C2^{2Nn}}{(1-\tau)^{N+2}}
    \Big(\frac{\varrho^2}{\vartheta}\Big)^\frac{N}{2}
    \biint_{Q_{\varrho,\vartheta}}|u|^{r}\,\dx\dt
    \bigg]^\frac{2}{\boldsymbol\lm_{r}},\\
    &\qquad\qquad\;\,\bigg[ C(N,p)
    \Big(\frac{\vartheta}{\varrho^2}\Big)^{p-\frac{N}2}
    \biint_{Q_{\varrho,\vartheta}}|\varrho F|^{2p}\,\dx\dt\bigg]^\frac{2}{\boldsymbol \lambda_{p(1+m)}},
    2^\frac{1}{1-m}\Big(
    \frac{\vartheta}{\rho^2}
    \Big)^{\frac{1}{1-m}}
    \Bigg\}.
\end{align*}
If for some $n\in \N_0$ the second, third, or fourth term in the maximum dominates the first, we immediately conclude the quantitative $L^\infty$-estimate. 
Otherwise, the first term is dominant for all $n\in\N_0$, so that for any $n\in\N_0$ we have
\begin{equation*}
    \boldsymbol  M_{n}
    \le
    \underbrace{
    \bigg[
     \frac{C}{(1-\tau)^{N+2}}
    \Big(\frac{\varrho^2}{\vartheta}\Big)^\frac{N}{2}
    \biint_{Q_{\varrho,\vartheta}}|u|^{r}\,\dx\dt
    \bigg]^\frac{1}{\frac12\boldsymbol\lm_{r}
    +\frac{r-mq}{\alpha}}}_{=:\mathfrak C}
    2^\frac{2Nn}{\frac12\boldsymbol\lm_{r}
    +\frac{r-mq}{\alpha}}
    \boldsymbol M_{n+1}^{1-\frac{\frac12\boldsymbol\lm_{r}}{\frac12\boldsymbol\lm_{r}
    +\frac{r-mq}{\alpha}}}.
\end{equation*}
Apply \cite[Chapter I, Lemma 4.3]{DiBe} with 
\begin{equation*}
    \alpha:= \frac{\frac12\boldsymbol\lm_{r}}{\frac12\boldsymbol\lm_{r}
    +\frac{r-mq}{\alpha}}
    \quad
    \mbox{and}
    \quad
    \boldsymbol b
    := 2^\frac{2N}{\frac12\boldsymbol\lm_{r}
    +\frac{r-mq}{\alpha}},
\end{equation*}
to conclude that
\begin{align*}
    \boldsymbol M_o
    &\le
    \bigg(\frac{2\mathfrak C}{\boldsymbol b^{1-\frac1{\alpha}}}
    \bigg)^\frac1{\alpha}
    =
    \bigg[
     \frac{C}{(1-\tau)^{N+2}}
    \Big(\frac{\varrho^2}{\vartheta}\Big)^\frac{N}{2}
    \biint_{Q_{\varrho,\vartheta}}|u|^{r}\,\dx\dt
    \bigg]^\frac{2}{\boldsymbol\lm_{r}}.
\end{align*}
Joining this with the first case,
we end up with the quantitative $L^\infty$-estimate for $r\ge1+3 m$. 
This finishes the proof.
\end{proof}

As a corollary we state the  sup-estimate from Proposition~\ref{prop:sup-est} on cylinders of the form $Q_{\rho}^{(\theta)}(z_o)$ as defined in \eqref{def-Q}. Here is a transition from backward cylinders to centered cylinders.

\begin{corollary}\label{cor:sup-est}
Let $m\in(0,m_c]$, $p>\frac{N+2}2$, and $r\ge 1$ satisfying \eqref{Eq:B:4:1}. Then there exists a constant
$C=C(N,m,\nu,L,p,r)$ such that, whenever $u$ is a locally bounded, local, weak solution to \eqref{por-med-eq} in $\Omega_T$ in the sense of Definition~{\upshape\ref{def:weak_solution}}  with $F\in L^{2p}_{\loc}(\Omega_T,\R^{kN})$, then for all cylinders 
$Q_{2\rho}^{(\theta)}(z_o)\Subset\Omega_T$, we have
\begin{align*}
    \sup_{Q_{\rho}^{(\theta)}(z_o)}|u| 
    &\le
    C\,
    \bigg[\big(\rho^{\frac{1}{m}}\theta^{\frac{2m}{1+m}}\big)^{\frac{N(m-1)}{2}}
    \biint_{Q_{2\rho}^{(\theta)}(z_o)} 
    |u|^r \,\dx\d\tau\bigg]^{\frac{2}{\boldsymbol{\lambda}_r}}\\
    &\phantom{\le\,}+
    C\,\rho^{\frac{1}{m}}
    \bigg[\theta^{\frac{mN(m-1)}{1+m}}
    \biint_{Q_{2\rho}^{(\theta)}(z_o)} |F|^{2p} \,\dx\d\tau\bigg]^{\frac{2}{\boldsymbol{\lambda}_{p(1+m)}}}
    +
    C\, \rho^{\frac{1}{m}}\theta^{\frac{2m}{1+m}}.
\end{align*}
\end{corollary}

\begin{proof}
We apply the sup-estimate from Proposition~\ref{prop:sup-est} with $\sigma =\frac12$ on  $Q_{2\rho}^{(\theta)}(z_o)$
instead of $Q_{\rho,\vartheta}(z_o)$. Then
\begin{equation*}
   \frac{\big(2\theta^{\frac{m(m-1)}{1+m}}\rho\big)^2}{(2\rho)^{\frac{1+m}{m}}}
    =2^{\frac{m-1}{m}}
    \big(\rho^{\frac{1}{m}}\theta^{\frac{2m}{1+m}}\big)^{m-1},
\end{equation*}
and $Q_{\rho}^{(\theta)}(z_o)\subset \tfrac12 Q_{2\rho}^{(\theta)}(z_o)$. The application yields
\begin{align*}
    \sup_{Q_{\rho}^{(\theta)}(z_o)}|u| 
    &\le
    C\, 
    \bigg[\big(\rho^{\frac{1}{m}}\theta^{\frac{2m}{1+m}}\big)^{\frac{N(m-1)}{2}}
    \biint_{Q_{2\rho}^{(\theta)}(z_o)} 
    |u|^r \,\dx\d\tau\bigg]^{\frac{2}{\boldsymbol{\lambda}_r}}\\
    &\phantom{\le\,}+
    C\,\bigg[\big(\rho^{\frac{1}{m}}\theta^{\frac{2m}{1+m}}\big)^{\frac{N}{2}(m-1)}\big(\rho^{\frac{1+m}{m}}\big)^{p}
    \biint_{Q_{2\rho}^{(\theta)}(z_o)} |F|^{2p} \,\dx\d\tau\bigg]^{\frac{2}{\boldsymbol{\lambda}_{p(m+1)}}}\\
    &\phantom{\le\,}+
    C\,\rho^{\frac{1}{m}}\theta^{\frac{2m}{1+m}} .
\end{align*}
For the radii in the second term we use the definition
of $\boldsymbol{\lambda}_{p(1+m)}= N(m-1)+2p(1+m)$  and compute 
$$
    \big(\rho^{\frac{1}{m}}\theta^{\frac{2m}{1+m}}\big)^{\frac{N}{2}(m-1)}\big(\rho^{\frac{1+m}{m}}\big)^{p}
    =
    \rho^{\frac{\boldsymbol{\lambda}_{p(1+m)}}{2m}}
    \theta^{\frac{mN(m-1)}{1+m}},
$$
which finishes the proof.
\end{proof}

\begin{remark}
\upshape In Definition~\ref{def:weak_solution} we require 
\[
    u\in C^0 \big((0,T); L^{1+m}(\Omega,\R^k)\big)
	\quad\mbox{with}\quad 
	\power{u}{m}\in L^2\big(0,T;W^{1,2}(\Omega,\R^k)\big)
\]
By the Sobolev embeddings this implies
$|u|^m\in L^q_{\rm loc}\big(\Omega_T,\R^k\big)$ with $q=2\frac{Nm+m+1}{Nm}$. Note that $q>\frac{m+1}m$ if $\frac{(N-2)_+}{N+2}<m<1$.  
Therefore, for $m$ in such a range, the condition $\boldsymbol \lambda_r>0$ is fulfilled with $1\le r\le\frac{2N}{N+2}$, so that we only have to deal with the case $r>\frac{2N}{N+2}$. This allows to use the entire proof from the sub-critical case $0<m\le\frac{(N-2)_+}{N+2}$ also in the super-critical range $\frac{(N-2)_+}{N+2}<m<1$, which per se is not what we are directly interested in here. The result is the quantitative $L^\infty$-estimate as stated in Proposition~\ref{prop:sup-est}, with the only condition that we have to choose $1\le r\le\frac{2N}{N+2}$; apart from this single difference, the final result is the same, that is, we have
\begin{align*}
    \sup_{Q_{\sigma\rho,\sigma\vartheta}}  |u|
    &\le
    \max\Bigg\{
    \bigg[
     \frac{C}{(1-\sigma)^{N+2}}
    \Big(\frac{\varrho^2}{\vartheta}\Big)^\frac{N}{2}
    \biint_{Q_{\varrho,\vartheta}}|u|^{r}\,\dx\dt
    \bigg]^\frac{2}{\boldsymbol\lm_{r}},\\
    &\qquad\qquad\;
    \bigg[
    \Big(\frac{\varrho^2}{\vartheta}\Big)^{\frac{N}2-p}
    \biint_{Q_{\varrho,\vartheta}}|\varrho F|^{2p}\,\dx\dt\bigg]^\frac{2}{\boldsymbol \lambda_{p(1+m)}},
    \bigg(
    \frac{\vartheta}{\rho^2}
    \bigg)^{\frac{1}{1-m}}
    \Bigg\}.
\end{align*}
\end{remark}

\section{Energy bounds}\label{sec:energy}
In this section we state without proof an energy inequality and a gluing lemma. Both follow by standard arguments
from the weak form \eqref{weak-solution} of the differential equation by testing with suitable testing functions.  Later on, they will be used in the proof of Sobolev-Poincar\'e and reverse H\"older-type inequalities.
The proof of the energy estimate is along
the lines of \cite[Lemma~3.1]{BDKS-higher-int}, taking into account
\cite[Lemma~2.3\,(i)]{BDKS-higher-int} or
\cite[Lemma~3.4]{BDKS-doubly} and the different definition of scaled
cylinders. The latter means  that the radii $\rho$ and $r$ in
\cite[Lemma~3.1]{BDKS-higher-int} have to be replaced by
$\theta^{\frac{m(m-1)}{1+m}}\rho\,$ and $\theta^{\frac{m(m-1)}{1+m}}r$.

\begin{lemma}\label{lem:energy}
Let $m>0$ and $u$ be a weak solution to \eqref{por-med-eq} in $\Omega_T$ in the sense of Definition~{\upshape\ref{def:weak_solution}}.
Then, on any cylinder $Q_{\rho}^{(\theta)}(z_o)\subseteq\Omega_T$ with
$\rho, \theta>0$, for any $r\in[\rho/2,\rho)$ and any $a \in\R^N$, we have 
\begin{align*}
	& \sup_{t \in \Lambda_r (t_o)} 
 	\mint_{B_r^{(\theta)} (x_o)} 
	\frac{\big|\power{u}{\frac{1+m}{2}}(t) - \power{a}{\frac{1+m}{2}}\big|^2}	
	{r^{\frac{1+m}{m}}} \dx +
	\biint_{Q_r^{(\theta)}(z_o)} |D\power{u}{m}|^2 \dx\dt \\
	&\qquad\leq 
	C\,\biint_{Q_\rho^{(\theta)}(z_o)} 
	\bigg[ 
	\frac{\big|\power{u}{\frac{1+m}{2}}-\power{a}{\frac{1+m}{2}}\big|^2}
	{\rho^{\frac{1+m}{m}}-r^{\frac{1+m}{m}}} + 
	\frac{\big|\power{u}{m} - \power{a}{m}\big|^2}
	{\theta^{\frac{2m(m-1)}{1+m}}(\rho-r)^2} + 
	|F|^{2} \bigg]\dx\dt ,
\end{align*}
where $C=C(N,m,\nu,L)$.
\end{lemma}

The following lemma compares the slice-wise mean values of a given weak solution 
taken at different times. It is often called {\em gluing lemma}. 
Such an assertion is necessary and very useful, since Poincar\'e's and Sobolev's inequality can only be applied slice-wise. 
The proof  is exactly as in \cite[Lemma~3.2]{BDKS-higher-int}, taking into account the different definitions of scaled cylinders.

\begin{lemma}\label{lem:time-diff}
Let $m>0$ and $u$ be a weak solution to \eqref{por-med-eq} in $\Omega_T$ in the sense of Definition~{\upshape\ref{def:weak_solution}}.
Then, for any cylinder $Q_{\rho}^{(\theta)}(z_o)\subseteq\Omega_T$ with
$\rho,\theta>0$ there exists $\hat\rho\in [\frac{\rho}{2},\rho]$ such that for all $t_1,t_2\in\Lambda_\rho(t_o)$ we have 
\begin{align*}
	\big|\langle u\rangle_{x_o;\hat\rho}^{(\theta)}(t_2) - 
	\langle u\rangle_{x_o;\hat\rho}^{(\theta)}(t_1)\big| 
	&\le
	C\,\theta^{\frac{m(1-m)}{1+m}}\rho^{\frac{1}{m}}
	\biint_{Q_{\rho}^{(\theta)}(z_o)} 
	\big[|D\power{u}{m}| + |F|\big] \dx\dt ,
\end{align*}
for a constant $C=C(N,L)$. 
\end{lemma}

\section{Sobolev-Poincar\'e-type inequality}\label{sec:poin}

In this section we consider cylinders $Q_{2\varrho}^{(\theta)}(z_o) \subseteq \Omega_T$, where $\varrho, \theta>0$ satisfy a sub-intrinsic coupling, 
in the sense that for some constants $r,K \geq 1$, we have
\begin{equation}\label{sub-intrinsic-poincare}
  	\biint_{Q_{2\rho}^{(\theta)}(z_o)} 
  	\frac{|u|^r}{(2\rho)^\frac{r}{m}}\,\dx\dt
  	\le 
   	K^\frac{r}{1+m} \theta^\frac{2r m}{1+m}.
\end{equation}
Furthermore, we assume that either of the following holds true: 
\begin{equation}\label{super-intrinsic-poincare}
  	\left\{
    \begin{array}{l}
    \displaystyle \theta^\frac{2r m}{1+m}
  	\le 
  	K^\frac{r}{1+m} \biint_{Q_\rho^{(\theta)}(z_o)} 
  	\frac{|u|^{r}}{\rho^\frac{r}{m}}\,\dx\dt
  	\quad\mbox{or} \\[25pt]
  	\displaystyle \theta^{2 m}
  	\le 
	K \Bigg[ \biint_{Q_{\rho}^{(\theta)}(z_o)} 
	|D\power{u}{m}|^2\dx\dt + 
    \bigg[ \biint_{Q_{\rho}^{(\theta)}(z_o)}
 |F|^{2p}\dx\dt\bigg]^\frac1p\Bigg].
    \end{array}\right.
\end{equation}
When both \eqref{sub-intrinsic-poincare} and \eqref{super-intrinsic-poincare}$_1$ hold, we refer to this coupling as intrinsic. Let us also recall the parameter 
\[
\boldsymbol\lambda_r:=N(m-1)+2r
\]
which will be used frequently.

The main result -- a Sobolev-Poicar\'e-type inequality -- will be proved in Lemma~\ref{lem:poin}. Prior to that we introduce three preliminary results. The first one examines the quantitative $L^{\infty}$-estimate derived in \S\,\ref{S:4} over subintrinsic cylinders.

\begin{lemma}\label{lem:sup-int}
Let $m\in(0,m_c]$, $p>\frac{N+2}2$, $F\in L^{2p}_{\loc}(\Omega_T,\R^{kN})$, and $u$ a local weak solution to \eqref{por-med-eq} in
$\Omega_T$ in the sense of
Definition~{\upshape\ref{def:weak_solution}} that satisfies $u\in L^r_{\rm loc}(\Omega_T, \R^N)$ for some $r$ with $\boldsymbol \lambda_r>0$. Then, on any cylinder
$Q_{2\rho}^{(\theta)}(z_o)\subseteq\Omega_T$ satisfying
\eqref{sub-intrinsic-poincare}, %and~\eqref{super-intrinsic-poincare}$_1$, 
with $\rho,\theta>0$ and some $K\ge1$, we have
\begin{equation*}
    \sup_{Q_{\rho}^{(\theta)}(z_o)}|u| 
    \le
    C \rho^{\frac{1}{m}} \theta^{\frac{2mN(m-1)}{(1+m)\boldsymbol\lambda_{p(1+m)}}} \bigg[
    \theta^{2pm} +
    \biint_{Q_{2\rho}^{(\theta)}(z_o)} |F|^{2p} \,\dx\dt\bigg]^{\frac{2}{\boldsymbol\lambda_{p(1+m)}}} ,
\end{equation*}
with a constant $C=C(N,m,\nu,L,p,r,K)$.
\end{lemma}

\begin{proof}
Throughout the proof we omit the center $z_o$ in our notation.
From Corollary~\ref{cor:sup-est}, we obtain 
\begin{align*}
    \sup_{Q_{\rho}^{(\theta)}}|u| 
    &\le
    C\bigg[\big(\rho^{\frac{1}{m}}\theta^{\frac{2m}{1+m}}\big)^{\frac{N(m-1)}{2}}
    \biint_{Q_{2\rho}^{(\theta)}} 
    |u|^{r} \,\dx\dt\bigg]^{\frac{2}{\boldsymbol \lambda_{r}}}\\
    &\quad+
    C\rho^{\frac{1}{m}}
    \bigg[\theta^{\frac{mN(m-1)}{1+m}}
    \biint_{Q_{2\rho}^{(\theta)}} |F|^{2p} \,\dx\dt\bigg]^{\frac{2}{\boldsymbol \lambda_{p(1+m)}}} +
    C \rho^{\frac{1}{m}}\theta^{\frac{2m}{1+m}},
\end{align*}
where $C=C(N,m,\nu,L,p,r)$. 
In view of assumption~\eqref{sub-intrinsic-poincare}, we find that 
\begin{align*}
    \bigg[\big(\rho^{\frac{1}{m}}\theta^{\frac{2m}{1+m}}\big)^{\frac{N(m-1)}{2}}
    \biint_{Q_{2\rho}^{(\theta)}} 
    |u|^{r} \,\dx\dt\bigg]^{\frac{2}{\boldsymbol \lambda_{r}}} 
    &\le 
    \Big[\big(\rho^{\frac{1}{m}}\theta^{\frac{2m}{1+m}}\big)^{\frac{N(m-1)}{2}}
    K^\frac{r}{1+m} (2\rho)^{\frac{r}m}\theta^\frac{2r m}{1+m}\Big]^{\frac{2}{\boldsymbol \lambda_{r}}} \\
    &=
    C\, \rho^{\frac{1}{m}}\theta^{\frac{2m}{1+m}},
\end{align*}
so that 
\begin{align*}
    \sup_{Q_{\rho}^{(\theta)}}|u| 
    &\le
    C \rho^{\frac{1}{m}}\theta^{\frac{2m}{1+m}} +
    C\rho^{\frac{1}{m}}
    \bigg[\theta^{\frac{mN(m-1)}{1+m}}
    \biint_{Q_{2\rho}^{(\theta)}} |F|^{2p} \,\dx\dt\bigg]^{\frac{2}{\boldsymbol \lambda_{p(1+m)}}} ,
\end{align*}
where $C=C(N,m,\nu,L,p,r,K)$. Re-writing this inequality yields the claim of the lemma.
\end{proof}

With the previous result at hand, we can give \eqref{super-intrinsic-poincare}$_1$ a new look.
\begin{lemma}\label{lem:alpha-super-int}
Let $m\in(0,m_c]$, $p>\frac{N+2}2$, $F\in L^{2p}_{\loc}(\Omega_T,\R^{kN})$, and $u$ a local
weak solution to \eqref{por-med-eq} in
$\Omega_T$ in the sense of
Definition~{\upshape\ref{def:weak_solution}} that satisfies $u\in L^r_{\rm loc}(\Omega_T, \R^N)$ for some $r$ with $\boldsymbol \lambda_r>0$. Then, on any cylinder
$Q_{2\rho}^{(\theta)}(z_o)\subseteq\Omega_T$ satisfying
\eqref{sub-intrinsic-poincare} and~\eqref{super-intrinsic-poincare}$_1$, with $\rho,\theta>0$ and some $K\ge1$, we have
\begin{equation*}
  	\theta^{2 m}
  	\le 
  	C\, \biint_{Q_\rho^{(\theta)}(z_o)} 
  	\frac{|u|^{1+m}}{\rho^{\frac{1+m}m}}\,\dx\dt +
    \bigg[\biint_{Q_{2\rho}^{(\theta)}(z_o)} |F|^{2p} \,\dx\dt\bigg]^{\frac{1}{p}}
\end{equation*}
with a constant $C=C(N,m,\nu,L,p,r,K)$.
\end{lemma}

\begin{proof}
Throughout the proof we omit the center $z_o$ in our notation.
In \eqref{super-intrinsic-poincare}$_1$ we use Lemma~\ref{lem:sup-int} in order to reduce the power of $|u|$ and obtain 
\begin{align*}
    \theta^\frac{2r m}{1+m}
  	&\le
  	K^\frac{r}{1+m}
    \bigg[\sup_{Q_\rho^{(\theta)}} 
    \frac{|u|}{\rho^{\frac{1}m}}\bigg]^{r-(1+m)}
    \biint_{Q_\rho^{(\theta)}} 
    \frac{|u|^{1+m}}{\rho^{\frac{1+m}m}}\,\dx\dt \\
    &\le 
    C\, \theta^{\frac{2mN(m-1)(r-(1+m))}{(1+m)\boldsymbol \lambda_{p(1+m)}}}\\
    &\phantom{\le C\,}\cdot\bigg[
    \theta^{2pm} +
    \biint_{Q_{2\rho}^{(\theta)}} |F|^{2p} \,\dx\dt\bigg]^{\frac{2(r-(1+m))}{\boldsymbol \lambda_{p(1+m)}}} 
    \biint_{Q_\rho^{(\theta)}} 
    \frac{|u|^{1+m}}{\rho^{\frac{1+m}m}}\,\dx\dt.
\end{align*}
We multiply both sides by
$ \theta^{\frac{2mN(1-m)(r-(1+m))}{(1+m)\boldsymbol \lambda_{p(1+m)}}}$ and compute %\textcolor{green}{($\lm_{sr}$ should be $\lm_{2pr}$.)} 
%\textcolor{red}{(I think it is correct, since by definition $\lm_{sr}$ contains a factor 2)}
\begin{align*}
    \frac{2r m}{1+m} + \frac{2mN(1-m)(r-(1+m))}{(1+m)\boldsymbol \lambda_{p(1+m)}}
    =
    2m\, \frac{\boldsymbol \lambda_{ pr}}{\boldsymbol\lambda_{p(1+m)}}.
\end{align*}
Hence, we obtain
\begin{align*}
    \theta^{2m \frac{\boldsymbol \lambda_{ pr}}{\boldsymbol\lambda_{p(1+m)}}} 
    \le
    C 
    \bigg[
    \theta^{2pm} +
    \biint_{Q_{2\rho}^{(\theta)}} |F|^{2p} \,\dx\dt\bigg]^{\frac{2(r-(1+m))}{\boldsymbol \lambda_{p(1+m)}}}
    \biint_{Q_\rho^{(\theta)}} 
    \frac{|u|^{1+m}}{\rho^{\frac{1+m}m}}\,\dx\dt  .
\end{align*}
Next, we raise both sides to the power $\frac{\boldsymbol \lambda_{p(1+m)}}{\boldsymbol\lambda_{pr}}$ and apply Young's inequality with the result
\begin{align*}
    \theta^{2m}
  	&\le 
    C 
    \Bigg[
    \theta^{2m} +
    \bigg[\biint_{Q_{2\rho}^{(\theta)}} |F|^{2p} \,\dx\dt\bigg]^{\frac{1}{p}}
    \Bigg]^{\frac{(r-(1+m))2p}{\boldsymbol\lambda_{pr}}}
    \bigg[\biint_{Q_\rho^{(\theta)}} 
    \frac{|u|^{1+m}}{\rho^{\frac{1+m}m}}\,\dx\dt 
    \bigg]^{\frac{\boldsymbol  \lambda_{p(1+m)}}{\boldsymbol\lambda_{pr}}} \\
    & \le
    \frac12 
    \Bigg[
    \theta^{2m} +
    \bigg[\biint_{Q_{2\rho}^{(\theta)}} |F|^{2p} \,\dx\dt\bigg]^{\frac{1}{p}}
    \Bigg] +
    C \biint_{Q_\rho^{(\theta)}} 
    \frac{|u|^{1+m}}{\rho^{\frac{1+m}m}}\,\dx\dt.
\end{align*}
Re-absorbing $\frac12\theta^{2m}$ into the left-hand side yields the claimed inequality.
\end{proof}

In case the considered cylinder is an intrinsic one, then the first term on the right-hand side of the energy inequality from Lemma~\ref{lem:energy} dominates the second one. This observation leads to the following improved energy inequality on intrinsic cylinders.

\begin{lemma}\label{lem:energy-intr}
Let $m\in (0,m_c]$, $F\in L^{2p}_{\loc}(\Omega_T,\R^{kN})$ with $p>\frac{N+2}{2}$, $K\ge 1$, and  $r>0$ with $\boldsymbol\lambda_r>0$. Then there exists a constant $C=C(N,m,\nu,L,p,r,K)$ such that whenever
$u$ is a local weak solution to \eqref{por-med-eq} in
$\Omega_T$ in the sense of Definition~{\upshape\ref{def:weak_solution}} that satisfies $u\in L^r_{\rm loc}(\Omega_T,\R^N)$, and \eqref{sub-intrinsic-poincare} and \eqref{super-intrinsic-poincare} on $Q_{2\rho}^{(\theta)}(z_o)\subseteq\Omega_T$ with parameters $K$ and $r$, then for any $\rho\le r_1<r_2\le 2\rho$ we have 
\begin{align*}
	& \sup_{t \in \Lambda_{r_1} (t_o)} 
 	\mint_{B_{r_1}^{(\theta)} (x_o)} 
	\frac{\big|\power{u}{\frac{1+m}{2}}(t) - (\power{u}{\frac{1+m}{2}})_{z_o; r_1}^{(\theta)}\big|^2}	
	{{r_1}^{\frac{1+m}{m}}} \dx +
	\biint_{Q_{r_1}^{(\theta)}(z_o)} |D\power{u}{m}|^2 \dx\dt \\
	&\qquad\leq 
	C\biint_{Q_{r_2}^{(\theta)}(z_o)} 
	\frac{\big|\power{u}{\frac{1+m}{2}}-(\power{u}{\frac{1+m}{2}})_{z_o; r_2}^{(\theta)}\big|^2}
	{({r_2}-{r_1})^{\frac{1+m}{m}}} \,\dx\dt + 
	C\bigg[\biint_{Q_{2\rho}^{(\theta)}(z_o)} |F|^{2p} \,\dx\dt\bigg]^{\frac{1}{p}} .
\end{align*}
%where $C=C(m,\nu,L)$.
\end{lemma}

\begin{proof}
We omit the reference to the center $z_o$ in the notation and consider radii $s_1,s_2$ with $r_1\le s_1<s_2\le r_2$. 
Note that hypotheses \eqref{sub-intrinsic-poincare} and \eqref{super-intrinsic-poincare} imply that analogous coupling conditions
are
satisfied on $Q_{s_2}^{(\theta)}$ with constant  $2^{N+2+\frac2m}K$
instead of $K$. 
From the energy estimate in Lemma \ref{lem:energy} 
%applied with $a=\power{\big[(u^{\frac{1+m}{2}})_{s_1}^{(\theta)}\big]}{\frac{2}{1+m}}$, 
we obtain  with a constant $C=C(m,\nu,L)$ that
\begin{align*}
	&\sup_{t \in \Lambda_{s_1}}
	\mint_{B_{s_1}^{(\theta)}} 
	\frac{\big|\power{u}{\frac{1+m}{2}}(t) - 
	(\power{u}{\frac{1+m}{2}})_{s_1}^{(\theta)}\big|^2}{{s_1}^{\frac{1+m}{m}}} \,\dx +
	\biint_{Q_{s_1}^{(\theta)}} |D\power{u}{m}|^2 \,\dx\dt \nonumber\\
	&\quad\le
	C\,\biint_{Q_{s_2}^{(\theta)}}
	\frac{\big|\power{u}{\frac{1+m}{2}} - 
	(\power{u}{\frac{1+m}{2}})_{s_1}^{(\theta)}\big|^2}
	{{s_2}^{\frac{1+m}{m}}-{s_1}^{\frac{1+m}{m}}}  \,\dx\dt +
	C\,\biint_{Q_{s_2}^{(\theta)}} 
	\frac{\big|\power{u}{m} - 
	\power{\big[(u^{\frac{1+m}{2}})_{s_1}^{(\theta)}\big]}{\frac{2m}{1+m}}
	\big|^2}
	{\theta^{\frac{2m(m-1)}{1+m}}({s_2}-{s_1})^2} \,\dx\dt
	\nonumber\\
	&\quad\quad +
	C\, \biint_{Q_{s_2}^{(\theta)}} |F|^{2} \,\dx\dt \nonumber\\
	&\quad =:
	\mbox{I} + \mbox{II} + \mbox{III},
\end{align*}
where the meaning of $\mbox{I}$, $\mbox{II}$ and $\mbox{III}$ is clear in this context. We let
$$
    \Phi(s)
    :=
    \biint_{Q_{s}^{(\theta)}} 
	\frac{\big|\power{u}{\frac{1+m}{2}} - 
	(\power{u}{\frac{1+m}{2}})_{s}^{(\theta)}\big|^2}
	{{s}^{\frac{1+m}{m}}}  \,\dx\dt
    \qquad\mbox{for $s\in[\rho,2\rho]$}
$$
and
$\mathcal R_{s_1,s_2}:=\frac{s_2}{s_2-s_1}$.
To estimate the term $\mbox{I}$ we first observe that $(s_2-s_1)^{\frac{1+m}{m}} \le s_2^{\frac{1+m}{m}}-s_1^{\frac{1+m}{m}}$. Together with an application of Lemma~\ref{lem:alphalemma}, this implies
\begin{align*}
	\mbox{I}
	&\le
	C\,\mathcal R_{s_1,s_2}^{\frac{1+m}m}\,
	\Phi(s_2) ,
\end{align*}
again with a constant $C$ depending on $m,\nu ,L$ only.
We now turn our attention to the term $\mbox{II}$, which we re-write as
\begin{align*}%\label{termII}
	\mbox{II}
	=
	C\,\mathcal R_{s_1,s_2}^{2} \theta^{\frac{2m(1-m)}{1+m}}
	\biint_{Q_{s_2}^{(\theta)}} 
	\frac{\big|\power{u}{m} - 
	\power{\big[(u^{\frac{1+m}{2}})_{s_1}^{(\theta)}\big]}{\frac{2m}{1+m}}
	\big|^2}
	{s_2^2} \,\dx\dt .
\end{align*}
If \eqref{super-intrinsic-poincare}$_2$ is satisfied, we apply Lemmas~\ref{lem:a-b} and \ref{lem:alphalemma}, and Young's inequality to obtain 
\begin{align*}
	\mbox{II}
	&\le 
	C\mathcal R_{s_1,s_2}^{2} \theta^{\frac{2m(1-m)}{1+m}}
	\Bigg[\biint_{Q_{s_2}^{(\theta)}} 
	\frac{\big|\power{u}{\frac{1+m}{2}} - 
	(\power{u}{\frac{1+m}{2}})_{s_1}^{(\theta)}\big|^2}
	{s_2^{\frac{1+m}{m}}} \,\dx\dt\Bigg]^{\frac{2m}{1+m}} \\
	&\le 
	C\mathcal R_{s_1,s_2}^{2} \theta^{\frac{2m(1-m)}{1+m}}
	\Phi(s_2)^{\frac{2m}{1+m}} \\
	&\le
	\frac1{2K}\theta^{2m} +
	C\,K^{\frac{1-m}{2m}}\mathcal R_{s_1,s_2}^{\frac{1+m}{m}}
	\Phi(s_2) \\
	&\le
	\frac12
	\Bigg[\biint_{Q_{s_2}^{(\theta)}} 
	|D\power{u}{m}|^2 \,\dx\dt + 
    \bigg[\biint_{Q_{s_2}^{(\theta)}} 
	|F|^{2p} \,\dx\dt\bigg]^{\frac1p} \Bigg] +
	C\,\mathcal R_{s_1,s_2}^{\frac{1+m}{m}}
	\Phi(s_2)
\end{align*}
with $C=C(m,\nu,L,K)$. 

Otherwise, if~\eqref{super-intrinsic-poincare}$_1$ is in force, we apply Lemma~\ref{lem:alpha-super-int} and obtain
\begin{align*}
  	\theta^{2m}
  	&\le
    C\biint_{Q_\rho^{(\theta)}} 
  	\frac{|u|^{1+m}}{\rho^{\frac{1+m}m}}\,\dx\dt +
    \bigg[\biint_{Q_{2\rho}^{(\theta)}} |F|^{2p} \,\dx\dt\bigg]^{\frac{1}{p}}\\ 
	&\le
	C
  	\biint_{Q_\rho^{(\theta)}} 
  	\frac{\big|\power{u}{\frac{1+m}{2}}-
  	(\power{u}{\frac{1+m}{2}})_{s_1}^{(\theta)}\big|^{2}}
  	{\rho^{\frac{1+m}m}}\,\dx\dt +
  	\frac{C\big|(\power{u}{\frac{1+m}{2}})_{s_1}^{(\theta)}\big|^2}
  	{\rho^{\frac{1+m}m}}
   +
    \bigg[\biint_{Q_{2\rho}^{(\theta)}} |F|^{2p} \,\dx\dt\bigg]^{\frac{1}{p}},
\end{align*}
where we have $C=C(N,m,\nu,L,p,r,K)$ from Lemma~\ref{lem:alpha-super-int}. This leads to
\begin{align*}
	\mbox{II} 
	\le
	c\,\mathcal R_{s_1,s_2}^{2} [\mathrm{II}_1 + \mathrm{II}_2 + \mathrm{II}_3],
\end{align*}
where we have set
\begin{align*}
	\mathrm{II}_1
	&:=
	\Bigg[\biint_{Q_\rho^{(\theta)}}
  	\frac{\big|\power{u}{\frac{1+m}{2}}-
  	(\power{u}{\frac{1+m}{2}})_{s_1}^{(\theta)}\big|^{2}}
  	{\rho^{\frac{1+m}m}} \,\dx\dt\Bigg]^{\frac{1-m}{1+m}} \\
	&\qquad\qquad\qquad\cdot
 \biint_{Q_{s_2}^{(\theta)}} 
	\frac{\big|\power{u}{m} - 
	\power{\big[(u^{\frac{1+m}{2}})_{s_1}^{(\theta)}\big]}{\frac{2m}{1+m}}\big|^2}
	{s_2^2} \,\dx\dt,
\end{align*}

\begin{align*}
	\mathrm{II}_2
	&:=
	\frac{\big|(\power{u}{\frac{1+m}{2}})_{s_1}^{(\theta)}\big|^{\frac{2(1-m)}{1+m}}}
	{\rho^{\frac{1-m}{m}}} 
	\biint_{Q_{s_2}^{(\theta)}} 
	\frac{\big|\power{u}{m} - 
	\power{\big[(u^{\frac{1+m}{2}})_{s_1}^{(\theta)}\big]}{\frac{2m}{1+m}}\big|^2}
	{s_2^2} \,\dx\dt, 
\end{align*}
and 
\begin{align*}
	\mathrm{II}_3
	&:=
	\bigg[\biint_{Q_{2\rho}^{(\theta)}} |F|^{2p} \,\dx\dt\bigg]^{\frac{1}{p}\frac{1-m}{1+m}} 
	\biint_{Q_{s_2}^{(\theta)}} 
	\frac{\big|\power{u}{m} - 
	\power{\big[(u^{\frac{1+m}{2}})_{s_1}^{(\theta)}\big]}{\frac{2m}{1+m}}\big|^2}
	{s_2^2} \,\dx\dt .
\end{align*}
Regarding term $\mathrm{II}_1$, we enlarge the domain of integration in the first integral from $Q_\rho^{(\theta)}$ to $Q_{s_2}^{(\theta)}$ and apply in turn  Lemma~\ref{lem:a-b}, H\"older's inequality and Lemma~\ref{lem:alphalemma}, and obtain
\begin{align*}
	\mathrm{II}_1
	\le
	C \biint_{Q_{s_2}^{(\theta)}} 
  	\frac{\big|\power{u}{\frac{1+m}{2}}-
  	(\power{u}{\frac{1+m}{2}})_{s_1}^{(\theta)}\big|^{2}}
  	{s_2^{\frac{1+m}m}}\,\dx\dt
	\le
	C\, \Phi(s_2);
\end{align*}
in a similar way, we have
\begin{align*}
	\mathrm{II}_3
	&\le
	C \bigg[\biint_{Q_{2\rho}^{(\theta)}} |F|^{2p} \,\dx\dt\bigg]^{\frac{1}{p}\frac{1-m}{1+m}}
    \Phi(s_2)^{\frac{2m}{1+m}} \\
    &\le 
    C \Phi(s_2) + \bigg[\biint_{Q_{2\rho}^{(\theta)}} |F|^{2p} \,\dx\dt\bigg]^{\frac{1}{p}};
\end{align*}
finally, to term $\mathrm{II}_2$ we apply Lemma~\ref{lem:Acerbi-Fusco} with $\alpha =\frac{1+m}{2m}$ and Lemma~\ref{lem:alphalemma}, which yield
\begin{align*}
	\mathrm{II}_2
	\le
	C \biint_{Q_{s_2}^{(\theta)}} 
	\frac{\big|\power{u}{\frac{1+m}{2}} - 
	(\power{u}{\frac{1+m}{2}})_{s_1}^{(\theta)}\big|^2}
	{s_2^{\frac{1+m}{m}}} \dx\dt 
	\le
	C \Phi(s_2) .
\end{align*}
Independently whether \eqref{super-intrinsic-poincare}$_1$ or \eqref{super-intrinsic-poincare}$_2$ are in force, combining both cases, we have
\begin{align*}
	\mbox{II}
	&\le
	\frac12 
	\biint_{Q_{s_2}^{(\theta)}} 
	|D\power{u}{m}|^2 \,\dx\dt +
	C\mathcal R_{s_1,s_2}^{\frac{1+m}{m}}
	\Phi(s_2) +
    C \bigg[\biint_{Q_{2\rho}^{(\theta)}} |F|^{2p} \,\dx\dt\bigg]^{\frac{1}{p}},
\end{align*}
with a constant $C=C(N,m,\nu,L,p,r,K)$. 
Inserting the estimates for I and II above, we find that 
\begin{align*}
	\sup_{t \in \Lambda_{s_1}}&
	\mint_{B_{s_1}^{(\theta)}} 
	\frac{\big|\power{u}{\frac{1+m}{2}}(t) - 
	(\power{u}{\frac{1+m}{2}})_{s_1}^{(\theta)}\big|^2}{s_1^{\frac{1+m}{m}}} \dx +
    \biint_{Q_{s_1}^{(\theta)}} |D\power{u}{m}|^2 \dx\dt \nonumber\\
    &\le
	\frac12 
	\biint_{Q_{s_2}^{(\theta)}} 
	|D\power{u}{m}|^2 \,\dx\dt  + 
    C\mathcal R_{s_1,s_2}^{\frac{1+m}{m}}
	\Phi(s_2) + C\bigg[\biint_{Q_{2\rho}^{(\theta)}} 
	|F|^{2p}\,\dx\dt\bigg]^{\frac1p} .
\end{align*}
To this inequality we apply the Iteration Lemma \ref{lem:tech-classical} to  re-absorb the term containing $|D\power{u}{m}|^2$ from the right-hand side into the left. This leads to the energy inequality on intrinsic cylinders and finishes the proof.
\end{proof}

The main objective of this section is to establish the following Sobolev-Poincaré-type inequality. 

\begin{lemma}\label{lem:poin}
Let $m\in (0,m_c]$, $p>\frac{N+2}{2}$, $F\in L^{2p}_{\loc}(\Omega_T,\R^{kN})$, $K\ge 1$, and  $r>0$ with
$\boldsymbol\lambda_r>0$. There exists a constant $C=C(N,m,\nu,L,p,r,K)$ such that, whenever
$u$ is a local, weak solution to \eqref{por-med-eq} in
$\Omega_T$ in the sense of Definition~{\upshape\ref{def:weak_solution}} that satisfies $u\in L^r_{\rm loc}(\Omega_T,\R^N)$,
\eqref{sub-intrinsic-poincare} and \eqref{super-intrinsic-poincare} on $Q_{2\rho}^{(\theta)}(z_o)\subseteq\Omega_T$ with parameters $K$ and $r$, then for any $\rho\le r_1<r_2\le 2\rho$, we have
\begin{align*}
	&\biint_{Q_{r_1}^{(\theta)}(z_o)} 
	\frac{\big|\power{u}{\frac{1+m}{2}}-
	(\power{u}{\frac{1+m}{2}})_{z_o;r_1}^{(\theta)}\big|^2}
	{r_1^{\frac{1+m}{m}}}  
	\,\dx\dt \\
	&\quad\le 
  	C \Bigg[
  	\Big(\frac{r_2}{r_2-r_1}\Big)^{\frac{1+m}{m}}
    \bigg[\biint_{Q_{r_2}^{(\theta)}(z_o)} 
	|D\power{u}{m}|^{2q}  \,\dx\dt 
  	\bigg]^{\frac{1}{q}} +
	\bigg[\biint_{Q_{2\rho}^{(\theta)}(z_o)} 
	|F|^{2p}  \,\dx\dt\bigg]^{\frac{1}{p}} \Bigg],
\end{align*}
where the integrability exponent $q$ is given by
\begin{equation}\label{def:q}
    q
    := 
    \frac{r N }{r (N+2)-N(1-m)}
    =
    \frac{r N }{rN+\boldsymbol\lambda_r}
    \in
    \Big(\frac{N}{N+2},1\Big) .
\end{equation}
\end{lemma}

\begin{proof}
We omit the reference to the center $z_o$ in the notation and consider radii $s_1,s_2$ with $r_1\le s_1<s_2\le r_2$. By $\hat{s} \in[\frac{1}{2} s_1, s_1]$ we denote the radius introduced in Lemma~\ref{lem:time-diff}. We start our considerations by estimating
\begin{align*}
	\biint_{Q_{s_1}^{(\theta)}} \!
	\frac{\big|\power{u}{\frac{1+m}{2}}-
	(\power{u}{\frac{1+m}{2}})_{s_1}^{(\theta)}\big|^2}
	{s_1^{\frac{1+m}{m}}}  
	\,\dx\dt
    &\le
	\biint_{Q_{s_1}^{(\theta)}} \!
	\frac{\big|\power{u}{\frac{1+m}{2}}-
	\power{\big[(u)_{\hat s}^{(\theta)}\big]}{\frac{1+m}{2}}\big|^2}
	{\rho^{\frac{1+m}{m}}}  
	\,\dx\dt 
    \le 
    2[\mathrm{I}+\mathrm{II}].
\end{align*}
Here we have abbreviated
\begin{align*}
	\mathrm{I}
    &:=
    \biint_{Q_{s_1}^{(\theta)}}
	\frac{\big|\power{u}{\frac{1+m}{2}}-
	\power{\big[\langle u\rangle_{\hat s}^{(\theta)}(t)\big]}{\frac{1+m}{2}}
	\big|^2}
    {\rho^{\frac{1+m}{m}}}\,\dx\dt,\\[7pt]
    \mathrm{II}
    &:=
    \bint_{\Lambda_{s_1}}
    \frac{\big|\power{\big[\langle u\rangle_{\hat s}^{(\theta)}(t)\big]}{\frac{1+m}{2}}-
    \power{\big[(u)_{\hat s}^{(\theta)}\big]}{\frac{1+m}{2}}\big|^2}
    {\rho^{\frac{1+m}{m}}}\,\dt.
\end{align*}
The first term can be estimated using Young's inequality with exponents $\frac{N+2}{N}$ and $\frac{N+2}2$, and Lemma~\ref{lem:alphalemma} with $p=2$, $\alpha =\frac{m+1}2$, $A=B_{\hat s}^{(\theta)}$, $B= B_{s_1}^{(\theta)}$, $a= \power{\big[\langle u^m\rangle_{s_1}^{(\theta)}(t)\big]}{\frac{1}{m}}$ and $a=[(\power{u}{\frac{1+m}{2}})_{s_1}^{(\theta)}]^{\frac{2}{1+m}}$. We obtain
\begin{align}\label{Est-I}
    \mathrm{I}
    &\le
    \frac{1}{\rho^{\frac{1+m}m}}
    \sup_{t\in\Lambda_{s_1}}
    \bigg[\bint_{B_{s_1}^{(\theta)}} 
    \big|\power{u}{\frac{1+m}{2}}-
	\power{\big[\langle u\rangle_{\hat s}^{(\theta)}(t)\big]}{\frac{1+m}{2}}
	\big|^2 \,\dx \bigg]^{\frac{2}{N+2}} \nonumber\\
	&\qquad\cdot
	\bint_{\Lambda_{s_1}}
    \bigg[\bint_{B_{s_1}^{(\theta)}} 
    \big|\power{u}{\frac{1+m}{2}}-
	\power{\big[\langle u\rangle_{\hat s}^{(\theta)}(t)\big]}{\frac{1+m}{2}}
	\big|^2 \,\dx
    \bigg]^{\frac{N}{N+2}} \dt \nonumber\\
    &\le
    \eps \sup_{t\in\Lambda_{s_1}}
    \bint_{B_{s_1}^{(\theta)}}
    \frac{\big|\power{u}{\frac{1+m}{2}}-
	(\power{u}{\frac{1+m}{2}})_{s_1}^{(\theta)}\big|^2}
	{\rho^{\frac{1+m}m}} \,\dx
    +
    \frac{C}{\eps^{\frac2N}\rho^{\frac{1+m}m}} \,
    \tilde{\mathrm{I}}^{\frac{N+2}{N}},
\end{align}
where $C=C(N, m)$ and
$$
	\tilde{\mathrm{I}}
	:=
	\bint_{\Lambda_{s_1}}
    \bigg[\bint_{B_{s_1}^{(\theta)}} 
    \big|\power{u}{\frac{1+m}{2}}-
	\power{\big[\langle u^m\rangle_{s_1}^{(\theta)}(t)\big]}{\frac{1+m}{2m}}
	\big|^2 \,\dx
    \bigg]^{\frac{N}{N+2}} \dt.
$$
We estimate the integral $\tilde{\mathrm{I}}$  by means of Lemma~\ref{lem:Acerbi-Fusco} applied with $\frac{1+m}{2 m}$ instead of $\alpha$ and Hölder's inequality in space with exponents 
$$
    \frac{r}{1-m}
    \quad\mbox{and}\quad 
    \beta:=\frac{r}{r-(1-m)}.
$$
This yields
\begin{align*}
	\tilde{\mathrm{I}}
    &\le
	C
    \bint_{\Lambda_{s_1}}
    \bigg[\bint_{B_{s_1}^{(\theta)}} 
    \big(|\power{u}{m}|+|\langle \power{u}{m}\rangle_{s_1}^{(\theta)}(t)|\big)^{\frac{1-m}{m}}
    \big|\power{u}{m}-\langle\power{u}{m}\rangle_{s_1}^{(\theta)}(t)\big|^2 \,\dx
    \bigg]^{\frac{N}{N+2}} \dt \\
     &\le
    C\bint_{\Lambda_{s_1}}
    \bigg[\bint_{B_{s_1}^{(\theta)}} 
    |u|^{r} 
    \,\dx\bigg]^{\frac{N(1-m)}{(N+2)r} } \mathrm S(t)^\frac{2N}{N+2}\,\dt
\end{align*}
for a constant $C=C(m)$. For the last line we introduced the short-hand notation
\begin{align*}
    \mathrm S(t)
    :=
    \bigg[\bint_{B_{s_1}^{(\theta)}}
    \big|\power{u}{m}-\langle\power{u}{m}\rangle_{s_1}^{(\theta)}(t)\big|^{2\beta} \,\dx\bigg]^\frac{1}{2\beta}.
\end{align*}
Recalling the definition of $q$, we have 
$$
    q=
    \frac{\beta N } {N+2\beta}
    =
    \frac{ Nr}{Nr +\boldsymbol\lambda_r}.
$$
Note that $2q>1$ is the Sobolev exponent of $2\beta$. 
Applying to the right-hand side integral of $\tilde{\mathrm{I}}$ the Hölder inequality in time with exponents
$$
    \frac{(N+2)r}{N(1-m)} 
    \quad \mbox{and}\quad 
    \frac{(N+2)r}{(N+2)r-N(1-m)}
    =
    \frac{N+2}{N}\cdot q,
$$
enlarging in the first integral the domain of integration from $ Q_{s_1}^{(\theta)}$ to $Q_{2\rho}^{(\theta)}$ and the sub-intrinsic coupling \eqref{sub-intrinsic-poincare}, we obtain
\begin{align*}
    \tilde{\mathrm{I}}
    &\le
    C 2^N \bigg[
    \biint_{Q_{2\rho}^{(\theta)}}
    |u|^{r} \,\dx\dt \bigg]^{\frac{N(1-m)}{(N+2)r}} 
    \bigg[\bint_{\Lambda_{s_1}}
    \mathrm S(t)^{2q}\,\dt
    \bigg]^{\frac{1}{2q}\frac{2N}{N+2}} \\
    &\le
    C 2^N \Big( K^\frac{r}{1+m}\theta^\frac{2r m}{1+m}\rho^\frac{r}{m}\Big)^\frac{N(1-m)}{(N+2)r} 
    \bigg[\bint_{\Lambda_{s_1}}
    \mathrm S(t)^{2q}\,\dt
    \bigg]^{\frac{1}{2q}\frac{2N}{N+2}} \\
    &=
    C \big(\theta^{\frac{m}{1+m}}\rho^{\frac{1}{2m}}\big)^{\frac{2N(1-m)}{N+2}}  
    \bigg[\bint_{\Lambda_{s_1}}
    \mathrm S(t)^{2q}\,\dt
    \bigg]^{\frac{1}{2q}\frac{2N}{N+2}}  
\end{align*}
where $C=C(N,m,K)$.
Using Sobolev's inequality on the time slices we get
\begin{align*}
    \mathrm S(t)
    &\le
    C \,\theta^\frac{m(m-1)}{1+m}\rho
    \bigg[\bint_{B_{s_1}^{(\theta)}}
    |D\power{u}{m}|^{2q} \,\dx\bigg]^\frac{1}{2q},
\end{align*}
for a constant $C=C(N,\beta)$.
Inserting this above, we get
\begin{align*}
	\tilde{\mathrm{I}}
     &\le
     C \big( \theta^{\frac{m}{1+m}}\rho^{\frac{1}{2m}}\big)^{\frac{2N(1-m)}{N+2}} 
     \big(\theta^\frac{m(m-1)}{1+m}\rho\big)^\frac{2N}{N+2}
    \bigg[\biint_{Q_{s_1}^{(\theta)}}
    |D\power{u}{m}|^{2q} \,\dx\dt\bigg]^{\frac{1}{2q}\frac{2N}{N+2}} \\
    &=
     C \rho^{\frac{1+m}{m}\frac{N}{N+2}}
     \bigg[\biint_{Q_{s_1}^{(\theta)}}
    |D\power{u}{m}|^{2q} \,\dx\dt\bigg]^{\frac{1}{2q}\frac{2N}{N+2}}.
\end{align*}
The condition on $r$, i.e.  $\boldsymbol\lambda_r>0$, also restricts the range of possible values of $\beta$ to the interval $\big[ 1, \frac{N}{N-2}\big]$ (recall that $N\ge 3$ in the sub-critical case). 
This allows us to eliminate the dependence on $\beta$ for the constant $C(N, \beta)$ from the Sobolev inequality. Consequently, there remain only the dependencies on $N$, $m$ and $K$ in the constant on the right-hand side of the last inequality. The combination of the last inequality with \eqref{Est-I} yields
\begin{align*}
	\mathrm{I}
    \le
    \eps \sup_{t\in\Lambda_{s_1}}
    \bint_{B_{s_1}^{(\theta)}}
    \frac{\big|\power{u}{\frac{1+m}{2}}(t)-
	(\power{u}{\frac{1+m}{2}})_{s_1}^{(\theta)}\big|^2}
	{\rho^{\frac{1+m}m}}\,\dx
  	+
  	\frac{C}{\epsilon^{\frac{2}{N}}}
  	\bigg[\biint_{Q_{s_1}^{(\theta)}} 
  	|D\power{u}{m}|^{2q} \,\dx\dt 
  	\bigg]^{\frac{1}{q}}.
\end{align*}
It remains to estimate II. To this end, we use the fact $\frac{1+m}{2} < 1$ in 
Lemma~\ref{lem:a-b} and the gluing Lemma~\ref{lem:time-diff} to deduce
\begin{align}\label{est-II-1}
     \mathrm{II}
     &\le
     C
     \bint_{\Lambda_{s_1}}
     \frac{\big|\langle u\rangle_{\hat s}^{(\theta)}(t)-(u)_{\hat s}^{(\theta)}\big|^{1+m}}
     {\rho^{\frac{1+m}{m}}}\,\dt \nonumber\\
     &\le
     \frac{C}{\rho^{\frac{1+m}m}}
     \bint_{\Lambda_{s_1}}\bint_{\Lambda_{s_1}}
     \big|\langle u\rangle_{\hat s}^{(\theta)}(t)-\langle u\rangle_{\hat s}^{(\theta)}(\tau)\big|^{1+m}
     \,\dt\d\tau \nonumber\\
     &\le
     C\,\theta^{m(1-m)}
     \bigg[\biint_{Q_{s_1}^{(\theta)}}\big[|D\power um|+|F|\big] \dx\dt
     \bigg]^{1+m},
\end{align}
for a constant $C=C(N,m,L)$. If \eqref{super-intrinsic-poincare}$_2$ is satisfied, we have
\begin{align*}
     \mathrm{II}
     &\le
     C\,K^\frac{1-m}{2}\bigg[\biint_{Q_{s_1}^{(\theta)}}\big[|D\power um|^2+|F|^2\big] \dx\dt
     \bigg]^{\frac{1-m}{2}}
     \bigg[\biint_{Q_{s_1}^{(\theta)}}\big[|D\power um|+|F|\big] \dx\dt
     \bigg]^{1+m} \\
     &\le
     C\,K^\frac{1-m}{2}\Bigg[
     \bigg[\biint_{Q_{s_1}^{(\theta)}}|D\power um|^2 \dx\dt
     \bigg]^{\frac{1-m}{2}}
     +
     \bigg[\biint_{Q_{s_1}^{(\theta)}}|F|^2 \dx\dt
     \bigg]^{\frac{1-m}{2}}\Bigg]
     \\
     &\qquad\qquad\qquad\qquad
     \cdot \Bigg[
     \bigg[\biint_{Q_{s_1}^{(\theta)}}|D\power um|^{2q}
     \dx\dt\bigg]^\frac1{q}
     +
     \biint_{Q_{s_1}^{(\theta)}}|F|^{2}
     \dx\dt
     \Bigg]^\frac{1+m}2\\
     &\le
     \epsilon\,\biint_{Q_{s_1}^{(\theta)}}|D\power um|^2 \dx\dt +
     \frac{C}{\epsilon^{\frac{1-m}{1+m}}}\Bigg[
     \bigg[\biint_{Q_{s_1}^{(\theta)}}|D\power um|^{2q} \dx\dt
     \bigg]^{\frac1q} +
     \biint_{Q_{s_1}^{(\theta)}} |F|^2 \dx\dt\Bigg],
\end{align*}
where the constant $C$ depends only on $N$, $m$, $L$, and $K$. To obtain the second line we used H\"older's inequality to increase the powers of
$|D\power um|$ and $|F|$, while for the last line we used Young's inequality with exponents $\frac{2}{1-m}$ and $\frac{2}{1+m}$. 
Note that $\frac{1-m}{1+m}\ge\frac{2}{N}$ since $m\le m_c$, which yields $\varepsilon^{-{\frac{1-m}{1+m}}}\ge\varepsilon^{-{\frac 2N}}$.

If \eqref{super-intrinsic-poincare}$_1$ is satisfied, we argue as follows. First, observe that by Lemma~\ref{lem:alpha-super-int} 
\begin{align*}
     \theta^{2m}
     &\le
     C\,\biint_{Q_\rho^{(\theta)}}
     \frac{\big|\power{u}{\frac{1+m}{2}}-
     \power{\big[(u)_{\hat s}^{(\theta)}\big]}{\frac{1+m}{2}}\big|^2}
     {\rho^{\frac{1+m}m}}\dx\dt
     +
     \frac{C\,\big|(u)_{\hat s}^{(\theta)}\big|^{1+m}}{\rho^{\frac{1+m}{m}}} \\
    &\quad +
    \bigg[\biint_{Q_{2\rho}^{(\theta)}} |F|^{2p} \,\dx\dt\bigg]^{\frac{1}{p}},
\end{align*}
where $C=C(N,m,\nu,L,p,r,K)$.
Therefore, we have
\begin{align*}
     \mathrm{II}
     &=
     \frac{(\theta^{2m})^\frac{1-m}{1+m}\, \mathrm{II}}{\theta^{\frac{2m(1-m)}{1+m}}}
     \le
     C\big[\mathrm{II}_1+\mathrm{II}_2+\mathrm{II}_3\big],
\end{align*}
with
\begin{align*}
     \mathrm{II}_1
     &:=
     \frac{1}{\theta^{\frac{2m(1-m)}{1+m}}}
     \Bigg[\biint_{Q_\rho^{(\theta)}}
     \frac{\big|\power{u}{\frac{1+m}{2}}-
     \power{\big[(u)_{\hat s}^{(\theta)}\big]}{\frac{1+m}{2}}\big|^2}
     {\rho^{\frac{1+m}m}}\,\dx\dt
     \Bigg]^{\frac{1-m}{1+m}}\cdot\mathrm{II}, \\[7pt]
     \mathrm{II}_2
     &:=
     \frac{1}{\theta^{\frac{2m(1-m)}{1+m}}}
     \frac{\big|(u)_{\hat s}^{(\theta)}\big|^{1-m}}{\rho^{\frac{1-m}{m}}}\cdot\mathrm{II} \\[7pt]
     \mathrm{II}_3
     &:=
     \frac{1}{\theta^{\frac{2m(1-m)}{1+m}}}
     \bigg[\biint_{Q_{2\rho}^{(\theta)}} |F|^{2p} \,\dx\dt\bigg]^{\frac{1-m}{p(1+m)}}
     \cdot\mathrm{II}.
\end{align*}
To estimate $\mathrm{II}_1$, we first note that enlarging the domain of integration from $Q_{\rho}^{(\theta)}$ to $Q_{2\rho}^{(\theta)}$, then applying Hölder's inequality to increase the  exponent $1+m$ by the factor 
$\frac{r}{1+m}$ (possible since $r>\frac{N(1-m)}2\ge 1+m$, because $0<m\le m_c$), and finally using \eqref{sub-intrinsic-poincare},
we can estimate the integral as 
\begin{align*}
    \biint_{Q_\rho^{(\theta)}}&
    \frac{\big|\power{u}{\frac{1+m}{2}}-
    \power{\big[(u)_{\rho}^{(\theta)}\big]}{\frac{1+m}{2}}\big|^2}
    {\rho^{\frac{1+m}m}}\,\dx\dt\\
    &\le
    C  \biint_{Q_{2\rho}^{(\theta)}}
    \frac{|u|^{1+m}}
    {(2\rho)^{\frac{1+m}m}}\,\dx\dt
    \le 
    C \bigg[ \biint_{Q_{2\rho}^{(\theta)}}
    \frac{|u|^{r}}{(2\rho)^{\frac{r}m}}\,\dx\dt\bigg]^\frac{1+m}{r}
    \le 
    C K\theta^{2m}
\end{align*}
for a constant $C=C(N,m)$.
Splitting the outer exponent of the integral $\mathrm{II}_1$ into $\frac{1-m}{1+m}=
\frac{(1-m)^2}{2(1+m)}+\frac{1-m}2$ creates a product with two factors. Estimating the first factor with the help of the preceding inequality yields 
\begin{align*}
	\mathrm{II}_1
    &\le
    \frac{C}{\theta^{m(1-m)}}
    \Bigg[\biint_{Q_\rho^{(\theta)}}
    \frac{\big|\power{u}{\frac{1+m}{2}}-
    \power{\big[(u)_{\hat s}^{(\theta)}\big]}{\frac{1+m}{2}}\big|^2}
    {\rho^{\frac{1+m}m}}\,\dx\dt
    \Bigg]^{\frac{1-m}{2}} 
    \cdot \mathrm{II}\\
     &\le
    \frac{C}{\theta^{m(1-m)}}
    \Bigg[\biint_{Q_\rho^{(\theta)}}
    \frac{\big|\power{u}{\frac{1+m}{2}}-
    (\power{u}{\frac{1+m}{2}})_{s_1}^{(\theta)}\big|^2}
    {\rho^{\frac{1+m}m}}\,\dx\dt
    \Bigg]^{\frac{1-m}{2}}\cdot \mathrm{II}.
\end{align*}
From the second-to-last line we used Lemma~\ref{lem:alphalemma} with $p=2$, $\al=\frac{m+1}{2}$,
$A=Q_{\hat s}^{(\theta)}$, $B=Q_{\rho}^{(\theta)}$, and $a=[(\power{u}{\frac{1+m}{2}})_{s_1}^{(\theta)}]^{\frac{2}{1+m}}$ in order to replace $\power{\big[(u)_{\hat s}^{(\theta)}\big]}{\frac{1+m}{2}}$ with $(\power{u}{\frac{1+m}{2}})_{s_1}^{(\theta)}$.
Finally, we apply \eqref{est-II-1} and Young's inequality with exponents $\frac{2}{1-m}$ and $\frac{2}{1+m}$ to obtain
\begin{align*}
    %C\,
    \mathrm{II}_1&\le 
    C\Bigg[\biint_{Q_\rho^{(\theta)}}
    \frac{\big|\power{u}{\frac{1+m}{2}}-
    (\power{u}{\frac{1+m}{2}})_{s_1}^{(\theta)}\big|^2}
    {\rho^{\frac{1+m}m}}\,\dx\dt
    \Bigg]^{\frac{1-m}{2}}
    \bigg[
    \biint_{Q_{s_1}^{(\theta)}}\big[|D\power um|+|F|\big]\dx\dt
    \bigg]^{1+m}\\
    &\le
    \frac1{4C}
    \biint_{Q_{s_1}^{(\theta)}}
    \frac{\big|\power{u}{\frac{1+m}{2}}-
    (\power{u}{\frac{1+m}{2}})_{s_1}^{(\theta)}\big|^2}
    {s_1^{\frac{1+m}{m}}}\,\dx\dt
    +
    C\bigg[\biint_{Q_{s_1}^{(\theta)}}\big[|D\power um|+|F|\big]\d x\d t
    \bigg]^2,
\end{align*}
with a constant $C=C(N,m, \nu,L,p,r, K)$.  
We point out that the constant $C$ on
the right-hand side is deliberately chosen to be the one from the decomposition of $\mathrm{II}$. 

For the term $\mathrm{II}_2$, we proceed as follows. We first insert the expression for the term II, then use Lemma~\ref{lem:Acerbi-Fusco}  with $\alpha:=\frac{2}{1+m}$, and finally apply the gluing Lemma~\ref{lem:time-diff}. This leads to
\begin{align*}
     \mathrm{II}_2
     &=
     \frac{1}{\theta^{\frac{2m(1-m)}{1+m}}\rho^{\frac2m}}\,
     \bint_{\Lambda_{s_1}}
     \big|(u)_{\hat s}^{(\theta)}\big|^{1-m}
     \Big|\power{\big[\langle u\rangle_{\hat s}^{(\theta)}(t)\big]}{\frac{1+m}{2}}-
     \power{\big[(u)_{\hat s}^{(\theta)}\big]}{\frac{1+m}{2}}\Big|^2
     \,\dt\\
     &\le
     \frac{C}{\theta^{\frac{2m(1-m)}{1+m}}\rho^{\frac2m}}\,
     \bint_{\Lambda_{s_1}}
     \big|\langle u\rangle_{\hat s}^{(\theta)}(t)-(u)_{\hat s}^{(\theta)}\big|^2
     \,\dt\\
     &\le
     \frac{C}{\theta^{\frac{2m(1-m)}{1+m}}\rho^{\frac2m}}\,
     \bint_{\Lambda_{s_1}}\bint_{\Lambda_{s_1}}
     \big|\langle u\rangle_{\hat s}^{(\theta)}(t)-\langle u\rangle_{\hat s}^{(\theta)}(\tau)\big|^2
     \,\dt\d\tau\\
     &\le
     C\bigg[
     \biint_{Q^{(\theta)}_{s_1}} \big[|D\power um|+|F|\big] \dx\dt\bigg]^2,
\end{align*}
with a constant $C=C(m, L)$. 

For the estimate of $\mathrm{II}_3$ we first apply Young's inequality with exponents $\frac{1+m}{1-m}$ and $\frac{1+m}{2m}$, and then use the identity $\frac{1+m}{2m}=\frac{(1-m)^2}{2m(1+m)}+\frac{2}{1+m}$ to re-write the result as
\begin{align*}
     \mathrm{II}_3
     &\le 
     \bigg[\biint_{Q_{2\rho}^{(\theta)}} |F|^{2p} \,\dx\dt\bigg]^{\frac{1}{p}} +
     \frac{1}{\theta^{1-m}} \mathrm{II}^{\frac{1+m}{2m}} \\
     &=
     \bigg[\biint_{Q_{2\rho}^{(\theta)}} |F|^{2p} \,\dx\dt\bigg]^{\frac{1}{p}} +
     \frac{1}{\theta^{1-m}}
     \mathrm{II}^{\frac{(1-m)^2}{2m(1+m)}}\mathrm{II}^{\frac{2}{1+m}} .
\end{align*}
Note that assumption~\eqref{sub-intrinsic-poincare} allows to estimate
\begin{align*}
    \mathrm{II}
    &\le
    C \biint_{Q_{2\rho}^{(\theta)}}
    \frac{|u|^{1+m}}
    {(2\rho)^{\frac{1+m}m}}\,\dx\dt
    \le 
    C \bigg[ \biint_{Q_{2\rho}^{(\theta)}}
    \frac{|u|^{r}}{(2\rho)^{\frac{r}m}}\,\dx\dt\bigg]^\frac{1+m}r
    \le 
    C K\theta^{2m}.
\end{align*}
In the estimate of $\mathrm{II}_3$ using the preceding inequality for the first term containing $\mathrm{II}$, and \eqref{est-II-1} for the second one yields
\begin{align*}
     \mathrm{II}_3
     \le & 
     \bigg[\biint_{Q_{2\rho}^{(\theta)}} |F|^{2p} \,\dx\dt\bigg]^{\frac{1}{p}}\\
     &+
     \frac{C}{\theta^{1-m}}
     \theta^{\frac{(1-m)^2}{1+m}}\Bigg[\theta^{m(1-m)}
     \bigg[\biint_{Q_{s_1}^{(\theta)}}\big[|D\power um|+|F|\big] \dx\dt
     \bigg]^{1+m}\Bigg]^{\frac{2}{1+m}} \\
     = &
     \bigg[\biint_{Q_{2\rho}^{(\theta)}} |F|^{2p} \,\dx\dt\bigg]^{\frac{1}{p}} +
     C \bigg[\biint_{Q_{s_1}^{(\theta)}}\big[|D\power um|+|F|\big] \dx\dt
     \bigg]^{2}.
\end{align*}
Now we collect the estimates for I and II; for the latter term we need to consider whether \eqref{super-intrinsic-poincare}$_2$ is satisfied or \eqref{super-intrinsic-poincare}$_1$ holds (in this second case, we rely on the upper bounds for II$_1$, II$_2$, and II$_3$); in any case, we arrive at
\begin{align*}
	\biint_{Q_{s_1}^{(\theta)}} &
    \frac{\big|\power{u}{\frac{1+m}{2}}-
    (\power{u}{\frac{1+m}{2}})_{s_1}^{(\theta)}\big|^2}
    {s_1^{\frac{1+m}{m}}} \,\dx\dt \\
    &\le
    \tfrac12 \biint_{Q_{s_1}^{(\theta)}}
	\frac{\big|\power{u}{\frac{1+m}{2}}-
    (\power{u}{\frac{1+m}{2}})_{s_1}^{(\theta)}\big|^2}
    {s_1^{\frac{1+m}{m}}} \,\dx\dt \\
    &\quad +
    \epsilon\Bigg[ 
    \sup_{t\in\Lambda_{s_1}} 
 	\bint_{B_{s_1}^{(\theta)}} \frac{\big|\power{u}{\frac{1+m}{2}}(t)-
    (\power{u}{\frac{1+m}{2}})_{s_1}^{(\theta)}\big|^2}
    {\rho^{\frac{1+m}{m}}} \,\dx +
    \biint_{Q_{s_1}^{(\theta)}} 
  	|D\power{u}{m}|^2 \,\dx\dt \Bigg] \\
  	&\quad +
  	\frac{C}{\epsilon^{\frac{1-m}{1+m}}}
  	\bigg[\biint_{Q_{s_1}^{(\theta)}} 
  	|D\power{u}{m}|^{2q} \dx \dt 
  	\bigg]^{\frac{1}{q}} +
    C\bigg[
    \biint_{Q^{(\theta)}_{s_1}} \big[|D\power um|+|F|\big] \dx\dt
	\bigg]^2 \\
    &\quad +
    C \bigg[\biint_{Q_{2\rho}^{(\theta)}} |F|^{2p} \,\dx\dt\bigg]^{\frac{1}{p}}.
\end{align*}
As already pointed out before, we used that $\varepsilon^{-{\frac{1-m}{1+m}}}\ge \varepsilon^{-{\frac 2N}}$, because  $m\le m_c$, and hence $\frac{1-m}{1+m}\ge\frac{2}{N}$.
Re-absorbing the first term on the right-hand side into the left-hand one, applying in turn Hölder's inequality 
to raise the power of $|D\power um|$ and $|F|$ and Lemma~\ref{lem:energy-intr}, and choosing $\varepsilon=\frac1{2C}\big(\frac{s_2-s_1}{s_2}\big)^{\frac{1+m}{m}}$ with $C$ from Lemma~\ref{lem:energy-intr}, we  obtain 
\begin{align*}
	\biint_{Q_{s_1}^{(\theta)}} &
    \frac{\big|\power{u}{\frac{1+m}{2}}-
    (\power{u}{\frac{1+m}{2}})_{s_1}^{(\theta)}\big|^2}
    {s_1^{\frac{1+m}{m}}} \,\dx\dt \\
    &\le
    \tfrac12 \biint_{Q_{s_2}^{(\theta)}} 
	\frac{\big|\power{u}{\frac{1+m}{2}}-
	(\power{u}{\frac{1+m}{2}})_{s_2}^{(\theta)}\big|^2}{s_2^{\frac{1+m}{m}}} \,\dx\dt \\ 
    &\qquad +
  	C\Big(\frac{s_2}{s_2-s_1}\Big)^{\frac{1-m}{m}} \bigg[\biint_{Q_{r_2}^{(\theta)} }
  	|D\power{u}{m}|^{2q} \dx \dt 
  	\bigg]^{\frac{1}{q}} +
    C \bigg[\biint_{Q_{2\rho}^{(\theta)}} |F|^{2p} \,\dx\dt\bigg]^{\frac{1}{p}}.
\end{align*}
Applying the Iteration Lemma \ref{lem:tech-classical}, we can re-absorb the first term on the right-hand side into the left-hand one, and conclude the asserted Sobolev-Poincaré inequality.
\end{proof}

\section{Reverse H\"older inequality}\label{sec:revholder}

The core of any proof of higher integrability of the gradient is a reverse H\"older inequality.
In this section, we establish such an inequality on intrinsic cylinders as introduced in Section~\ref{sec:poin}. 

\begin{proposition}\label{prop:revhoelder}
Let $m\in(0,m_c]$, $F\in L^{2p}_{\loc}(\Omega_T,\R^{kN})$ with 
$p>\frac{N+2}2$, $r>0$ with $\boldsymbol\lambda_r>0$, $K\ge 1$, and 
$u$ a  local, weak solution to \eqref{por-med-eq} in
$\Omega_T$ in the sense of Definition~{\upshape\ref{def:weak_solution}} that satisfies $u\in L^r_{\rm loc}(\Omega_T,\R^N)$. 
Then, on any  cylinder $Q_{2\rho}^{(\theta)}(z_o)\subseteq\Omega_T$ with $\rho,\theta>0$ satisfying \eqref{sub-intrinsic-poincare} and \eqref{super-intrinsic-poincare} with parameters $r$ and $K$, we have
\begin{align*}
	\biint_{Q_{\rho}^{(\theta)}(z_o)} |D\power{u}{m}|^2 \,\dx\dt 
	\le
	C\bigg[\biint_{Q_{2\rho}^{(\theta)}(z_o)} 
	|D\power{u}{m}|^{2q} \,\dx\dt \bigg]^{\frac{1}{q}} +
	C \bigg[\biint_{Q_{2\rho}^{(\theta)}(z_o)} |F|^{2p} \,\dx\dt \bigg]^{\frac{1}{p}},
\end{align*}
for a constant $C=C(N,m,\nu,L,p,r,K)$. Here, $0<q<1$ 
is the integrability exponent from \eqref{def:q}.
\end{proposition}

\begin{proof} 
We omit the reference to the center $z_o$ in the notation. 
In turn, we apply the energy estimate on intrinsic cylinders in Lemma \ref{lem:energy-intr} and the Poincaré-type inequality from Lemma~\ref{lem:poin} to obtain
\begin{align*}
	&
	\biint_{Q_{\rho}^{(\theta)}} |D\power{u}{m}|^2 \,\dx\dt \\
	&\quad\le
	C\,\biint_{Q_{\frac32\rho}^{(\theta)}}
	\frac{\big|\power{u}{\frac{1+m}{2}} - 
	(\power{u}{\frac{1+m}{2}})_{\frac32\rho}^{(\theta)}\big|^2}
	{\rho^{\frac{1+m}{m}}}  \,\dx\dt +
	C \bigg[\biint_{Q_{2\rho}^{(\theta)}} |F|^{2p} \,\dx\dt\bigg]^{\frac1p} \\
    &\quad\le 
  	C \Bigg[
  	\bigg[\biint_{Q_{2\rho}^{(\theta)}} 
	|D\power{u}{m}|^{2q}  \,\dx\dt 
  	\bigg]^{\frac{1}{q}} +
	\bigg[\biint_{Q_{2\rho}^{(\theta)}} 
	|F|^{2p}  \,\dx\dt\bigg]^{\frac{1}{p}} \Bigg],
\end{align*}
where $C=C(N,m,\nu,L,p,r,K)$. This proves the claimed reverse H\"older inequality.
\end{proof}

At the end of this section we provide a technical auxiliary result, which is essentially a direct consequence
of Lemma~\ref{lem:poin} and the energy inequality.

\begin{lemma}\label{lem:theta}
Let $m\in(0,m_c]$, $F\in L^{2p}_{\loc}(\Omega_T,\R^{kN})$ with $p>\frac{N+2}2$, $r>0$ with $\boldsymbol\lambda_r >0$, and 
$u$ a local weak solution to \eqref{por-med-eq} in
$\Omega_T$ in the sense of
Definition~{\upshape\ref{def:weak_solution}}
that satisfies $u\in L^r_{\rm loc}(\Omega_T,\R^n)$.
Then, on any  cylinder $Q_{2\rho}^{(\theta)}(z_o)\subseteq\Omega_T$ with $\rho,\theta>0$ satisfying \eqref{sub-intrinsic-poincare} and \eqref{super-intrinsic-poincare}$_1$ with $K=1$, and $r$, we have
\begin{align*}
  	\theta^m
  	&\le
    \tfrac1{\sqrt2}
    \Bigg[\biint_{Q_{\rho/2}^{(\theta)}(z_o)}
    \frac{|u|^{1+m}}{(\rho/2)^{\frac{1+m}{m}}} \,\dx\dt 
    \Bigg]^{\frac12} \\
    &\quad +
  	C\,\Bigg[\biint_{Q_{2\rho}^{(\theta)}(z_o)} 
	|D\power{u}{m}|^{2}  \,\dx\dt +
	\bigg[\biint_{Q_{2\rho}^{(\theta)}(z_o)} 
	|F|^{2p}  \,\dx\dt\bigg]^{\frac{1}{p}}\Bigg]^{\frac12} ,
\end{align*}
where $C=C(N,m,\nu,L,p,r)$.
\end{lemma}
Notice that since we assume $K=1$, in this case the constant $C$ does not depend explicitly on $K$.

\begin{proof}
We omit the reference to the center $z_o$ in the notation. 
We use \eqref{super-intrinsic-poincare}$_1$ with $K=1$ and Minkowski's inequality and Lemma~\ref{lem:alphalemma} to deduce 
\begin{align}\label{est-theta-2}
    \theta^{m}
    &\le
    \bigg[\biint_{Q_{\rho}^{(\theta)}} 
    \frac{|u|^{r}}{\rho^{\frac{r}m}} \,\dx\dt\bigg]^{\frac{1+m}{2r}} \nonumber\\
    &\le 
    \Bigg[\biint_{Q_{\rho}^{(\theta)}}  
    \frac{\big|\power{u}{\frac{1+m}{2}} -
    (\power{u}{\frac{1+m}{2}})_{\rho/2}^{(\theta)}\big|^\frac{2r}{1+m}}
    {\rho^{\frac{r}{m}}}
    \,\dx\dt
    \Bigg]^{\frac{1+m}{2r}} +
    \frac{\babs{(\power{u}{\frac{1+m}{2}})_{\rho/2}^{(\theta)}}}
    {\rho^{\frac{1+m}{2m}}} \nonumber\\
    &\le 
    \mathrm{A}^{\frac{1+m}{2r}} +
    \frac{\babs{(\power{u}{\frac{1+m}{2}})_{\rho/2}^{(\theta)}}}
    {\rho^{\frac{1+m}{2m}}},
\end{align}
where 
\begin{equation*}
    \mathrm{A}
    := 
    \biint_{Q_{\rho}^{(\theta)}}  
    \frac{\big|\power{u}{\frac{1+m}{2}} -
    (\power{u}{\frac{1+m}{2}})_{\rho}^{(\theta)}\big|^\frac{2r}{1+m}}
    {\rho^{\frac{r}{m}}}
    \,\dx\dt.
\end{equation*}
We now concentrate on the estimate of the first term on the right-hand side. 
In view of assumption~\eqref{sub-intrinsic-poincare} we have
\begin{align}\label{est-A-1}
    \mathrm{A}
    \le 
    C\biint_{Q_{2\rho}^{(\theta)}} 
  	\frac{|u|^{r}}{(2\rho)^{\frac{r}m}}\,\dx\dt
  	\le 
   	C \theta^\frac{2r m}{1+m},
\end{align}
where $C=C(N,m,r)$. 
On the other hand, in view of the sup-estimate from Lemma~\ref{lem:sup-int} and the Sobolev-Poincaré inequality from Lemma~\ref{lem:poin}, we can estimate $\mathrm{A}$ as
\begin{align}\label{est-A-2}
    \mathrm{A}
    &\le  
    C\bigg[\sup_{Q_\rho^{(\theta)}} 
    \frac{|u|}{\rho^{\frac{1}m}}\bigg]^{r -(1+m)}
   \cdot \biint_{Q_{\rho}^{(\theta)}}  
    \frac{\big|\power{u}{\frac{1+m}{2}} -
    (\power{u}{\frac{1+m}{2}})_{\rho}^{(\theta)}\big|^{2}}
    {\rho^{\frac{1+m}{m}}}
    \,\dx\dt \nonumber\\
    &\le  
    C\Bigg[\theta^{-\frac{2mN(1-m)}{(1+m)\boldsymbol \lambda_{p(1+m)}}} \bigg[
    \theta^{2pm} +
    \biint_{Q_{2\rho}^{(\theta)}} |F|^{2p} \,\dx\dt\bigg]^{\frac{2}{\boldsymbol \lambda_{p(1+m)}}}\Bigg]^{r-(1+m)} \cdot \mathrm{B} \nonumber\\
    &\le  
    C\bigg[\theta^{-\frac{2mN(1-m)}{(1+m)\boldsymbol \lambda_{p(1+m)}}} \big[
    \theta^{2m} + \mathrm{B} \big]^{\frac{2p}{\boldsymbol \lambda_{p(1+m)}}}\bigg]^{r-(1+m)} \cdot\mathrm{B},
\end{align}
where $C=C(N,m,\nu,L,p,r,K)$ and %we abbreviated 
\begin{equation*}
    \mathrm{B}
    := 
    \biint_{Q_{2\rho}^{(\theta)}} 
	|D\power{u}{m}|^{2}  \,\dx\dt +
	\bigg[\biint_{Q_{2\rho}^{(\theta)}} 
	|F|^{2p}  \,\dx\dt\bigg]^{\frac{1}{p}} .
\end{equation*}
We now resume to estimate the first term on the right-hand side of~\eqref{est-theta-2}. 
For a parameter $\beta\in(0,1)$ to be chosen later, we first decompose $\mathrm{A}$ into two factors, and then use~\eqref{est-A-1} to estimate the first resulting term, and~\eqref{est-A-2} for the second one. In this way we obtain 
\begin{align*}
    \mathrm{A}^{\frac{1+m}{2r}}
    &=
    \big[\mathrm{A}^{1-\beta} \mathrm{A}^\beta\big]^{\frac{1+m}{2r}} \\
    &\le 
    C\, \Bigg[\theta^\frac{2r m(1-\beta)}{1+m} 
    \bigg[\theta^{-\frac{2m N(1-m)}{(1+m)\boldsymbol \lambda_{p(1+m)}}} \big[
    \theta^{2m} + \mathrm{B} \big]^{\frac{2p}{\boldsymbol\lambda_{p(1+m)}}}\bigg]^{(r-(1+m))\beta} \mathrm{B}^\beta\Bigg]^{\frac{1+m}{2r}}.
\end{align*}
To cancel the first two factors of $\theta$, we now choose $\beta$ as to satisfy 
\begin{equation*}
    \frac{r(1-\beta)}{1+m}
    =
    \frac{N(1-m)}{\boldsymbol\lambda_{p(1+m)}}\frac{r-(1+m)}{1+m}
    \beta,
\end{equation*}
which is equivalent to 
\begin{equation*}
    \beta=
    \frac{r\boldsymbol \lambda_{p(1+m)}}{r\boldsymbol \lambda_{p(1+m)} + N(1-m)(r-(1+m))}
    \in(0,1).
\end{equation*}
We insert the definition of $\beta$, cancel the first two factors of $\theta$, and subsequently apply Young's inequality with exponents $\frac{r}{r-\beta(1+m)}$, $\frac{r}{\beta(1+m)}$ to obtain
\begin{align*}
    \mathrm{A}^{\frac{1+m}{2r}}
    &\le 
    C\, \bigg[ 
    \Big[\big[
    \theta^{2m} + \mathrm{B} \big]^{\frac{2p}{\boldsymbol\lambda_{p(1+m)}}}\Big]^{(r-(1+m))\beta} \mathrm{B}^\beta\bigg]^{\frac{1+m}{2r}} \\
    &\le 
    \big(1-2^{-\frac{1}{2m}} \big)
    \big[\theta^{2m} + \mathrm{B} \big]^{\frac{p(r-(1+m))}{\boldsymbol \lambda_{p(1+m)}}\cdot\frac{\beta(1+m)}{r -\beta (1+m)}} +
    C\, \mathrm{B}^{\frac12} \\
    &= 
    \big(1-2^{-\frac{1}{2m}}\big) 
    \big[\theta^{2m} + \mathrm{B} \big]^{\frac12} +
    C\, \mathrm{B}^{\frac12} \\
    &\le
    \big(1-2^{-\frac{1}{2m}}\big) \theta^{m} +
    C\, \mathrm{B}^{\frac12},
\end{align*}
where from the second to the third line we used the computation
\begin{align*}
    &\frac{p(r-(1+m))}{\boldsymbol \lambda_{p(1+m)}}\cdot\frac{\beta(1+m)}{r-\beta(1+m)}\\
    &\quad
    =
    \frac{p(r-(1+m))}{\boldsymbol \lambda_{p(1+m)}}\cdot
    \frac{\boldsymbol \lambda_{p(1+m)}(1+m)}{r\boldsymbol \lambda_{p(1+m)}+N(1-m)(r-(1+m))- \boldsymbol \lambda_{p(1+m)}(1+m) }
 \\
    &\quad=
    \frac{p(1+m)}{\boldsymbol\lambda_{p(1+m)} + N(1-m)} 
    =
    \frac12.
\end{align*}
Inserting this inequality into~\eqref{est-theta-2} and absorbing $\theta^{m}$ into the left-hand side, we find
\begin{align*}
    2^{-\frac{1}{2m}} \theta^{m}
    &\le 
    C\, \mathrm{B}^{\frac12} +
    \frac{\babs{(\power{u}{\frac{1+m}{2}})_{\rho/2}^{(\theta)}}}
    {\rho^{\frac{1+m}{2m}}} \\
    &\le 
    C\, \mathrm{B}^{\frac12} +
    \bigg[\biint_{Q_{\rho/2}^{(\theta)}}
    \frac{|u|^{1+m}}{\rho^{\frac{1+m}{m}}} \,\dx\dt 
    \bigg]^{\frac12},
\end{align*}
which yields the asserted inequality, once we take into account the definition of $\mathrm{B}$. 
\end{proof}

\section{Proof of the higher integrability}\label{sec:hi}

We consider a fixed cylinder
$$
	Q_{8R}(y_o,\tau_o)
	\equiv
	B_{8R}(y_o)\times \big(\tau_o-(8R)^\frac{1+m}{m},\tau_o+(8R)^\frac{1+m}{m}\big)
	\subseteq\Omega_T
$$
with $R>0$. Again, we omit the center in the notation and write
$Q_{\rho}:=Q_{\rho}(y_o,\tau_o)$ for short, for any radius
$\rho\in(0,8R]$. We consider  a parameter
\begin{equation}\label{lambda-0}
  \lambda_o
  \ge
  1+\bigg[\biint_{Q_{4R}} 
  \frac{\abs{u}^{r}}{(4 R)^{\frac{r}{m}}}\d x\d t
  \bigg]^{\frac{m+1}{ m\boldsymbol \lambda_r}},
\end{equation}
which will be fixed later. Recall that the {\it scaling deficit} 
$$
    d
    = 
    \frac{2r}{2r+N(m-1)}=\frac{2r}{\boldsymbol \lambda_r}
$$ 
is defined in \eqref{def:d}.
For $z_o\in Q_{2R}$, $\rho\in(0,R]$, and $\theta\ge 1$, we consider
space-time cylinders $Q_\rho^{(\theta)}(z_o)$ as defined in
\eqref{def-Q}. Note that these cylinders depend monotonically on $\theta$
in the sense that $Q_\rho^{(\theta_2)}(z_o)\subset
Q_\rho^{(\theta_1)}(z_o)$ whenever $1\le \theta_1<\theta_2$,
and that $Q_\rho^{(\theta)}(z_o)\subset Q_{4R}$ for $z_o\in Q_{2R}$, $\rho\in(0,R]$, and $\theta\ge 1$.

\subsection{Construction of a non-uniform system of cylinders}\label{sec:cylinders}
The following construction of a non-uniform system of cylinders has its origin in \cite{Gianazza-Schwarzacher, Schwarzacher}. Let $z_o\in Q_{2R}$. For a radius $\rho\in (0,R]$ we define
$$
	\widetilde\theta_\rho
	\equiv
	\widetilde\theta_{z_o;\rho}
	:=
	\inf\bigg\{\theta\in[\lambda_o,\infty):
	\frac{1}{|Q_\rho|}
	\iint_{Q^{(\theta)}_\rho(z_o)} 
	\frac{\abs{u}^{r}}{\rho^{\frac{r}m}} \dx\dt 
	\le 
	\theta^{\frac{ m\boldsymbol\lambda_r}{1+m}} \bigg\}.
$$
We note  that
$\widetilde\theta_\rho$ is well defined, since the
infimum in the definition is taken over a non-empty set. In fact,
in the limit $\theta\to\infty$ the integral on the left-hand side
converges to zero (and is constant in the case $m=1$, respectively), while the right-hand side grows 
with speed $\theta^{\frac{m\boldsymbol\lambda_r}{1+m}}$.
The choice of the exponent on the right-hand side becomes more clear after taking means in the
integral condition, since then the condition takes the form 
$$
  \biint_{Q^{(\theta)}_\rho(z_o)}
  \frac{\abs{u}^{r}}{\rho^{\frac{r}m}}\dx\dt 
  \le 
  \theta^\frac{2r m}{1+m}
$$
(compare to Sections~\ref{sec:poin} and~\ref{sec:revholder}). As an immediate consequence of the definition of
$\widetilde\theta_\rho$, we either have 
$$
	\widetilde\theta_\rho=\lambda_o
	\quad\mbox{and}\quad
	\biint_{Q_{\rho}^{(\widetilde\theta_\rho)}(z_o)}
    \frac{\abs{u}^{r}}{\rho^{\frac{r}m}} \dx\dt
	\le
	\widetilde\theta_\rho^\frac{2r m}{1+m}
	=
	\lambda_o^\frac{2r m}{1+m},
$$
or
\begin{equation}\label{theta>lambda}
	\widetilde\theta_\rho>\lambda_o
	\quad\mbox{and}\quad
	\biint_{Q_{\rho}^{(\widetilde\theta_\rho)}(z_o)} 
	\frac{\abs{u}^{r}}{\rho^{\frac{r}m}}\dx\dt
	=
	\widetilde\theta_\rho^\frac{2r m}{1+m}.
\end{equation}
In the case $\rho=R$,
we have $\widetilde \theta_{R}\ge \lambda_o\ge 1$.
Moreover, in the case $\widetilde\theta_{R}>\lambda_o$, property
\eqref{theta>lambda},
the inclusion $Q_{R}^{(\widetilde\theta_{R})}(z_o)\subset Q_{4R}$ and \eqref{lambda-0} yield that
\begin{align*}
	\widetilde\theta_{R}^{\frac{ m\boldsymbol\lambda_r}{1+m}}
	&=
	\frac{1}{|Q_{R}|}
	\iint_{Q_{R}^{(\widetilde\theta_{R})}(z_o)}
	\frac{\abs{u}^{r}}{R^{\frac{r}m}} \dx\dt\\ 
	&\le
	\frac{4^{\frac{r}m}}{|Q_{R}|}
	\iint_{Q_{4R}} \frac{\abs{u}^{r}}{(4R)^{\frac{r}m}}\d x\d t
	\le 
	4^{N+1+\frac{r+1}{m}}\lambda_o^{\frac{ m\boldsymbol \lambda_r}{1+m}},
\end{align*}
from which we infer  the bound
\begin{align}\label{bound-theta-R}
  \widetilde\theta_{R} 
  \le
  % 4^{\frac{d[m(n+2)+2]}{2m^2}}
  % \lambda_o
  % =
  4^{\frac{1+m}{ m\boldsymbol\lambda_r}(N+1+\frac{r+1}{m})}\lambda_o.
\end{align}     
Our next goal is to prove the continuity of the mapping $(0,R]\ni\rho\mapsto
\widetilde\theta_\rho$.
For $\rho\in(0,R]$ and  $\eps>0$, we abbreviate 
$\theta_+:=\widetilde\theta_\rho+\eps$. We first observe
that 
there exists $\delta=\delta(\eps,\rho)>0$ such that
for all radii ${\mathfrak r}\in(0,R]$ with $|\radius-\rho|<\delta$ there holds
\begin{equation}\label{claim-cont}
  \frac{1}{|Q_\radius|}
  \iint_{Q_{\radius}^{(\theta_+)}(z_o)} \frac{\abs{u}^{r}}{\radius^{\frac{r}m}} \dx\dt 
  <
  \theta_+^{\frac{ m\boldsymbol \lambda_r}{1+m}}.
\end{equation}
In fact, if 
$\radius=\rho$, this is a consequence of the definition of
$\widetilde\theta_\rho$, since $\theta_+^{\frac{ m\boldsymbol \lambda_r}{1+m}}>\widetilde\theta_\rho^{\frac{ m\boldsymbol \lambda_r}{1+m}}$.
By the absolute continuity of the integral, the inequality
\eqref{claim-cont} continues to hold for radii $\radius$ sufficiently close
to $\rho$. Hence, the definition of $\widetilde\theta_\radius$ implies that
$\widetilde\theta_\radius<\theta_+=\widetilde\theta_\rho+\eps$, provided $|\radius-\rho|<\delta$. 
For the corresponding lower bound
$\widetilde\theta_\radius>\theta_-:=\widetilde\theta_\rho-\eps$,
we proceed similarly.
First, we note that we can assume $\theta_-\ge\lambda_o$ and hence $\widetilde\theta_\rho>\lambda_o$, since
otherwise the claim immediately follows  from the property
$\widetilde\theta_\radius\ge\lambda_o$.
Then, we claim that 
\begin{equation}\label{claim-cont-2}
  	\frac{1}{|Q_\radius|}
	\iint_{Q_{\radius}^{(\theta_-)}(z_o)} \frac{\abs{u}^{r}}{\radius^{\frac{r}m}} \dx\dt 
  	>
  	\theta_-^{\frac{ m\boldsymbol \lambda_r}{1+m}}
\end{equation}
for all $\radius\in(0,R]$ with $|\radius-\rho|<\delta$, after diminishing
$\delta=\delta(\eps,\rho)>0$ if necessary.
Again, we first consider the case $\radius=\rho$, in which the claim follows
from the definition of $\widetilde\theta_\rho$. In fact, if the claim
did not hold, we
would arrive at the contradiction $\widetilde\theta_\rho\le\theta_-$.
Now, for radii
$\radius$ with $|\radius-\rho|<\delta$ the assertion 
follows from the continuous dependence
of the left-hand side upon $\radius$.
Having established \eqref{claim-cont-2}, we can conclude from 
the definition of $\widetilde\theta_\radius$ that
$\widetilde\theta_\radius>\theta_-=\widetilde\theta_\rho-\eps$. Altogether
we have shown that $\widetilde\theta_\rho-\eps<
\widetilde\theta_\radius< \widetilde\theta_\rho+\eps$  for all radii $\radius\in(0,R]$
with $|\radius-\rho|<\delta$, which
completes the proof of the continuity of $(0,R]\ni\rho\mapsto
\widetilde\theta_\rho$.

Unfortunately, the mapping $(0,R]\ni\rho\to \widetilde\theta_\rho$ might not be
decreasing.
For this reason we work with a modified version of
$\widetilde\theta_\rho$, which we denote by $\theta_\rho$. 
More precisely, we
define 
$$
	\theta_\rho
	\equiv
	\theta_{z_o;\rho}
	:=
	\max_{\radius\in[\rho,R]} \widetilde\theta_{z_o;\radius}\,. 
$$
As an immediate consequence of the construction, the mapping
$(0,R]\ni\rho\mapsto \theta_\rho$ is
continuous and monotonically decreasing. 
In general, the modified cylinders $Q_{\rho}^{(\theta_\rho)}(z_o)$ cannot 
be expected to be intrinsic in the sense of \eqref{theta>lambda}.  
However, we can show that the cylinders $Q_{s}^{(\theta_\rho)}(z_o)$
are sub-intrinsic for all radii $s\ge\rho$. More precisely, we we have 
\begin{align}\label{sub-intrinsic-2}
	\biint_{Q_{s}^{(\theta_{\rho})}(z_o)} 
	\frac{\abs{u}^{r}}{s^{\frac{r}m}} \dx\dt
	\le 
	\theta_\rho^\frac{2r m}{1+m}
	\quad\mbox{for any $0<\rho\le s\le R$.}
\end{align}
For the proof of this inequality, we use the chain of inequalities 
$\widetilde\theta_s\le \theta_{s}\le \theta_{\rho}$, which implies
$Q_{s}^{(\theta_{\rho})}(z_o)\subseteq
Q_{s}^{(\widetilde\theta_{s})}(z_o)$, and the fact that the latter
cylinder is sub-intrinsic (i.e., it satisfies \eqref{sub-intrinsic-poincare} with $K=1$). Subsequently, we recall that $\frac{2r m}{1+m} - \frac{Nm(1-m)}{1+m}=\frac{m\boldsymbol\lambda_r}{1+m}>0$ and again use $\widetilde\theta_s\le \theta_{\rho}$.
In this way, we deduce 
\begin{align*}
	\biint_{Q_{s}^{(\theta_{\rho})}(z_o)} 
	\frac{\abs{u}^{r}}{s^{\frac{r}m}} \dx\dt
	&\le 
	\Big(\frac{\theta_{\rho}}{\widetilde\theta_{s}}\Big)^{\frac{Nm(1-m)}{1+m}}
	\biint_{Q_{s}^{(\widetilde\theta_{s})}(z_o)}
    \frac{\abs{u}^{r}}{s^{\frac{r}m}} \dx\dt \\
	&\le 
	\Big(\frac{\theta_{\rho}}{\widetilde\theta_{s}}\Big)^{\frac{Nm(1-m)}{1+m}}
	\widetilde\theta_{s}^\frac{2r m}{1+m}
	=
	\theta_\rho^{\frac{Nm(1-m)}{1+m}}\,\widetilde\theta_s^{\frac{2r m}{1+m} - \frac{Nm(1-m)}{1+m}} \\
	&\le 
	\theta_\rho^\frac{2r m}{1+m},
\end{align*}
which is exactly assertion~\eqref{sub-intrinsic-2}. 

Next, we define 
\begin{equation}\label{rho-tilde}
	\widetilde\rho
	:=
	\left\{
	\begin{array}{cl}
	R, &
	\quad\mbox{if $\theta_\rho=\lambda_o$,} \\[5pt]
	\min\big\{s\in[\rho, R]: \theta_s=\widetilde \theta_s \big\}, &
	\quad\mbox{if $\theta_\rho>\lambda_o$.}
	\end{array}
	\right.
\end{equation}
By definition, for any $s\in [\rho,\widetilde\rho]$ we have
$\theta_s=\widetilde\theta_{\widetilde\rho}$.
% ; see again Figure~\ref{fig:sunrise}.
Our next goal is the proof of the upper bound 
\begin{align}\label{bound-theta}
  \theta_\rho 
  \le
  \Big(\frac{s}{\rho}\Big)^{\frac{1+m}{ m\boldsymbol\lambda_r}(N+1+\frac{r+1}{m})}
  \theta_{s}  
  \quad\mbox{for any $s\in(\rho,R]$.}
\end{align}
In the case $\theta_\rho=\lambda_o$ this is immediate since
$\theta_s\ge\lambda_o$. Another easy case is that of radii
$s\in(\rho,\widetilde\rho]$, since then we have
$\theta_s=\widetilde\theta_{\widetilde\rho}=\theta_\rho$. Therefore, it
only remains to prove \eqref{bound-theta}  
for the case $\theta_\rho>\lambda_o$ and radii $s\in(\widetilde\rho,R]$.
To this end, we use the monotonicity of $\rho\mapsto\theta_\rho$, \eqref{theta>lambda} and
\eqref{sub-intrinsic-2} to conclude 
\begin{align*}
	\theta_\rho^{\frac{ m\boldsymbol\lambda_r}{1+m}}
	&=
	\widetilde \theta_{\widetilde\rho}^{\frac{ m\boldsymbol\lambda_r}{1+m}}
	=
	\frac{1}{|Q_{\widetilde\rho}|}
	\iint_{Q_{\widetilde\rho}^{(\widetilde\theta_{\widetilde\rho})}(z_o)}
	\frac{\abs{u}^{r}}{\widetilde\rho^{\frac{r}m}} \dx\dt \\
	&\le
	\Big(\frac{s}{\widetilde\rho}\Big)^{N+1+\frac{r+1}{m}}
	\frac{1}{|Q_{s}|} 
	\iint_{Q_{s}^{(\theta_{s})}(z_o)} 
	\frac{\abs{u}^{r}}{s^{\frac{r}m}} \dx\dt\\ 
	&\le
	\Big(\frac{s}{\rho}\Big)^{N+1+\frac{r+1}{m}}
	\theta_{s}^{\frac{ m\boldsymbol\lambda_r}{1+m}} .
\end{align*}
This yields the claim~\eqref{bound-theta} also  in the remaining case.
We now apply \eqref{bound-theta} with $s=R$. Using the fact
$\theta_{R}=\widetilde\theta_{R}$ and
estimate \eqref{bound-theta-R} for $\widetilde\theta_{R}$, we deduce 
\begin{align}\label{bound-theta-2}
	\theta_\rho 
	\le
	\Big(\frac{R}{\rho}\Big)^{\frac{1+m}{ m\boldsymbol\lambda_r}(N+1+\frac{r+1}{m})}
	\theta_{R} 
	\le
	\Big(\frac{4R}{\rho}\Big)^{\frac{1+m}{ m\boldsymbol\lambda_r}(N+1+\frac{r+1}{m})}
	\lambda_o
\end{align}
for every $\rho\in(0,R]$.       
In summary, for every $z_o\in Q_{2R}$, we have constructed a
system of concentric sub-intrinsic cylinders
$Q_{\rho}^{(\theta_{z_o;\rho})}(z_o)$ with radii $\rho\in (0,R]$. As a
consequence of the monotonicity of $\rho\mapsto \theta_{z_o;\rho}$, these
cylinders are nested, in the sense that
$$
	Q_{\radius}^{(\theta_{z_o;\radius})}(z_o)
	\subset
	Q_{s}^{(\theta_{z_o;s})}(z_o),\,\,\text{ whenever }\,\,0<\radius<s\le R.
$$
However, keep in mind that in general these cylinders are not intrinsic but only sub-intrinsic.

\subsection{Covering property}%\label{sec:covering}
Our next goal is to establish the following Vitali-type covering property for the cylinders constructed in the previous section. 

\begin{lemma}\label{lem:vitali}
There exists a constant $\hat c=\hat c(N,m,r)\ge 20$ 
such that, whenever $\mathcal F$ is a collection of cylinders
$Q_{4\radius}^{(\theta_{z;\radius})}(z)$,
where $Q_{\radius}^{(\theta_{z;\radius})}(z)$ is a cylinder as constructed in Section~{\upshape\ref{sec:cylinders}} with radius $\radius\in(0,\tfrac{R}{\hat c})$, there exists a countable subfamily $\mathcal G$ of disjoint cylinders in $\mathcal F$ such that  
\begin{equation}\label{covering}
	\bigcup_{Q\in\mathcal F} Q
	\subseteq 
	\bigcup_{Q\in\mathcal G} \widehat Q,
\end{equation}
where $\widehat Q$ denotes the $\frac{\hat c}{4}$-times enlarged cylinder $Q$, i.e.~if $Q=Q_{4\radius}^{(\theta_{z;\radius})}(z)$, then $\widehat Q=Q_{\hat c \radius}^{(\theta_{z;\radius})}(z)$.
\end{lemma}

\begin{proof}
For the moment, suppose that $\hat c\ge20$ has already been determined. For $j\in\N$ we subdivide $\mathcal{F}$ into the subfamilies 
$$
  \mathcal F_j
  :=
  \big\{Q_{4\radius}^{(\theta_{z;\radius})}(z)\in \mathcal F: 
  \tfrac{R}{2^j\hat c}<\radius\le \tfrac{R}{2^{j-1}\hat c} \big\}.
$$
Then, we choose finite subfamilies $\mathcal G_j\subset \mathcal F_j$
according to the following scheme.
We start by choosing $\mathcal G_1$ as an arbitrary maximal disjoint collection of cylinders in
$\mathcal F_1$. The subfamily $\mathcal G_1$ is finite, since
\eqref{bound-theta-2} and the definition of $\mathcal F_1$ imply a
lower bound on the volume of each cylinder in $\mathcal G_1$.
Now, assuming that the subfamilies $\mathcal G_1, \mathcal G_2,
\dots, \mathcal G_{k-1}$ have already been constructed for some integer
$k\ge 2$, we choose $\mathcal G_k$ to be any maximal disjoint subcollection of 
$$
	\Bigg\{Q\in \mathcal F_k: 
	Q\cap Q^\ast=\emptyset\, \mbox{ for any }\, 
    Q^\ast\in \bigcup_{j=1}^{k-1} \mathcal G_j 
	\Bigg\}.
$$
For the same reason as above, the collection $\mathcal G_k$ is
finite. Hence, the family 
$$
	\mathcal G
	:=
	\bigcup_{j=1}^\infty \mathcal G_j\subseteq\mathcal{F}
$$
defines a countable collection of disjoint cylinders.
It remains to prove that
for each cylinder  $Q\in\mathcal F$ there exists a cylinder $Q^\ast\in\mathcal
G$ with  $Q\subset \widehat {Q}^\ast$.  
To this end, we fix a cylinder
$Q=Q_{4\radius}^{(\theta_{z;\radius})}(z)\in\mathcal F$.
Let $j\in\N$ be such that $Q\in\mathcal F_j$.
The maximality of $\mathcal G_j$ ensures the existence of a cylinder
$Q^\ast=Q_{4\radius_\ast}^{(\theta_{z_\ast;\radius_\ast})}(z_\ast)\in
\bigcup_{i=1}^{j} \mathcal G_i$ with $Q\cap Q^\ast\not=\emptyset$.
We will show that this cylinder has the desired property $Q\subset \widehat {Q}^\ast$.
First, we observe that the properties $\radius\le\tfrac{R}{2^{j-1}\hat c}$
and $\radius_\ast>\tfrac{R}{2^j\hat c}$ imply $\radius\le 2\radius_\ast$, which ensures
$\Lambda_{4\radius}(t)\subseteq \Lambda_{\hat c\radius_\ast}(t_\ast)$ since $\hat c\ge20$. For the proof of the corresponding spatial inclusion
$B^{(\theta_{z,\radius})}_{4\radius}(x)\subseteq B_{\hat c \radius_\ast}^{(\theta_{z_\ast,\radius_\ast})}(x_\ast)$, we
first shall derive the bound 
\begin{equation}\label{control-theta-1}
  \theta_{z_\ast;\radius_\ast}
  \le
  52^{\frac{1+m}{ m\boldsymbol\lambda_r}(N+1+\frac{r+1}{m})}\,
  \theta_{z;\radius}\,.
\end{equation}
We recall the definition \eqref{rho-tilde} of the radius $\widetilde
\radius_\ast\in [\radius_\ast,R]$ which is 
associated to the cylinder
$Q_{\radius_\ast}^{(\theta_{z_\ast;\radius_\ast})}(z_\ast)$. According to the
definition, we either have that $Q_{\widetilde
  \radius_\ast}^{(\theta_{z_\ast;\radius_\ast})}(z_\ast)$ is intrinsic
or  that $\widetilde \radius_\ast=R$ and $\theta_{z_\ast;\radius_\ast}=\lambda_o$. In
the second alternative, the claim~\eqref{control-theta-1} is
immediate, since 
$$
	\theta_{z_\ast;\radius_\ast}
	=
	\lambda_o
	\le 
	\theta_{z;\radius}\,.
$$
Therefore, it remains to consider the case where $Q_{\widetilde
  \radius_\ast}^{(\theta_{z_\ast;\radius_\ast})}(z_\ast)$ is intrinsic in the
sense that 
\begin{align}\label{control-theta-2}
	\theta_{z_\ast;\radius_\ast}^{\frac{ m\boldsymbol\lambda_r}{1+m}}
	=
	\frac{1}{|Q_{\widetilde \radius_\ast}|}
	\iint_{Q_{\widetilde \radius_\ast}^{(\theta_{z_\ast;\radius_\ast})}(z_\ast)}
	\frac{\abs{u}^{r}}{\widetilde \radius_\ast^{\frac{r}m}} \dx\dt.
\end{align}
We distinguish between the cases  $\widetilde \radius_\ast\le \frac{R}{\mu}$ and $\widetilde \radius_\ast> \frac{R}{\mu}$, where $\mu:= 13$. 
We start with the latter case. Using \eqref{control-theta-2} and the
definition of $\lambda_o$ and $\theta_{z;\radius}$, we estimate
\begin{align*}
	\theta_{z_\ast;\radius_\ast}^{\frac{ m\boldsymbol\lambda_r}{1+m}}
	&\le
	\Big(\frac{4R}{\widetilde \radius_\ast}\Big)^{\frac{r}m} \frac{1}{|Q_{\widetilde \radius_\ast}|}
	\iint_{Q_{4R}}
	\frac{\abs{u}^{r}}{(4R)^{\frac{r}m}} \dx\dt \\
	&\le
	\Big(\frac{4 R}{\widetilde \radius_\ast}\Big)^{N+1+\frac{r+1}{m}} 
	\lambda_o^{\frac{ m\boldsymbol\lambda_r}{1+m}}\\
	&\le
	(4\mu)^{N+1+\frac{r+1}{m}} \theta_{z;\radius}^{\frac{ m\boldsymbol\lambda_r}{1+m}},
\end{align*}
which can be rewritten in the form
\begin{align*}
	\theta_{z_\ast;\radius_\ast}
	\le
	(4\mu)^{\frac{1+m}{m\boldsymbol\lambda_r}(N+1+\frac{r+1}{m})}\,\theta_{z;\radius}\,.
\end{align*}
This yields \eqref{control-theta-1} in the second case, and it only
remains to consider the first case $\widetilde \radius_\ast\le \frac{R}{\mu}$. Here, the
key step is to prove the inclusion
\begin{equation}\label{incl-cyl}
  Q_{\widetilde \radius_\ast}^{(\theta_{z_\ast;\radius_\ast})}(z_\ast)
  \subseteq 
  Q_{\mu\widetilde \radius_\ast}^{(\theta_{z;\mu\widetilde \radius_\ast})}(z).
\end{equation}
We first observe
that $\widetilde \radius_\ast\ge \radius_\ast$ and
$|t-t_\ast|<(4\radius)^{\frac{1+m}m}+(4\radius_\ast)^{\frac{1+m}m}\le
(12\radius_\ast)^{\frac{1+m}m}$ implies $\Lambda_{\widetilde \radius_\ast}(t_\ast)\subseteq \Lambda_{\mu\widetilde
  \radius_\ast}(t)$.  In addition, we have 
\begin{equation}\label{x-x_ast}
	|x-x_\ast|\le \theta_{z;\radius}^{\frac{m(m-1)}{1+m}}4\radius
	+
	\theta_{z_\ast;\radius_\ast}^{\frac{m(m-1)}{1+m}}4\radius_\ast.
\end{equation}
At this point, we may assume that
$\theta_{z;\radius}\le\theta_{z_\ast;\radius_\ast}$, since \eqref{control-theta-1}
clearly is satisfied in the alternative  case.
Then, the monotonicity of $\rho\mapsto \theta_{z;\rho}$  and 
$\radius\le 2\radius_\ast\le2\widetilde \radius_\ast\le\mu\widetilde \radius_\ast$ imply
$$
	\theta_{z_\ast;\radius_\ast}
	\ge 
	\theta_{z;\radius}
	\ge
	\theta_{z;\mu \widetilde \radius_\ast}.
$$
Combining this with \eqref{x-x_ast}, we conclude  that
\begin{align*}
	\theta_{z_\ast;\radius_\ast}^{\frac{m(m-1)}{1+m}}\widetilde \radius_\ast + |x-x_\ast|
	\le
    \theta_{z_\ast;\radius_\ast}^{\frac{m(m-1)}{1+m}}5\widetilde \radius_\ast+
    \theta_{z;\radius}^{\frac{m(m-1)}{1+m}}4\radius 
	\le
	\theta_{z;\mu\widetilde \radius_\ast}^{\frac{m(m-1)}{1+m}}\mu\widetilde \radius_\ast,
\end{align*}
from which we deduce the inclusion 
$$
  B_{\widetilde \radius_\ast}^{(\theta_{z_\ast;\radius_\ast})}(x_\ast)
  \subseteq 
  B_{\mu\widetilde \radius_\ast}^{(\theta_{z;\mu\widetilde \radius_\ast})}(x).
$$
This completes the proof of \eqref{incl-cyl}.

Using \eqref{control-theta-2}, \eqref{incl-cyl}, and
\eqref{sub-intrinsic-2} with $\rho=s=\mu\tilde \radius_\ast$, we estimate 
\begin{align*}
	\theta_{z_\ast;\radius_\ast}^{\frac{ m\boldsymbol\lambda_r}{1+m}}
	\le
	\frac{\mu^{\frac{r}{m}}}{|Q_{\widetilde \radius_\ast}|} 
	\iint_{Q_{\mu\widetilde \radius_\ast}^{(\theta_{z;\mu\widetilde \radius_\ast})}(z)}
	\frac{\abs{u}^{r}}{(\mu\tilde \radius_\ast)^{\frac{r}{m}}} \dx\dt
	\le
	\mu^{N+1+\frac{r+1}{m}} \theta_{z;\radius}^{\frac{ m\boldsymbol\lambda_r}{1+m}},
\end{align*}
which implies 
\begin{align*}
	\theta_{z_\ast;\radius_\ast}
	\le
	\mu^{\frac{1+m}{m\boldsymbol\lambda_r}(N+1+\frac{r+1}{m})}\,\theta_{z;\radius}.
\end{align*}
This yields \eqref{control-theta-1} also in the last 
case.

Having \eqref{control-theta-1} at hand, we resume to prove $Q=Q_{4\radius}^{(\theta_{z;\radius})}(z)\subseteq
Q_{\hat c \radius_\ast}^{(\theta_{z_\ast;\radius_\ast})}(z_\ast)=\widehat Q^\ast$, which will
give \eqref{covering} and complete the proof. First, recall
that we have shown
$\Lambda_{4\radius}(t)\subseteq\Lambda_{\hat c
\radius_\ast}(t_\ast)$, and it only remains to compare the spatial balls. Moreover, from \eqref{x-x_ast} and
\eqref{control-theta-1}, and taking into account that $\radius\le 2\radius_\ast$,  we conclude
\begin{align*}
	\theta_{z;\radius}^{\frac{m(m-1)}{1+m}}4\radius + |x-x_\ast|
	&\le
	2\theta_{z;\radius}^{\frac{m(m-1)}{1+m}}4\radius + 
	\theta_{z_\ast;\radius_\ast}^{\frac{m(m-1)}{1+m}}4\radius_\ast \\
	&\le
	4\Big[4\cdot 52^{\frac{1-m}{\boldsymbol\lambda_r}(N+1+\frac{r+1}{m})}+1\Big]
	\theta_{z_\ast;\radius_\ast}^{\frac{m(m-1)}{1+m}}\radius_\ast \\
	&\le
	\theta_{z_\ast;\radius_\ast}^{\frac{m(m-1)}{1+m}}\hat c \radius_\ast,
\end{align*}
for a suitable constant $\hat c=\hat c(N,m,r)$. This implies the spatial inclusion
$B_{4\radius}^{(\theta_{z;\radius})}(x)\subseteq B_{\hat c
\radius_\ast}^{(\theta_{z_\ast;\radius_\ast})}(x_\ast)$, which yields the desired cylinder inclusion $Q\subset \widehat Q^\ast$.
\end{proof}

\subsection{Stopping time argument}
Now, we fix the parameter $\lambda_o$ by letting
\begin{equation}\label{Eq:def-lm0}
\begin{aligned}
 	\lambda_o
 	&:=
 	1+\Bigg[
    \bigg[\biint_{Q_{4R}} 
 	\frac{\abs{u}^{r}}{(4 R)^{\frac{r}m}} \,\dx\dt \bigg]^{\frac{1+m}{r}}\\
    &\qquad\qquad + 
    \biint_{Q_{4R}}
    |D\power{u}{m}|^2 \,\dx\dt +
 	\bigg[ \biint_{Q_{4R}} |F|^{2p}\,\dx\dt \bigg]^{\frac1p} 
    \Bigg]^{\frac{r}{m\boldsymbol\lambda_r}}.
\end{aligned}
\end{equation}
For $\lambda>\lambda_o$ and $\radius\in(0,2R]$, we define the super-level set
of the function $|D\power{u}{m}|$ by 
\[
	\boldsymbol E(\radius,\lambda)
	:=
	\Big\{z\in Q_{\radius}:
	\mbox{$z$ is a Lebesgue point of $|D\power{u}{m}|$ and 
	$|D\power{u}{m}|(z) > \lambda^{m}$}\Big\}.
\]
In the definition of $\boldsymbol E(\radius,\lambda)$, the notion of Lebesgue points is to be understood
with regard to the system of cylinders built in
Section~\ref{sec:cylinders}. We point out that also with respect to
these cylinders, $\mathcal L^{N+1}$-a.e. point is a Lebesgue
point. This follows from \cite[2.9.1]{Federer}, since we have already 
verified the Vitali-type covering Lemma~\ref{lem:vitali}. 

Now, we fix radii $R\le R_1<R_2\le 2R$. Note that for any
$z_o\in Q_{R_1}$, $\kappa\ge1$ and $\rho\in(0,R_2-R_1]$ we have
\begin{equation*}
  	Q_\rho^{(\kappa)}(z_o)
	\subseteq 
	Q_{R_2}
	\subseteq 
	Q_{2R}.
\end{equation*}
For the following argument we restrict ourselves to levels $\lambda$
with 
\begin{equation}\label{choice_lambda}
	\lambda
	>
	B\lambda_o,
	\quad\mbox{where}
	\quad
	B
	:=
	\Big(\frac{4\hat c R}{R_2-R_1}\Big)^{\frac{r(N+2)(1+m)}{2m^2\boldsymbol\lambda_r}}
	>1,
\end{equation}
and $\hat c=\hat c(N,m,r)$ is the constant from the Vitali-type
covering Lemma \ref{lem:vitali}.
We fix $z_o\in \boldsymbol E(R_1,\lambda)$ and abbreviate
$\theta_s\equiv \theta_{z_o;s}$ for $s\in(0,R]$ throughout this
section. By definition of $\boldsymbol E(R_1,\lambda)$, we have
\begin{align}\label{larger-lambda}
	\liminf_{s\downarrow 0} 
	\Bigg[ 
    \biint_{Q_{s}^{(\theta_{s})}(z_o)} 
	|D\power{u}{m}|^2 \,\dx\dt + 
    \bigg[ \biint_{Q_{s}^{(\theta_{s})}(z_o)} |F|^{2p} \,\dx\dt\bigg]^{\frac1p} \Bigg] 
	&\ge
	|D\power{u}{m}|^2(z_o) \nonumber \\
	&>
	\lambda^{2m}.
\end{align}
On the other hand, for any radius $s$ with
\begin{align}\label{radius-s}
	\frac{R_2-R_1}{\hat c}\le s\le R
\end{align}
the definition of $\lambda_o$, estimate
\eqref{bound-theta-2}, assumption \eqref{radius-s} and the definition
of $d$ imply 
\begin{align}\label{smaller-lambda}
	&\biint_{Q_{s}^{(\theta_{s})}(z_o)} 
	|D\power{u}{m}|^2 \,\dx\dt + 
    \bigg[ \biint_{Q_{s}^{(\theta_{s})}(z_o)} |F|^{2p} \,\dx\dt\bigg]^{\frac1p} \nonumber\\
	&\qquad\le
	\frac{|Q_{4R}|}{|Q_{s}^{(\theta_s)}|}
	\biint_{Q_{4R}} |D\power{u}{m}|^2 \,\dx\dt + 
    \bigg[ \frac{|Q_{4R}|}{|Q_{s}^{(\theta_s)}|}
    \biint_{Q_{4R}} |F|^{2p} \,\dx\dt\bigg]^{\frac1p} \nonumber\\
	&\qquad\le
	\frac{|Q_{4R}|}{|Q_{s}^{(\theta_s)}|}
	\Bigg[ 
    \biint_{Q_{4R}} |D\power{u}{m}|^2 \,\dx\dt + 
    \bigg[ 
    \biint_{Q_{4R}} |F|^{2p} \,\dx\dt\bigg]^{\frac1p} \Bigg]\nonumber\\
	&\qquad\le
	\frac{|Q_{4R}|}{|Q_{s}|}\,
	\theta_s^{\frac{Nm(1-m)}{1+m}}
	\lambda_o^{\frac{m\boldsymbol\lambda_r}{1+m}} 
     \nonumber\\
	&\qquad\le
	\Big(\frac{4R}{s}\Big)^{N+1+\frac{1}m+\frac{1+m}{ m\boldsymbol\lambda_r}(N+1+\frac{r+1}{m})\cdot\frac{Nm(1-m)}{1+m}} 
	\lambda_o^{\frac{Nm(1-m)}{1+m}+\frac{m\boldsymbol\lambda_r}{1+m}} \nonumber\\
    &\qquad=
	\Big(\frac{4R}{s}\Big)^{\frac{r(N+2)(1+m)}{m\boldsymbol\lambda_r}} 
	\lambda_o^\frac{2r m}{1+m} \nonumber\\
	&\qquad\le
%	\Big(\frac{4\hat c R}{R_2-R_1}\Big)^{\frac{d(n+2)(1+m)}{2m}} \lambda_o^{2m} \nonumber\\
%	&=
	B^{2m} \lambda_o^{2m} 
	<
	\lambda^{2m},
\end{align}
since
\begin{align*}
    N+1+\frac{1}m&+\frac{1+m}{ m\boldsymbol\lambda_r}\Big(N+1+\frac{r+1}{m}\Big)\cdot\frac{Nm(1-m)}{1+m}\\
    &=
    N+1+\frac{1}m+\frac{N(1-m)}{\boldsymbol\lambda_r}\Big(N+1+\frac{r+1}{m}\Big)\\
    &=
    \Big( N+1+\frac{1}m\Big) 
    \Big( 1+\frac{N(1-m)}{\boldsymbol\lambda_r}\Big)
    +\frac{N(1-m)r}{m\boldsymbol\lambda_r}\\
    &=
    \Big( N+1+\frac{1}m\Big)\frac{2r}{\boldsymbol\lambda_r}
    +\frac{N(1-m)r}{m\boldsymbol\lambda_r}\\
    &=
    \frac{r}{m\boldsymbol\lambda_r}\big(2Nm +2m +2+N-Nm\big)\\
    &= \frac{r(N+2)(1+m)}{m\boldsymbol\lambda_r}.
\end{align*}
By the continuity of the mapping $s\mapsto\theta_s$ and the absolute
continuity of the integral, the left-hand side of
\eqref{smaller-lambda} depends continuously on $s$. Therefore, in view
of \eqref{larger-lambda} and \eqref{smaller-lambda}, there exists
a maximal radius $0<\rho_{z_o} < \tfrac{R_2-R_1}{\hat c}$ for which
the above inequality becomes an equality, i.e.~$\rho_{z_o}$ is the
maximal radius with 
\begin{align}\label{=lambda}
	\biint_{Q_{\rho_{z_o}}^{(\theta_{\rho_{z_o}})}(z_o)} 
	|D\power{u}{m}|^2 \,\dx\dt + \bigg[\biint_{Q_{\rho_{z_o}}^{(\theta_{\rho_{z_o}})}(z_o)} |F|^{2p} \,\dx\dt\bigg]^{\frac1p}
	=
	\lambda^{2m}.
\end{align}
The maximality of the radius $\rho_{z_o}$ implies in particular that 
\begin{align*}
	\biint_{Q_{s}^{(\theta_{s})}(z_o)} 
	|D\power{u}{m}|^2 \,\dx\dt + \bigg[\biint_{Q_{s}^{(\theta_{s})}(z_o)} |F|^{2p} \,\dx\dt\bigg]^{\frac1p} 
	<
	\lambda^{2m}
	\quad
	\mbox{for any $s\in (\rho_{z_o}, R]$.}
\end{align*}
Due to the monotonicity of the mapping $\rho\mapsto \theta_\rho$ and \eqref{bound-theta} we have 
\begin{align*}
	\theta_\sigma
	\le 
	\theta_{s}
	\le
	\Big(\frac{\sigma}{s}\Big)^{\frac{1+m}{ m\boldsymbol\lambda_r}(N+1+\frac{r+1}{m})}
  	\theta_{\sigma}
	\quad
	\mbox{for any $\rho_{z_o}\le s<\sigma\le R$,}
\end{align*}
so that 
\begin{align}\label{<lambda}
	&\biint_{Q_{\sigma}^{(\theta_{s})}(z_o)} 
	|D\power{u}{m}|^2 \,\dx\dt + \bigg[\biint_{Q_{\sigma}^{(\theta_{s})}(z_o)} |F|^{2p} \,\dx\dt\bigg]^{\frac1p} \nonumber\\
	&\quad\le
	\Big(\frac{\theta_{s}}{\theta_\sigma}\Big)^{\frac{Nm(1-m)}{1+m}}
    \Bigg[ 
	\biint_{Q_{\sigma}^{(\theta_{\sigma})}(z_o)} 
	|D\power{u}{m}|^2 \,\dx\dt + \bigg[ \biint_{Q_{\sigma}^{(\theta_{\sigma})}(z_o)} |F|^{2p} \,\dx\dt\bigg]^{\frac1p} 
    \Bigg] \nonumber\\
	&\quad<
	\Big(\frac{\sigma}{s}\Big)^{\frac{N(1-m)}{\boldsymbol\lambda r}(N+1+\frac{r+1}{m})}\,
	\lambda^{2m}
\end{align}
for any $\rho_{z_o}\le s<\sigma\le R$. 
Finally, we recall that the cylinders are constructed in such a
way that 
\begin{equation*}
  Q_{\hat c\rho_{z_o}}^{(\theta_{\rho_{z_o}})}(z_o)
  \subseteq
  Q_{\hat c\rho_{z_o}}(z_o)
  \subseteq
  Q_{R_2}.
\end{equation*}

\subsection{A Reverse H\"older Inequality}
For a level  $\lambda$ as in \eqref{choice_lambda} and a point $z_o\in
\boldsymbol E(R_1,\lambda)$, we consider the radius
$\widetilde\rho_{z_o}\in[\rho_{z_o},R]$ as defined in
\eqref{rho-tilde}. In the sequel
we write $\theta_{\rho_{z_o}}$ instead of $\theta_{z_o;\rho_{z_o}}$.
We recall that $\widetilde\rho_{z_o}$ has been defined in such a way
that for any $s\in [\rho_{z_o},
\widetilde\rho_{z_o}]$ we have
$\theta_s=\theta_{\rho_{z_o}}$,  and, in particular,
$\theta_{\widetilde\rho_{z_o}}=\theta_{\rho_{z_o}}$.

The aim of this section is the proof of a reverse
H\"older inequality on 
$Q_{2\rho_{z_o}}^{(\theta_{\rho_{z_o}})}(z_o)$. To this end, we need to
verify the assumptions of Proposition \ref{prop:revhoelder}. 
First, we note that \eqref{sub-intrinsic-2} with
$s=2\rho_{z_o}$ implies
\begin{equation*}
  \biint_{Q_{2\rho_{z_o}}^{(\theta_{\rho_{z_o}})}(z_o)} 
  \frac{\abs{u}^{r}}{(2\rho_{z_o})^{\frac{r}m}} \dx\dt
  \le 
  \theta_{\rho_{z_o}}^\frac{2r m}{1+m},
\end{equation*}
which means that assumption \eqref{sub-intrinsic-poincare} is fulfilled for the cylinder
$Q_{2\rho_{z_o}}^{(\theta_{\rho_{z_o}})}(z_o)$ with $K=1$.
For the estimate of $\theta_{\rho_{z_o}}^{2m}$ from above,  
we distinguish between the cases $\widetilde\rho_{z_o}\le 2\rho_{z_o}$
and $\widetilde\rho_{z_o}> 2\rho_{z_o}$.
In the former case, we use the fact
$\theta_{\rho_{z_o}}=\theta_{\widetilde\rho_{z_o}}=\widetilde\theta_{\widetilde\rho_{z_o}}$, which implies
that 
$Q_{\widetilde\rho_{z_o}}^{(\theta_{\rho_{z_o}})}(z_o)$ is intrinsic, 
and then the bound $\widetilde\rho_{z_o}\le 2\rho_{z_o}$, with the
result that
\begin{equation*}
  \theta_{\rho_{z_o}}^\frac{2r m}{1+m}
  =
  \biint_{Q_{\widetilde\rho_{z_o}}^{(\theta_{\rho_{z_o}})}(z_o)}
  \frac{\abs{u}^{r}}{\widetilde\rho_{z_o}^{\frac{r}m}}
  \dx\dt
  \le
  2^{N+\frac{r}m}
  \biint_{Q_{2\rho_{z_o}}^{(\theta_{\rho_{z_o}})}(z_o)}
  \frac{\abs{u}^{r}}{(2\rho_{z_o})^{\frac{r}m}}
  \dx\dt.
\end{equation*}
This means that in this case, assumption \eqref{super-intrinsic-poincare}$_1$
is satisfied with $K\equiv  2^{N+\frac{r}m}$.

Next, we consider the remaining case 
$\widetilde\rho_{z_o}>2\rho_{z_o}$.
First, we claim that
\begin{equation}
  \label{claim-case-2}
  \theta_{\rho_{z_o}}\le C(N,m,\nu,L,p,r)\lambda.
\end{equation}
For the proof we treat the cases 
$\widetilde\rho_{z_o}\in(2\rho_{z_o},\frac R2]$ and
$\widetilde\rho_{z_o}\in(\frac R2,R]$ separately. 
In the latter case, we use
\eqref{bound-theta-2} with $\rho=\widetilde\rho_{z_o}$ and  
the bound $\widetilde\rho_{z_o}>\frac R2$ in order to estimate
$\theta_{\rho_{z_o}}=\theta_{\widetilde\rho_{z_o}}\le C\lambda_o \le C\lambda$, 
which yields \eqref{claim-case-2}. 

In the alternative case $\widetilde\rho_{z_o}\in(2\rho_{z_o},\frac
R2]$, the cylinder
$Q_{\widetilde\rho_{z_o}}^{(\theta_{\rho_{z_o}})}(z_o)$ is intrinsic
by definition of $\widetilde\rho_{z_o}$, and the two times enlarged  cylinder
$Q_{2\widetilde\rho_{z_o}}^{(\theta_{\rho_{z_o}})}(z_o)$ is
sub-intrinsic by \eqref{sub-intrinsic-2}.
Therefore, assumptions \eqref{sub-intrinsic-poincare} and \eqref{super-intrinsic-poincare}$_1$ of Lemma~\ref{lem:theta} are satisfied for $Q_{\widetilde\rho_{z_o}}^{(\theta_{\rho_{z_o}})}(z_o)$. The application of the lemma yields
\begin{align*}
  	\theta_{\rho_{z_o}}^{m}
  	&\le
	\tfrac1{\sqrt2} \Bigg[
	\biint_{Q_{\widetilde\rho_{z_o}/2}^{(\theta_{\rho_{z_o}})}(z_o)}\!
    \frac{|u|^{1+m}}{(\widetilde\rho_{z_o}/2)^{\frac{1+m}{m}}} \dx\dt 
    \Bigg]^{\frac1{2}} \\
    &\quad +
  	C \Bigg[\biint_{Q_{2\widetilde\rho_{z_o}}^{(\theta_{\rho_{z_o}})}(z_o)} 
   	|D\power{u}{m}|^{2} \,\dx\dt +
    \bigg[\biint_{Q_{2\widetilde\rho_{z_o}}^{(\theta_{\rho_{z_o}})}(z_o)} |F|^{2p} \,\dx\dt\bigg]^{\frac1p} \Bigg]^{\frac12}
\end{align*}
with  a constant $C=C(N,m,\nu,L,p,r)$. For the first term, we exploit the
sub-intrinsic coupling \eqref{sub-intrinsic-2} with radii
$\rho=\rho_{z_o}$ and $s=\frac12\widetilde\rho_{z_o}>\rho_{z_o}$.
For the estimate of the last integral, we
recall that $\theta_{\widetilde\rho_{z_o}}=\theta_{\rho_{z_o}}$,
which allows to use \eqref{<lambda} with $s=\widetilde\rho_{z_o}$ and
$\sigma=2\widetilde\rho_{z_o}\in(s,R]$. This leads to the upper bound
\begin{align*}
  	\theta_{\rho_{z_o}}^{m}
  	\le
	\tfrac1{\sqrt2} \theta_{\rho_{z_o}}^m +
  	C\,\lambda^m .
\end{align*}
Here, we re-absorb $\tfrac1{\sqrt2} \theta_{\rho_{z_o}}^m$ into
the left and obtain the claim
\eqref{claim-case-2} in any case.

Combining this with the identity~\eqref{=lambda}, we obtain the bound  
\begin{align*}
  	\theta_{\rho_{z_o}}^{2m}
  	\le
  	C \Bigg[\biint_{Q_{2\widetilde\rho_{z_o}}^{(\theta_{\rho_{z_o}})}(z_o)} 
   	|D\power{u}{m}|^{2} \,\dx\dt +
    \bigg[\biint_{Q_{2\widetilde\rho_{z_o}}^{(\theta_{\rho_{z_o}})}(z_o)} |F|^{2p} \,\dx\dt\bigg]^{\frac1p} \Bigg].
\end{align*}
This means that in the case $\widetilde\rho_{z_o}>2\rho_{z_o}$
assumption \eqref{super-intrinsic-poincare}$_2$ is satisfied  on
$Q_{2\rho_{z_o}}^{(\theta_{\rho_{z_o}})}$ with a constant $K\equiv K(N,m,\nu,L,p,r)$.
In any case, we have eventually shown that all hypotheses of Proposition~\ref{prop:revhoelder} 
are satisfied. Consequently, the proposition yields the desired reverse
H\"older inequality 
\begin{align}\label{rev-hoelder}
	\biint_{Q_{2\rho_{z_o}}^{(\theta_{\rho_{z_o}})}(z_o)} & 
	|D\power{u}{m}|^2 \dx\dt \nonumber\\
	&\le
	C\bigg[\biint_{Q_{4\rho_{z_o}}^{(\theta_{\rho_{z_o}})}(z_o)} 
	|D\power{u}{m}|^{2q} \dx\dt \bigg]^{\frac{1}{q}} +
	C\bigg[ \biint_{Q_{4\rho_{z_o}}^{(\theta_{\rho_{z_o}})}(z_o)} |F|^{2p} \,\dx\dt\bigg]^{\frac1p},
\end{align}
for a constant $C=C(N,m,\nu,L,,p,r)$ and with   $q=\frac{rN}{rN+{\boldsymbol\lambda}_r}<1$, as given in \eqref{def:q}.

\subsection{Estimate on super-level sets}
So far we have shown that 
for every $\lambda$ as in  \eqref{choice_lambda} and every $z_o\in
\boldsymbol E(R_1,\lambda)$, there exists a  cylinder
$Q_{\rho_{z_o}}^{(\theta_{z_o;\rho_{z_o}})}(z_o)$ with 
$Q_{\hat c\rho_{z_o}}^{(\theta_{z_o;\rho_{z_o}})}(z_o)\subseteq Q_{R_2}$, 
for
which the properties \eqref{=lambda} and \eqref{<lambda},  and
the reverse H\"older-type estimate \eqref{rev-hoelder} are satisfied. 
This allows us to establish
a reverse H\"older inequality for the distribution
function of $|D\power{u}{m}|^2$
by a Vitali-type covering argument. The precise argument is as follows.
We define the super-level set of the inhomogeneity $F$ by 
$$
	\boldsymbol F(\radius,\lambda)
	:=
	\Big\{z\in Q_{\radius}: 
	\mbox{$z$ is a Lebesgue point of $F$ and 
	$|F|(z)>\lambda^{m}$}\Big\},
$$
where again, the Lebesgue points have to be understood with respect to the cylinders built in  Section~\ref{sec:cylinders}.
Using \eqref{=lambda} and \eqref{rev-hoelder}, we estimate
\begin{align*}
	\lambda^{2m}
	&=
	\biint_{Q_{\rho_{z_o}}^{(\theta_{\rho_{z_o}})}(z_o)} 
	|D\power{u}{m}|^2 \,\dx\dt + \bigg[ \biint_{Q_{\rho_{z_o}}^{(\theta_{\rho_{z_o}})}(z_o)} |F|^{2p} \,\dx\dt\bigg]^{\frac1p} \\
	&\le
	C\bigg[\biint_{Q_{4\rho_{z_o}}^{(\theta_{\rho_{z_o}})}(z_o)} 
	|D\power{u}{m}|^{2q} \dx\dt \bigg]^{\frac{1}{q}} +
	C\bigg[  \biint_{Q_{4\rho_{z_o}}^{(\theta_{\rho_{z_o}})}(z_o)} |F|^{2p} \,\dx\dt\bigg]^{\frac1p} \\
	&\le
	C\,\eta^{2m}\lambda^{2m} +
	C\Bigg[\frac{1}{\big|Q_{4\rho_{z_o}}^{(\theta_{\rho_{z_o}})}(z_o)\big|}
	\iint_{Q_{4\rho_{z_o}}^{(\theta_{\rho_{z_o}})}(z_o)\cap \boldsymbol E(R_2,\eta\lambda)} 
	|D\power{u}{m}|^{2q} \dx\dt \Bigg]^{\frac{1}{q}} \\
	&\quad+
	C\Bigg[ 
    \frac{1}{\big|Q_{4\rho_{z_o}}^{(\theta_{\rho_{z_o}})}(z_o)\big|}
	\iint_{Q_{4\rho_{z_o}}^{(\theta_{\rho_{z_o}})}(z_o)\cap \boldsymbol F(R_2,\eta\lambda)} 
	|F|^{2p} \dx\dt\Bigg]^{\frac1p} ,
\end{align*}
for a constant $C=C(N,m,\nu,L,C_o,p,r)$ and any $\eta\in(0,1)$.
For the estimate of the second last term, we apply 
H\"older's inequality and \eqref{<lambda}, with the result
\begin{align*}
	\Bigg[\frac{1}{\big|Q_{4\rho_{z_o}}^{(\theta_{\rho_{z_o}})}(z_o)\big|} &
	\iint_{Q_{4\rho_{z_o}}^{(\theta_{\rho_{z_o}})}(z_o)\cap \boldsymbol E(R_2,\eta\lambda)} 
	|D\power{u}{m}|^{2q} \dx\dt \Bigg]^{\frac{1}{q}-1} \\
	&\le
	\bigg[\biint_{Q_{4\rho_{z_o}}^{(\theta_{\rho_{z_o}})}(z_o)} 
	|D\power{u}{m}|^{2} \dx\dt \bigg]^{1-q}
	\le
	C\,\lambda^{2m(1-q)};
\end{align*}
for the last term we use Young's inequality with the result
\begin{align*}
	C\Bigg[ 
    \frac{1}{\big|Q_{4\rho_{z_o}}^{(\theta_{\rho_{z_o}})}(z_o)\big|} &
	\iint_{Q_{4\rho_{z_o}}^{(\theta_{\rho_{z_o}})}(z_o)\cap \boldsymbol F(R_2,\eta\lambda)} 
	|F|^{2p} \dx\dt\Bigg]^{\frac1p} \\
    &\le 
    \tfrac14 \lambda^{2m} +
    \frac{C\,\lambda^{-2m(p-1)}}{\big|Q_{4\rho_{z_o}}^{(\theta_{\rho_{z_o}})}(z_o)\big|} 
	\iint_{Q_{4\rho_{z_o}}^{(\theta_{\rho_{z_o}})}(z_o)\cap \boldsymbol F(R_2,\eta\lambda)} 
	|F|^{2p} \dx\dt .
\end{align*}
We now choose the parameter $\eta\in(0,1)$ in dependence on $N$, $m$,
$\nu$, $L$, $p$, $r$
in such a way that $\eta^{2m}=\frac{1}{4C}$. This allows us
to re-absorb the terms containing $\lambda^{2m}$ into the
left-hand side. 
The  resulting inequality is then
multiplied  by
$\big|Q_{4\rho_{z_o}}^{(\theta_{\rho_{z_o}})}(z_o)\big|$. In this way,
we arrive at
\begin{align*}
	\lambda^{2m}\big|Q_{4\rho_{z_o}}^{(\theta_{\rho_{z_o}})}(z_o)\big|
	&\le
	C\iint_{Q_{4\rho_{z_o}}^{(\theta_{\rho_{z_o}})}(z_o)\cap \boldsymbol E(R_2,\eta\lambda)} 
	\lambda^{2m(1-q)}|D\power{u}{m}|^{2q} \dx\dt \\
	&\quad+
	C
	\iint_{Q_{4\rho_{z_o}}^{(\theta_{\rho_{z_o}})}(z_o)\cap \boldsymbol F(R_2,\eta\lambda)} 
	\lambda^{-2m(p-1)} |F|^{2p} \dx\dt,
\end{align*}
where $C=C(N,m,\nu,L,p,r)$. 
On the other hand, we bound the left-hand side
from below by the use of \eqref{<lambda}. This leads to the
inequality
\begin{align*}
	\lambda^{2m}
    &\ge
    C
	\biint_{Q_{\hat c\rho_{z_o}}^{(\theta_{\rho_{z_o}})}(z_o)} 
	|D\power{u}{m}|^2 \dx\dt,
\end{align*}
where we relied on the fact that $\hat c$ is a universal constant depending only on $N$, $m$, and $r$.       
Combining the two preceding estimates and using  again that $\hat
c=\hat c(N,m,r)$, we arrive at
\begin{align}\label{level-est}
	\iint_{Q_{\hat c\rho_{z_o}}^{(\theta_{\rho_{z_o}})}(z_o)} 
	|D\power{u}{m}|^2 \d x\d t 
	&\le
    C
    \iint_{Q_{4\rho_{z_o}}^{(\theta_{\rho_{z_o}})}(z_o)\cap \boldsymbol E(R_2,\eta\lambda)} 
	\lambda^{2m(1-q)}|D\power{u}{m}|^{2q} \dx\dt \nonumber\\
	&\quad +
	C
	\iint_{Q_{4\rho_{z_o}}^{(\theta_{\rho_{z_o}})}(z_o)\cap \boldsymbol  F(R_2,\eta\lambda)} 
	\lambda^{-2m(p-1)} |F|^{2p} \dx\dt
\end{align}
for a constant $C=C(N,m,\nu,L,p,r)$.
Since the preceding inequality holds for every center $z_o\in \boldsymbol
E(R_1,\lambda)$, we conclude that it is possible to cover the
super-level set
$\boldsymbol E(R_1,\lambda)$ by a  family
$\mathcal F\equiv\big\{Q_{4\rho_{z_o}}^{(\theta_{z_o;\rho_{z_o}})}(z_o)\big\}$ of parabolic cylinders 
with center $z_o\in \boldsymbol E(R_1,\lambda)$, such that each of the
cylinders is contained in $Q_{R_2}$, and
on each cylinder estimate \eqref{level-est} is valid.
An application of the Vitali-type covering Lemma \ref{lem:vitali}
provides us with a countable disjoint subfamily
$$
	\Big\{Q_{4\rho_{z_i}}^{(\theta_{z_i;\rho_{z_i}})}(z_i)\Big\}_{i\in\N}
	\subseteq 
	\mathcal F
$$
with the property
$$
	\boldsymbol E(R_1,\lambda)
	\subseteq 
	\bigcup_{i=1}^\infty Q_{\hat c\rho_{z_i}}^{(\theta_{z_i;\rho_{z_i}})}(z_i)
	\subseteq
	Q_{R_2}.
$$
We apply \eqref{level-est} for each of the cylinders
$Q_{4\rho_{z_i}}^{(\theta_{z_i;\rho_{z_i}})}(z_i)$ and add the
resulting inequalities. Since the cylinders
$Q_{4\rho_{z_i}}^{(\theta_{z_i;\rho_{z_i}})}(z_i)$ are pairwise
disjoint, we obtain 
\begin{align*}
	\iint_{\boldsymbol E(R_1,\lambda)} 
	|D\power{u}{m}|^2 \dx\dt
	&\le
	C\iint_{\boldsymbol E(R_2,\eta\lambda)} 
	\lambda^{2m(1-q)}|D\power{u}{m}|^{2q} \dx\dt \\
    &\quad +
	C
    \iint_{\boldsymbol F(R_2,\eta\lambda)} 
	\lambda^{-2m(p-1)}|F|^{2p} \dx\dt,
\end{align*}
with $C=C(N,m,\nu,L,p,r)$. In order to compensate for the fact that
super-level sets of different levels appear on both sides of the
preceding estimate, we need an
estimate on the difference $\boldsymbol E(R_1,\eta\lambda)\setminus
\boldsymbol E(R_1,\lambda)$.
However, on this set we can simply estimate
$|D\power{u}{m}|^2\le\lambda^{2m}$, which leads to the bound 
\begin{align*}
	\iint_{\boldsymbol E(R_1,\eta\lambda)\setminus \boldsymbol E(R_1,\lambda)} 
	|D\power{u}{m}|^2 \d x\d t
	&\le
	\iint_{\boldsymbol E(R_2,\eta\lambda)} 
	\lambda^{2m(1-q)}|D\power{u}{m}|^{2q} \dx\dt.
\end{align*}
Adding the last two inequalities, we obtain a reverse
H\"older-type inequality for the distribution function of $|D\power{u}{m}|^2$ for levels $\eta\lambda$. 
In this inequality we replace  $\eta\lambda$ by $\lambda$ and recall that $\eta\in(0,1)$
was chosen as a universal constant depending only on $N$, $m$, $\nu$, $L$, $p$, $r$. In this way, we obtain for any $\lambda\ge \eta B\lambda_o
=:\lambda_1$ that
\begin{align}\label{pre-1}
	\iint_{\boldsymbol E(R_1,\lambda)} &
	|D\power{u}{m}|^2 \d x\d t \nonumber\\
	&\le
	C\iint_{\boldsymbol E(R_2,\lambda)} 
	\lambda^{2m(1-q)}|D\power{u}{m}|^{2q} \dx\dt +
	C \iint_{\boldsymbol F(R_2,\lambda)} \lambda^{-2m(p-1)}|F|^{2p} \dx\dt
\end{align}
with a constant $C=C(N,m,\nu,L,p,r)$. This is the desired
reverse H\"older-type inequality for the distribution function of $D\power{u}{m}$.

\subsection{Proof of the gradient estimate of Theorem~\ref{thm:higherint}}

For  $k> \lambda_1$ we define 
the {\it truncation} of $|D\power{u}{m}|$ by
$$
	|D\power{u}{m}|_k
	:=
	\min\big\{|D\power{u}{m}|, k^m\big\},
$$
and for $\radius\in(0,2R]$ the corresponding super-level set
$$
	\boldsymbol E_k(\radius,\lambda)
	:=
	\big\{z\in Q_\radius: |D\power{u}{m}|_k>\lambda^{m}\big\}.
$$
Note that $|D\power{u}{m}|_k\le |D\power{u}{m}|$ a.e., as well as $\boldsymbol E_k(\radius,\lambda)=\emptyset$ for $k\le\lambda$ and $\boldsymbol E_k(\radius,\lambda)=\boldsymbol E(\radius,\lambda)$ for $k>\lambda$. Therefore, 
it follows from \eqref{pre-1} that
\begin{align*}
	&\iint_{\boldsymbol E_k(R_1,\lambda)} 
	|D\power{u}{m}|_k^{2-2q}|D\power{u}{m}|^{2q} \d x\d t\\
	&\qquad\le
	C\iint_{\boldsymbol E_k(R_2,\lambda)} 
	\lambda^{2m(1-q)}|D\power{u}{m}|^{2q} \dx\dt +
	C \iint_{\boldsymbol F(R_2,\lambda)} \lambda^{-2m(p-1)}|F|^{2p} \dx\dt,
\end{align*}
whenever $k>\lambda$. Since $\boldsymbol E_k(\radius,\lambda )=\emptyset$ for $k\le \lambda$,
the last inequality  also holds in this case.
Now, we multiply the preceding inequality by $\lambda^{2m\epsilon -1}$, where $\epsilon\in (0,1]$ will be chosen later
in a universal way, and  integrate  the result  with respect to $\lambda$ over the interval $(\lambda_1,\infty)$. This gives
\begin{align*}
	\int_{\lambda_1}^\infty \lambda^{2m\epsilon-1} &
	\bigg[\iint_{\boldsymbol E_k(R_1,\lambda)} 
	|D\power{u}{m}|_k^{2-2q}
	|D\power{u}{m}|^{2q} \,\dx\dt \bigg] \d\lambda
	\\
	&\le
	C\int_{\lambda_1}^\infty \lambda^{2m(1-q+\epsilon)-1} 
	\bigg[
	\iint_{\boldsymbol E_k(R_2,\lambda)} 
	|D\power{u}{m}|^{2q} \,\dx\dt\bigg]\d\lambda \\
	&\quad +
	C \int_{\lambda_1}^\infty \lambda^{2m(1-p+\epsilon)-1}
	\bigg[
	\iint_{\boldsymbol F(R_2,\lambda)} |F|^{2p} \,\dx\dt\bigg]\d\lambda.
\end{align*}
Here we exchange the order of integration with the help of Fubini's theorem.
For the integral on the left-hand side Fubini's theorem implies
\begin{align*}
	\int_{\lambda_1}^\infty &\lambda^{2m\epsilon-1} 
	\bigg[
	\iint_{\boldsymbol E_k(R_1,\lambda)} 
	|D\power{u}{m}|_k^{2-2q}
	|D\power{u}{m}|^{2q} \,\dx\dt \bigg]\,\d\lambda
	\\
	&=
	\iint_{\boldsymbol E_k(R_1,\lambda_1)}
	|D\power{u}{m}|_k^{2-2q}
	|D\power{u}{m}|^{2q}
	\bigg[
	\int_{\lambda_1}^{|D\powerexp{u}{m}|_k^{\frac{1}{m}}} 
	\lambda^{2m\epsilon-1} \,\d\lambda\bigg]
	 \,\dx\dt \\
	&=
	\tfrac{1}{2m\epsilon} \iint_{\boldsymbol E_k(R_1,\lambda_1)}
	\Big[|D\power{u}{m}|_k^{2-2q+2\epsilon}
	|D\power{u}{m}|^{2q}%\\
%	&\phantom{=\,}
	-
	\lambda_1^{2m\epsilon} |D\power{u}{m}
	|_k^{2-2q}|D\power{u}{m}|^{2q} \Big]
	\,\dx\dt ,
\end{align*}
while for the first integral on the right-hand side we recall $q<1$ and find that
\begin{align*}
	\int_{\lambda_1}^\infty &\lambda^{2m(1-q+\epsilon)-1} 
	\bigg[
	\iint_{\boldsymbol E_k(R_2,\lambda)} 
	|D\power{u}{m}|^{2q} \,\dx\dt\bigg]\,\d\lambda \\
	&=
	\iint_{\boldsymbol E_k(R_2,\lambda_1)} |D\power{u}{m}|^{2q}
	\bigg[
	\int_{\lambda_1}^{|D\powerexp{u}{m}|_k^{\frac{1}{m}}} 
	\lambda^{2m(1-q+\epsilon)-1} \,\d\lambda\bigg]
	\,\dx\dt \\
	&\le 
	\tfrac{1}{2m(1-q+\epsilon)} \iint_{\boldsymbol E_k(R_2,\lambda_1)} 
	|D\power{u}{m}|_k^{2-2q+2\epsilon} 
	|D\power{u}{m}|^{2q}
	\,\dx\dt \\
	&\le 
	\tfrac{1}{2m(1-q)} \iint_{\boldsymbol E_k(R_2,\lambda_1)} 
	|D\power{u}{m}|_k^{2-2q+2\epsilon} 
	|D\power{u}{m}|^{2q}
	\,\dx\dt.
\end{align*}
Finally, for the second integral on the right-hand side we note that $1-p+\epsilon< -\frac{N}2+\epsilon\le-\frac12<0$, since $p>\frac{N+2}2$ and $N\ge3$ and obtain
\begin{align*}
	\int_{\lambda_1}^\infty & \lambda^{2m(1-p+\epsilon)-1}
	\bigg[\iint_{\boldsymbol F(R_2,\lambda)} |F|^{2p} \,\dx\dt\bigg]\,\d\lambda \\
	&=
	\iint_{\boldsymbol F(R_2,\lambda_1)} |F|^{2p}
	\bigg[\int_{\lambda_1}^{ |F|^{\frac{1}{m}}} 
	\lambda^{2m(1-p+\epsilon)-1} \,\d\lambda\bigg]
	 \,\dx\dt \\
	&\le
	\frac{\lambda_1^{-2m(p-1-\epsilon)}}{2m(p-1-\epsilon)} \iint_{\boldsymbol F(R_2,\lambda_1)}
	 |F|^{2p} \,\dx\dt \\
	 &\le
	\frac{\lambda_1^{-2m(p-1-\epsilon)}}{m}
    \iint_{Q_{2R}}
	 |F|^{2p} \,\dx\dt. 
\end{align*}
We insert these estimates above and multiply by $2m\epsilon$. This leads to 
\begin{align*}
	 \iint_{\boldsymbol E_k(R_1,\lambda_1)}&
	|D\power{u}{m}|_k^{2-2q+2\epsilon}
	|D\power{u}{m}|^{2q} \,\dx\dt \\
	&\le
	\lambda_1^{2m\epsilon} 
	\iint_{\boldsymbol E_k(R_1,\lambda_1)}
	|D\power{u}{m}|_k^{2-2q}
	|D\power{u}{m}|^{2q} \,\dx\dt \\
	&\phantom{\le\,}+
	\frac{C\, \epsilon}{1-q} \iint_{\boldsymbol E_k(R_2,\lambda_1)} 
	|D\power{u}{m}|_k^{2-2q+2\epsilon} 
	|D\power{u}{m}|^{2q}
	\,\dx\dt \\
	&\phantom{\le\,}+
	\frac{C\,\epsilon}{\lambda_1^{2m(p-1-\epsilon)}}
    \iint_{Q_{2R}} |F|^{2p} \,\dx\dt .	
\end{align*}
The last inequality is now combined with the corresponding inequality on the complement $Q_{R_1}\setminus \boldsymbol E_k(R_1,\lambda_1)$,
i.e.~with the inequality
\begin{align*}
	\iint_{Q_{R_1}\setminus \boldsymbol E_k(R_1,\lambda_1)}&
	|D\power{u}{m}|_k^{2-2q+2\epsilon}
	|D\power{u}{m}|^{2q} \,\dx\dt \\
	&\le
	\lambda_1^{2m\epsilon} 
	\iint_{Q_{R_1}\setminus \boldsymbol E_k(R_1,\lambda_1)}
	|D\power{u}{m}|_k^{2-2q}
	|D\power{u}{m}|^{2q} \,\dx\dt .
\end{align*}
We also  take into account that $|D\power{u}{m}|_k\le |D\power{u}{m}|$.  All together this gives the inequality
\begin{align*}
	 \iint_{Q_{R_1}}&
	|D\power{u}{m}|_k^{2-2q+2\epsilon}
	|D\power{u}{m}|^{2q} \,\dx\dt \\
	&\le
	\frac{C_\ast \epsilon}{1-q} \iint_{Q_{R_2}} 
	|D\power{u}{m}|_k^{2-2q+2\epsilon} 
	|D\power{u}{m}|^{2q}
	\,\dx\dt \\
	&\phantom{\le\ }+
	\lambda_1^{2m\epsilon} 
	\iint_{Q_{2R}} |D\power{u}{m}|^{2} \dx\dt +
	\frac{C\,\epsilon}{\lambda_1^{2m(p-1-\epsilon)}}
    \iint_{Q_{2R}} |F|^{2p} \,\dx\dt,
\end{align*}
where $C_\ast=C_\ast(N,m,\nu,L,p,r)\ge 1$. Now, we choose
$$
	0
	<
	\epsilon
	\le
	\epsilon_o,
	\quad\mbox{where }
	\epsilon_o
	:=
	\frac{1-q}{2C_\ast}
	<\frac12.
$$
Note that $\epsilon_o$ depends only on $N, m, \nu, L, p, r$. Moreover, observe that $\lambda_1^{p-1-\epsilon}\equiv (\eta B\lambda_o)^{p-1-\epsilon} \ge \eta^{2p-2} \lambda_o^{p-1-\epsilon}$ and $\lambda_1^\epsilon\equiv (\eta B\lambda_o)^\epsilon\le B \lambda_o^\epsilon$, since $\eta\le1$, $B\ge 1$ and $0<\epsilon\le 1$. Therefore, from the previous inequality we conclude that for any pair of radii $R_1$, $R_2$ with  $R\le R_1<R_2\le 2R$ there holds
\begin{align*}
	 \iint_{Q_{R_1}}&
	|D\power{u}{m}|_k^{2-2q+2\epsilon}
	|D\power{u}{m}|^{2q} \,\dx\dt \\
	&\le
	\tfrac{1}{2} \iint_{Q_{R_2}} 
	|D\power{u}{m}|_k^{2-2q+2\epsilon} 
	|D\power{u}{m}|^{2q}
	\,\dx\dt \\
	&\phantom{\le\ }+
	C\,\Big(\frac{R}{R_2{-}R_1}\Big)^{\frac{r(N+2)(1+m)}{m{\boldsymbol\lambda_r}}}
	\lambda_o^{2m\epsilon} 
	\iint_{Q_{2R}} |D\power{u}{m}|^{2} \,\dx\dt \\
    &\quad +
	\frac{C}{\lambda_o^{2m(p-1-\epsilon)}}
    \iint_{Q_{2R}} |F|^{2p} \,\dx\dt.
\end{align*}
We can now apply the Iteration Lemma~\ref{lem:tech-classical} to the last inequality, which yields
\begin{align*}
	\iint_{Q_{R}}
	|D\power{u}{m}|_k^{2-2q+2\epsilon}
	|D\power{u}{m}|^{2q} \,\dx\dt 
	&\le
	C\, 
	\lambda_o^{2m\epsilon} 
	\iint_{Q_{2R}}
	|D\power{u}{m}|^{2} \,\dx\dt \\
	&\quad +
	\frac{C}{\lambda_o^{2m(p-1-\epsilon)}}
    \iint_{Q_{2R}} |F|^{2p} \,\dx\dt.
\end{align*}
On the left side we apply Fatou's Lemma and pass to the limit $k\to\infty$.
In the result, we go over to  means on both sides.  This gives 
\begin{align*}
	\biint_{Q_{R}} 
	|D\power{u}{m}|^{2+2\epsilon} \,\dx\dt 
	&\le
	C\, 
	\lambda_o^{2m\epsilon} 
	\biint_{Q_{2R}}
	|D\power{u}{m}|^{2} \,\dx\dt +
	\frac{C}{\lambda_o^{2m(p-1-\epsilon)}}
    \biint_{Q_{2R}} |F|^{2p} \,\dx\dt .
\end{align*}
At this stage it remains to bound  $\lambda_o$ from above and from below. From its definition in \eqref{Eq:def-lm0} we infer
\begin{equation*}
 	\lambda_o
 	\ge 
 	\max\Bigg\{1,
 	\bigg[ \biint_{Q_{4R}} |F|^{2p}\,\dx\dt \bigg]^{\frac{d}{2 mp}} \Bigg\},
\end{equation*}
so that 
\begin{align*}
	&\frac{1}{\lambda_o^{2m(p-1-\epsilon)}}
    \biint_{Q_{2R}} |F|^{2p} \,\dx\dt \\
    &\qquad=
    \frac{1}{\lambda_o^{2m(p-1-\epsilon)}}
    \bigg[\biint_{Q_{2R}} |F|^{2p} \,\dx\dt \bigg]^{\frac{p-1-\epsilon}{p}} 
    \bigg[\biint_{Q_{2R}} |F|^{2p} \,\dx\dt \bigg]^{\frac{1+\epsilon}{p}} \\
    &\qquad\le 
    \frac{C}{\lambda_o^{2m(p-1-\epsilon)(1-\frac1d)}} 
    \bigg[\biint_{Q_{2R}} |F|^{2p} \,\dx\dt \bigg]^{\frac{1+\epsilon}{p}} \\
    &\qquad\le  
    C\bigg[\biint_{Q_{2R}} |F|^{2p} \,\dx\dt \bigg]^{\frac{1+\epsilon}{p}}.
\end{align*}
On the other hand, by an application  of the energy estimate from Lemma \ref{lem:energy} with $\theta=1$
and $a=0$ and Young's inequality, we have
\begin{equation*}
 	\lambda_o
 	\le 
 	C\Bigg[1+
    \bigg[ \biint_{Q_{8R}} 
 	\frac{\abs{u}^{r}}{(4 R)^{\frac{r}m}}  \,\dx\dt\bigg]^{\frac{1+m}{r}} +
 	\bigg[ \biint_{Q_{8R}} |F|^{2p}\,\dx\dt \bigg]^{\frac1p} 
    \Bigg]^{\frac{r}{m\boldsymbol\lambda_r}}.
\end{equation*}
where $C=C(N,m,\nu,L)$.
Plugging this into the preceding estimate, we arrive at 
\begin{align*}
	&\biint_{Q_{R}} 
	|D\power{u}{m}|^{2+2\epsilon} \,\dx\dt \\
	&\quad\le
	C
	\Bigg[
	1+\bigg[ \biint_{Q_{8R}} 
 	\frac{\abs{u}^{r}}{R^{\frac{r}m}}  \,\dx\dt\bigg]^{\frac{1+m}{r}}  +
 	\bigg[ \biint_{Q_{8R}} |F|^{2p}\,\dx\dt \bigg]^{\frac1p}
	\Bigg]^{\epsilon \frac{2r}{\boldsymbol\lambda_r}}
	\biint_{Q_{2R}}
	|D\power{u}{m}|^{2} \dx\dt\\
        &\qquad+
	C\bigg[ \biint_{Q_{2R}}  |F|^{2p} \dx\dt \bigg]^{\frac{1+\epsilon}{p}},
\end{align*}
for a constant $C=C(N,m,\nu,L,p,r)$. The asserted reverse H\"older inequality
\eqref{eq:higher-int} now follows by a standard covering argument. More
precisely, we cover $Q_R$ by finitely many cylinders of radius $\frac
R8$, apply the preceding estimate on each of the smaller cylinders and
take the sum. 
This yields the same estimate as above, but with integrals over
$Q_{2R}$ instead of $Q_{8R}$. This completes the proof of
Theorem~\ref{thm:higherint}.
\qed

\subsection{Proof of Corollary \ref{cor:higher-int}}
We consider a standard parabolic
cylinder $C_{2R}(z_o):=
B_{2R}(x_o)\times(t_o-(2R)^2,t_o+(2R)^2)\subseteq\Omega_T$ and rescale
the problem via
\begin{equation*}
	\left\{
	\begin{array}{c}
  	v(x,t):=u(x_o+Rx,t_o+R^2t) \\[5pt]
	\boldsymbol B(x,t,u,\xi):=R\,\boldsymbol A\big(x_o+Rx,t_o+R^2t,u,\tfrac1{R}\xi\big) \\[5pt]
	G(x,t):=R\,F(x_o+Rx,t_o+R^2t),
	\end{array}
	\right.
\end{equation*}
whenever $(x,t)\in C_{2}$ and $(u,\xi)\in\R^N\times\R^{Nk}$.
The rescaled function $v$ is a weak solution to the differential equation
\begin{equation*}
    \partial_tv-\mathrm{div}\,\boldsymbol B(x,t,v,D\power{v}{m})=\mathrm{div}\,G
    \qquad\mbox{in $\widetilde Q:=Q_{2^\frac{2m}{1+m}}\subseteq C_{2}$}
  \end{equation*}
  in the sense of Definition~\ref{def:weak_solution}. Moreover, the rescaled
  vector-field $\mathbf B$ satisfies
  assumptions~\eqref{growth}. Consequently, we can apply estimate
  \eqref{eq:higher-int} to $v$ on the cylinder $\widetilde Q$, which gives 
   \begin{align*}
	\biint_{\frac12\widetilde Q} &
          |D\power{v}{m}|^{2+2\epsilon} \d x\d t \\
	&\le
	C\Bigg[
	1+\bigg[\biint_{\widetilde Q} \abs{v}^{r}\dx\dt\bigg]^\frac{1+m}{r} 
    + 
    \bigg[
    \biint_{\widetilde Q}|G|^{2p}\d x\d t 
	\bigg]^\frac1p\Bigg]^{\epsilon\frac{2r}{\boldsymbol\lambda_r}}
	\biint_{\widetilde Q} |D\power{v}{m}|^{2} \dx\dt\\
    &\phantom{\le\,}+
	C\bigg[\biint_{\widetilde Q} |G|^{2p} \dx\dt\bigg]^\frac{1+\epsilon}{p},    
  \end{align*}
  for every $\eps\in(0,\eps_o]$, where $C=C(N,m,\nu,L,p,r)$, and we
  abbreviated $\frac12\widetilde Q:=Q_{2^{\frac{m-1}{1+m}}}$.   
  Scaling back and using the fact that the cylinder $C_{\gamma R}$
  with $\gamma:=2^\frac{m-1}{2m}$ is contained in the re-scaled
  version of the cylinder $\frac12\widetilde Q$ 
  we deduce
   \begin{align*}
	R^{2+2\eps} &
	\biint_{C_{\gamma R}(z_o)}
    |D\power{u}{m}|^{2+2\epsilon} \d x\d t \\
	&\le
	C\,R^2\Bigg[
	1+\bigg[ \biint_{C_{2R}(z_o)} \abs{u}^{r}\dx\dt\bigg]^\frac{1+m}{r}
    +
    R^2\bigg[\biint_{C_{2R}(z_o)}
    |F|^{2p} \d x\d t \bigg]^\frac1p
	\Bigg]^{\epsilon\frac{2r}{\boldsymbol\lambda_r}}\\
 &\qquad\qquad\qquad\qquad\qquad\qquad\cdot
	\biint_{C_{2R}(z_o)}
	|D\power{u}{m}|^{2} \dx\dt \\
	&\quad +
	C\, R^{2+2\epsilon}
    \bigg[
	\biint_{C_{2R}(z_o)} |F|^{2p} \dx\dt\bigg]^\frac{1+\epsilon}{p} .
   \end{align*}
The asserted estimate on the pair $C_R$, $C_{2R}$ of standard parabolic cylinders now follows with a standard covering argument.  This 
finishes  the proof of Corollary~\ref{cor:higher-int}. \qed

\appendix
\appendixpage

\section{Significance of Hypothesis~\eqref{growth}$_3$}\label{sec:app-dia}

\begin{remark}\upshape\label{rem:C_o}
Concerning the structure condition \eqref{growth}$_3$ some remarks are in order. First of all, it must be noted that even in the case $m=1$ conditions \eqref{growth}$_1$ and \eqref{growth}$_2$ are not sufficient for the boundedness of weak solutions, as classical counterexamples in the stationary case show; cf.~\cite {DeGiorgi-1, Giusti-Miranda, Necas-Stara}. In a certain sense, one needs a diagonal structure of the system, which is encoded in hypothesis \eqref{growth}$_3$. In fact, in \cite{Meier} Meier pointed out that the indicator function
\begin{align*}
    \boldsymbol I \big(x,t,u,\xi \big) 
    &=
    \sum_{\al=1}^N
    \boldsymbol A_\alpha(x,t,u,\xi)\cdot \frac{u}{|u|}\, \xi_\al\cdot \frac{u}{|u|}
\end{align*}
is of crucial importance for a precise understanding of the local behavior of
weak solutions to elliptic systems. In \cite[Theorem~1]{Meier} he assumes a lower bound for $\boldsymbol I$ for all $(x,t,u,\xi)$ of the form
$(x,t,u(x,t), Du(x,t))$, where
$(x,t)\in \{ |u|>L\}$. In our setting  this would mean that the indicator function $\boldsymbol I (x,t,u(x,t), D\power{u}{m}(x,t))$ is non-negative on the solution $u$
for large values of the solution itself. %$u(x,t)$. 
However, we do not pursue generality this direction, which would certainly be possible, and
instead we assume the non-negativity of the indicator function for all values of $u$ (and not only for its large values).
We are able to prove by a Moser iteration scheme a qualitative local $L^\infty$-bound for weak solutions to the porous medium-type system in the sub-critical range under the additional assumption that $u$ is integrable enough. 

However, the qualitative, local $L^\infty$-estimate is not sufficient to establish the higher integrability of the gradient of solutions; a quantitative local $L^\infty$-estimate is required for that. The proof is based on an energy inequality for $(|u|^m-k^m)_+^2$, which is essential for the De Giorgi approach for equations of the porous medium type. For this we refer to the monograph~\cite{DBGV-book}, in which the use of such energy inequalities is extensively discussed. 
They represent the core of most regularity results, such as the expansion of positivity, H\"older continuity of weak solutions, Harnack's inequality, to name just a few.

The difference from the scalar case with non-negative solutions, where it is essentially possible to test with $(u^m-k^m)_+$ or powers of $u$, is that in the case of systems $\power{u}{m}|u|^\alpha$
respectively $\power{u}{m}(|u|^m-k^m)_+^2$ have to be implemented as testing functions. As a consequence, such testing functions naturally produce an integral whose integrand contains the indicator function. If \eqref{growth}$_3$ is fulfilled, this term has a sign and can be discarded. 
\end{remark}

\section{Energy estimates for solutions to porous medium systems}\label{sec:app-energy}

In the following, unless otherwise stated, summation over repeated indices is assumed. Moreover, we use the one-sided cylinders introduced in~\eqref{cy-infty}. 

In this appendix, we collect two types of energy estimates, the first of which is designed for Moser's iteration.
\begin{lemma}[energy estimate I]\label{lem:energy-est}
Let $m>0$. There exist constants $C_1=C_1(\nu,L)$
and $C_2=C_2(\nu)$, such that whenever $u$ is a weak solution in $\Omega_T$ to \eqref{por-med-eq}
under the assumptions \eqref{growth},  $\Phi\colon
\R_{\ge 0}\to \R_{\ge 0}$ is a non-negative, bounded, non-decreasing, Lipschitz continuous function that satisfies
\begin{equation}\label{condition-Phi}
    C_\Phi:= \sup_{\tau\in \R_{\ge 0}\cap \{ \Phi>0\} }
    \frac{\tau\Phi'(\tau)}{\Phi (\tau)}<\infty,
\end{equation}
and $Q_{R,S}(z_o)\subset \Omega_T$, then for any cut-off function $\zeta\in C^\infty\big(Q_{R,S}(z_o),[0,1]\big)$ that vanishes on the parabolic boundary $\partial_P Q_{R,S}(z_o)$ of $Q_{R,S}(z_o)$, we have 
\begin{align}\label{En-Est-Phi}\nonumber
    \tfrac1{m+1}&\sup_{\tau\in (t_o-S,t_o]}\int_{B_R(x_o)\times\{\tau\}}  v\, \zeta^2\,\dx 
    + \tfrac12\nu 
    \iint_{Q_{R,S}(z_o)} \left|D\power{u}{m}\right|^2 \Phi(\left|u\right|^{2m}) \zeta^2 \,\dx\dt \\\nonumber
    & \le
    C_1  \iint_{Q_{R,S}(z_o)}|\power{u}{m}|^{2 } \Phi(|\power{u}{m}|^{2})|\nabla \zeta|^2\,\dx\dt
    +
    \tfrac2{m+1}\iint_{Q_{R,S}(z_o)} v\, |\partial_{t}\zeta^2|\,\dx\dt\\
    &\phantom{\le\,}
    +C_2(1+C_\Phi)\iint_{Q_{R,S}(z_o)} |F|^2\Big[ \Phi(|\power{u}{m}|^{2})
    +|\power{u}{m}|^2\Phi^{\prime}(|\power{u}{m}|^{2})\Big]\zeta^2\,\dx\dt,
\end{align}
where
\begin{equation*}
    v:=\int_0^{|u|^{1+m}} \Phi\big(s^{\frac{2 m}{m+1}}\big)\, \ds.
\end{equation*}
\end{lemma}
\begin{proof}
Following the proof of \cite[Lemma~3.1]{BDS-sing}, for $w\in L^1(\Omega_T,\R^k)$, we define the following mollification in time
\[
\llbracket w\rrbracket_h(x,t):=\frac1h \int_0^t e^{\frac{s-t}h}w(x,s)\,\dx\ds,
\]
where $h>0$. 
Since we can assume without loss of generality that $u\in C^0([0,T];L^2_{\loc}(\Omega,\R^k))$, from \eqref{weak-solution} we deduce its mollified version
\begin{equation}\label{Eq-Sys_Moll}
\begin{aligned}
    \iint_{ \Omega_T} & 
    \big[\partial_t\llbracket u\rrbracket_h \cdot \varphi+
    \llbracket{\boldsymbol A}(x,t,  u,D\power{u}{m})\rrbracket_h \cdot D \varphi\big] \,\dx\dt\\
    &=
    -\iint_{\Omega_T} \llbracket F\rrbracket_h \cdot D\varphi\,\dx\dt+\frac1h\int_\Omega u(0)\cdot\int_0^T e^{-\frac sh}\varphi\,\ds\dx, 
\end{aligned}
\end{equation}
for any $\varphi\in L^2(0,T;W^{1,2}_0(\Omega,\R^k))$ and $h>0$. 

Given $z_o\in\Omega_T$ and a cylinder $Q_{R,S}(z_o)\Subset\Omega_T$, fix a time $\tau\in (t_o-S,t_o)$ and $\delta \in (0,t_o-\tau)$, let $\chi(t)$ be a function which is equal to $1$ for $t\in (t_o-S,\tau)$, vanishes for $t\in (\tau+\delta, t_o)$, and is linearly interpolated for $t\in [\tau,\tau +\delta]$, and $\zeta$ as in the statement of the lemma.
%=B_R(x_o)\times (t_o-S,t_o]$. 
To abbreviate and simplify the notation, we write $Q_{R,S}$ instead of $Q_{R,S}(z_o)$, $B_R$ instead of $B_R(x_o)$. We choose
\begin{equation*}
   \varphi_h
   = 
   \power{u}{m}\Phi(|\power{\llbracket u\rrbracket_h}{m}|^{2 }) \zeta^2 \chi.
\end{equation*} 
as testing function in \eqref{Eq-Sys_Moll}. 

For the part in \eqref{Eq-Sys_Moll} containing the time derivative, we get 
\begin{equation}\label{Eq:Sys:1}
\begin{aligned}
    \iint_{\Omega_T} \partial_t\llbracket u\rrbracket_h \cdot \varphi_h \,\dx\dt
    &=
    \iint_{Q_{R,S}} \partial_t\llbracket u\rrbracket_h \cdot \power{u}{m} \Phi(|\power{\llbracket u\rrbracket_h}{m}|^{2 }) \zeta^2 \chi \,\dx\dt\\
    &=
    \iint_{Q_{R,S}} \partial_t\llbracket u\rrbracket_h \cdot (\power{u}{m}- \power{\llbracket u\rrbracket_h}{m}) \Phi(|\power{\llbracket u\rrbracket_h}{m}|^{2 }) \zeta^2 \chi \,\dx\dt\\
    &\quad +
    \iint_{Q_{R,S}} \partial_t\llbracket u\rrbracket_h \cdot \power{\llbracket u\rrbracket_h}{m}\Phi(|\power{\llbracket u\rrbracket_h}{m}|^{2 }) \zeta^2 \chi \,\dx\dt\\
    &\ge
    \iint_{Q_{R,S}} \partial_t\llbracket u\rrbracket_h \cdot \power{\llbracket u\rrbracket_h}{m}\Phi(|\power{\llbracket u\rrbracket_h}{m}|^{2 }) \zeta^2 \chi \,\dx\dt,
\end{aligned}
\end{equation}
where we took into account that
\begin{equation}\label{moly-sign}
\partial_t\llbracket u\rrbracket_h=\frac1h(u-\llbracket u\rrbracket_h)\quad\Rightarrow\quad\partial_t\llbracket u\rrbracket_h \cdot (\power{u}{m}- \power{\llbracket u\rrbracket_h}{m})\ge0,
\end{equation}
due to the monotonicity of the map $u\mapsto\power{u}{m}$.
Note that
\begin{align*}
    \partial_t |\llbracket u\rrbracket_h|^{m+1}
    = 
    (m+1) |\llbracket u\rrbracket_h|^{m-1}\llbracket u\rrbracket_h\cdot\partial_t \llbracket u\rrbracket_h
    = 
    (m+1) \power{\llbracket u\rrbracket_h}{m} \cdot\partial_t \llbracket u\rrbracket_h .
\end{align*}
If this is substituted back into \eqref{Eq:Sys:1}, the result is
\begin{align*}
    \iint_{\Omega_T} \partial_t\llbracket u\rrbracket_h \cdot \varphi_h \,\dx\dt
    &\ge
    \tfrac{1}{m+1} \iint_{Q_{R,S}} \partial_t|\llbracket u\rrbracket_h|^{m+1} \Phi\Big(|\llbracket u\rrbracket_h|^{\frac{2 m}{m+1}(m+1)}\Big) \zeta^2 \chi \,\dx\dt \\
    &=
    \tfrac{1}{m+1} \iint_{Q_{R,S}} \partial_t \bigg[\int_0^{|\llbracket u\rrbracket_h|^{m+1}} \Phi\big(s^{\frac{2 m}{m+1}}\big) \,\ds\bigg] \zeta^2 \chi \,\dx\dt \\
    &= 
    -\tfrac{1}{m+1} \iint_{Q_{R,S}} \bigg[\int_0^{|\llbracket u\rrbracket_h|^{m+1}} \Phi\big(s^{\frac{2 m}{m+1}}\big) \,\ds\bigg] \,\partial_t(\zeta^2 \chi) \,\dx\dt.
\end{align*}
We may then pass to the limit $h\downarrow0$, and conclude that 
\[
    \liminf_{h\downarrow0}\iint_{\Omega_T} \partial_t\llbracket u\rrbracket_h \cdot \varphi_h \,\dx\dt
    \ge
    -\tfrac1{m+1}\iint_{Q_{R,S}} v \,\partial_t(\zeta^2 \chi) \,\dx\dt.
\]

Next, we consider the diffusion term in \eqref{Eq-Sys_Moll}. After passing to the limit $h\downarrow0$, in order to estimate the contribution of such a term in the weak formulation, we first calculate the spatial derivative of the testing function $\varphi=\power{u}{m} \Phi(|\power{u}{m}|^{2}) \zeta^2 \chi$. For simplicity, we write ${\boldsymbol A}_\alpha$, $\alpha\in\{1,\dots,N\}$ instead of ${\boldsymbol A}_\alpha (x,t,u,D\power{u}{m})$, and we obtain
\begin{align*}
    D_\alpha\varphi
    =
    \big[D_\alpha \power{u}{m} \Phi(|\power{u}{m}|^{2}) 
    +2\power{u}{m}  
    (D_\alpha \power{u}{m} \cdot \power{u}{m})
    \Phi^{\prime}(|\power{u}{m}|^{2}) \big]\zeta^2 \chi
    +
    \power{u}{m} \Phi(|\power{u}{m}|^{2 }) D_\alpha \zeta^2\, \chi ,
 \end{align*}
 which implies that 
 \begin{align*}
    {\boldsymbol A}_\alpha \cdot D_\alpha\varphi
    &=
    \big[{\boldsymbol A}_\alpha\cdot D_\alpha \power{u}{m}\, \Phi(|\power{u}{m}|^{2}) +
    2 ({\boldsymbol A}_\alpha \cdot \power{u}{m}) (D_\alpha \power{u}{m} \cdot \power{u}{m}) \Phi^\prime(|\power{u}{m}|^{2 }) \big]\zeta^2 \chi \\
    &\quad+
    {\boldsymbol A}_\alpha \cdot \power{u}{m} \Phi(|\power{u}{m}|^{2}) D_\alpha \zeta^2\, \chi \\
    &\ge 
    {\boldsymbol A}_\alpha\cdot D_\alpha \power{u}{m}\, \Phi(|\power{u}{m}|^{2}) \zeta^2 \chi +
    {\boldsymbol A}_\alpha \cdot \power{u}{m} \Phi(|\power{u}{m}|^{2}) D_\alpha \zeta^2\, \chi,
\end{align*}
where we used hypothesis \eqref{growth}$_3$ and the fact that $\Phi'\ge 0$ almost everywhere on $\R_{\ge 0}$. 
As for the right-hand side, once more we pass to the limit $h\downarrow0$, and we obtain
\begin{align*}
    F_\alpha\cdot D_\alpha\varphi
    &=
    \chi\Big[D_\alpha\power{u}{m}\Phi(|\power{u}{m}|^{2}) \zeta^2 +
    2\power{u}{m} 
    (D_\alpha \power{u}{m} \cdot \power{u}{m})
    \Phi^\prime(|\power{u}{m}|^{2 })\zeta^2 \\
    &\qquad\qquad\qquad\qquad\qquad\qquad\qquad\qquad 
    +
    \power{u}{m}\Phi(|\power{u}{m}|^{2 }) D_\alpha\zeta^2\Big]\cdot  F_\alpha.
\end{align*}
Substituting everything back into \eqref{Eq-Sys_Moll}, we get
\begin{align}\label{Eq:pre-energy}\nonumber
    -\tfrac{1}{m+1}
    & 
    \iint_{Q_{R,S}} v\, \partial_t(\zeta^2 \chi)\,\dx\dt  +
    \iint_{Q_{R,S}} 
    \boldsymbol{A}_\alpha \cdot
    D_\alpha \power{u}{m}\,\Phi(|\power{u}{m}|^{2}) \zeta^2 \chi\,\dx\dt \\\nonumber
    &\le
    -2 \iint_{Q_{R,S}} {\boldsymbol A} \cdot (\power{u}{m} \otimes \nabla\zeta)\, \Phi(|\power{u}{m}|^{2})\zeta \chi \,\dx\dt\\\nonumber
    &\quad-
    \iint_{Q_{R,S}} D_\alpha\power{u}{m} \cdot F_\alpha\,\Phi(|\power{u}{m}|^{2})\zeta^2\chi \,\dx\dt\\\nonumber
    &\quad-
    2 \iint_{Q_{R,S}} (\power{u}{m}\cdot F_\alpha)(D_\alpha \power{u}{m} \cdot \power{u}{m})\,\Phi^{\prime}(|\power{u}{m}|^{2})\zeta^2\chi \,\dx\dt\\\nonumber
    &\quad-
    2\iint_{Q_{R,S}} F \cdot(\power{u}{m}
    \otimes \nabla\zeta)\, \Phi(|\power{u}{m}|^{2}) \zeta\chi \,\dx\dt\\
    &=:
    \mathbf{I}_1+\mathbf{I}_2+\mathbf{I}_3+\mathbf{I}_4.
\end{align}
In the following we deal with the individual terms one by one. We start with the second integral  on the left side. Using the lower bound 
\eqref{growth}$_1$, we get
\begin{align*}
    \nu \iint_{Q_{R,S}}  |D\power{u}{m}|^{2}\Phi(|\power{u}{m}|^{2}) 
    \zeta^2 \chi \,\dx\dt
    \le
   \iint_{Q_{R,S}} 
   \boldsymbol{A}_\alpha \cdot
    D_\alpha \power{u}{m}\Phi(|\power{u}{m}|^{2}) \zeta^2 \chi\,\dx\dt.
\end{align*}
The term $\mathbf{I}_1$ is estimated by \eqref{growth}$_2$ and Young's inequality. This gives
\begin{align*}
    \mathbf{I}_1
    &\le 
    2L \iint_{Q_{R,S}} |D\power{u}{m}| |\power{u}{m}|  |\nabla\zeta|\Phi(|\power{u}{m}|^{2}) \zeta \chi \,\dx\dt \\
    &\le 
    \tfrac{\nu}6
    \iint_{Q_{R,S}}\left|D \power{u}{m}\right|^2 \Phi(|\power{u}{m}|^{2}) \zeta^2 \chi \,\dx\dt +
    \tfrac{36 L^2}{\nu} \iint_{Q_{R,S}}   |\power{u}{m}|^{2}  \Phi(|\power{u}{m}|^{2})\,  |\nabla\zeta|^2\chi  \,\dx\dt.
\end{align*}
The remaining terms $\mathbf{I}_2$, $\mathbf{I}_3$, and $\mathbf{I}_4$
are estimated by Young's inequality. 
Namely, for the integral $\mathbf{I}_2$ we have
\begin{align*}
    \mathbf{I}_2
    &\le 
    \tfrac16 \nu\iint_{Q_{R,S}} |D \power{u}{m}|^2 \Phi(|\power{u}{m}|^{2}) \zeta^2 \chi \,\dx\dt +
    \tfrac{3}{2\nu}\iint_{Q_{R,S}} \Phi(|\power{u}{m}|^{2}) |F|^2\zeta^2\chi \,\dx\dt.
\end{align*}
For the third term we use Young's inequality in the form
\begin{align*}
    2(\power{u}{m}&\cdot F_\alpha)
    (D_\alpha \power{u}{m}\cdot \power{u}{m})\Phi^\prime( |\power{u}{m}|^2)\\
    &\le
    2|D \power{u}{m}|
    \Phi^\frac12( |\power{u}{m}|^2)
    |F| \frac{ |\power{u}{m}|^2\Phi^\prime( |\power{u}{m}|^2)}
    {\Phi^\frac12( |\power{u}{m}|^2)}\\
    &\le
    \eps
    | D \power{u}{m}|^2
    \Phi( |\power{u}{m}|^2) 
    +
    \eps^{-1} 
    |F|^2 |\power{u}{m}|^2\Phi^\prime( |\power{u}{m}|^2)
    \underbrace{\frac{|\power{u}{m}|^2\Phi^\prime( |\power{u}{m}|^2)}
    {\Phi( |\power{u}{m}|^2)}}_{\le\, C_\Phi}
    \\
    &\le
    \eps
    | D \power{u}{m}|^2
    \Phi( |\power{u}{m}|^2) 
    +
    \frac{C_\Phi}{\eps} |F|^2 |\power{u}{m}|^2\Phi^\prime( |\power{u}{m}|^2).
\end{align*}
With $\eps=\frac16\nu$ this leads to
\begin{align*}
    \mathbf{I}_3
    &\le 
    \tfrac16\nu \iint_{Q_{R,S}} 
    | D \power{u}{m}|^2
    \Phi( |\power{u}{m}|^2) 
    \zeta^2\chi\,\dx\dt\\
    &\phantom{\le\,}
    +
    \tfrac{6C_\Phi}{\nu}
    \iint_{Q_{R,S}} |F|^2|\power{u}{m}|^{2}\Phi^{\prime}(|\power{u}{m}|^{2}) \zeta^2\chi \,\dx\dt.
\end{align*}
Finally, the last integral $\mathbf{I}_4$ is estimated by
\begin{align*}
    \mathbf{I}_4 
    &\le  
    \iint_{Q_{R,S}} |\power{u}{m}|^{2} |\nabla\zeta|^2 \Phi(|\power{u}{m}|^{2}) \chi\,\dx\dt
    + 
    \iint_{Q_{R,S}} |F|^2 \Phi(|\power{u}{m}|^{2 }) \zeta^2\chi \,\dx\dt.
\end{align*}
Combining these estimates with \eqref{Eq:pre-energy}, we arrive at
\begin{align*}
    -\tfrac{1}{m+1}
    \iint_{Q_{R,S}}& v\,\zeta^2 \partial_t\chi\,\dx\dt
    +
    \tfrac12 \nu
    \iint_{Q_{R,S}} \left|D \power{u}{m}\right|^2 \Phi(|\power{u}{m}|^{2}) \zeta^2 \chi\,\dx\dt \\
    &\le 
     \big(1+\tfrac{36 L^2}{\nu}\big) \iint_{Q_{R,S}}
    |\power{u}{m}|^{2 }\Phi(|\power{u}{m}|^{2 })|\nabla \zeta|^2 \chi \,\dx\dt\\
    &\phantom{\le\,}
    +
    \big(1+\tfrac{3}{2\nu}\big)
    \iint_{Q_{R,S}} |F|^2\Phi(|\power{u}{m}|^{2}) \zeta^2\chi\,\dx\dt
    \\
    &\phantom{\le\,}
    +
    \tfrac{6C_\Phi}{\nu}
    \iint_{Q_{R,S}} |F|^2|\power{u}{m}|^{2} \Phi^{\prime}(|\power{u}{m}|^{2})\zeta^2\chi\,\dx\dt\\
    &\phantom{\le\,}
    +
    \tfrac{1}{m+1}
    \iint_{Q_{R,S}} v\, \chi\partial_t\zeta^2 \,\dx\dt.
\end{align*}
Arrived at this point, we let $\delta\downarrow 0$ in $\chi\equiv \chi_\delta$. In the limit we have
$\chi\downarrow\mathbf 1_{(t_o-S,\tau]}$. The first integral on the left-hand side
converges to $ \tfrac1{m+1}\int_{B_{R}\times\{\tau\}}v\zeta^2\,\dx $,  while the second one can be estimated from below by $\frac12\nu \iint_{B_R\times (t_o-S,\tau]}\zeta^{2} |D\power{u}{m}|^2 \Phi(|\power{u}{m}|^2)\,\dx\dt$. In the right-hand side integrals we bound the cut-off function $\chi$ by 1. In this way, we obtain
\begin{align*}
    \tfrac{1}{m+1}
    \int_{B_{R}\times\{\tau\}}& v\zeta^2\,\dx
    +
    \tfrac12 \nu
    \iint_{B_{R}\times (t_o-S,\tau]} |D \power{u}{m}|^2 \Phi(|\power{u}{m}|^{2}) \zeta^2 \,\dx\dt \\
    &\le 
     \big(1+\tfrac{36 L^2}{\nu}\big) \iint_{Q_{R,S}}
     |\power{u}{m}|^{2 }\Phi(|\power{u}{m}|^{2 })|\nabla \zeta|^2 \,\dx\dt\\
    &\phantom{\le\,}
    +
    \big(1+\tfrac{3}{2\nu}\big)
    \iint_{Q_{R,S}} |F|^2\Phi(|\power{u}{m}|^{2}) \zeta^2\,\dx\dt
    \\
    &\phantom{\le\,}
    +
    \tfrac{6C_\Phi}{\nu}\iint_{Q_{R,S}}|\power{u}{m}|^{2} |F|^2\Phi^{\prime}(|\power{u}{m}|^{2})\zeta^2\,\dx\dt\\
    &\phantom{\le\,}
    +
    \tfrac{1}{m+1}
    \iint_{Q_{R,S}} v\, |\partial_t\zeta^2| \,\dx\dt.
\end{align*}
After that, we take the supremum over $\tau\in (t_o-S,t_o)$ in the first integral on the left and let $\tau\uparrow t_o$ in the second one to deduce the claimed energy inequality. 
\end{proof}

In particular,  for the function $\Phi:=\phi_{\alpha,k,\ell}$ defined 
in \eqref{def:phi-alpha-l-k} in the proof of the qualitative $L^\infty$-bound of Proposition \ref{prop:sup-est-qual}, we have $C_{\phi_{\alpha,k,\ell}}=\alpha$.  Therefore, we  can use the previously shown energy estimate and run a Moser iteration scheme obtaining a \emph{qualitative} local $L^\infty$-bound.

Now we present the second energy estimate which will be used in De Giorgi's iteration.
\begin{lemma}[energy estimate II]\label{lem:new-energy-est}
Let $m>0$ and $F\in L^2_{\rm loc}(\Omega_T,\R^{kN})$. There exists a constant $C=C(m,\nu,L)$, such that whenever $u$ is a locally bounded weak solution in $\Omega_T$ to \eqref{por-med-eq}
under the assumptions \eqref{growth}, then for any cut-off function $\zeta\in C^\infty\big(Q_{R,S}(z_o),[0,1]\big)$ that vanishes on the parabolic boundary $\partial_P Q_{R,S}(z_o)$ of $Q_{R,S}(z_o)$ and any $k\ge 0$ we have
\begin{align*}
    &  \sup_{\tau\in (t_o-S,t_o)} 
    \int_{B_R \times\{\tau\}} 
    \big(|u|^m-
    |k|^m\big)_+^{3+\frac{1}{m}}
    \, \zeta^2\,\dx
    +
    \iint_{Q_{R, S}}\big|\nabla\big(|u|^m-k^m\big)_{+}^2\big|^2 \zeta^2\,\dx\dt\\
    &\qquad\le 
    C\iint_{Q_{R,S}}|u|^{4m}\mathbf{1}_{\{|u|>k\}}|\nabla \zeta|^2
    \,\dx\dt +
    C \iint_{Q_{R,S}} 
    |u|
    \big(|u|^{m}-k^m\big)_+^3
    |\partial_t \zeta^2|\,\dx\dt \\
    &\qquad\phantom{\le\,}
    +
    C\iint_{Q_{R,S}}|F|^2|u|^{2m}\mathbf{1}_{\{|u|>k\}}
    \zeta^2 \,\dx\dt.
\end{align*}
\end{lemma}

\begin{proof}

Let us refer back to the mollified equation \eqref{Eq-Sys_Moll}. Given $z_o\in\Omega_T$ and a cylinder $Q_{R,S}(z_o)\Subset\Omega_T$, fix a time $\tau\in (t_o-S,t_o)$ and $\delta \in (0,t_o-\tau)$, let $\chi(t)$ be a function which is equal to $1$ for $t\in (t_o-S,\tau)$, vanishes for $t\in (\tau+\delta, t_o)$, and is linearly interpolated for $t\in [\tau,\tau +\delta]$, and $\zeta$ as in the statement of the lemma. To abbreviate and simplify the notation, we write $Q_{R,S}$ instead of $Q_{R,S}(z_o)$, $B_R$ instead of $B_R(x_o)$. Moreover, we assume $u$ to be locally bounded (this is a consequence of the qualitative $L^\infty$-estimate obtained with Moser's technique in Proposition~\ref{prop:sup-est-qual}), and we choose in the mollified equation \eqref{Eq-Sys_Moll} as testing function
\[
    \varphi
    =
    \power{u}{m} \left(|\power{\llbracket u\rrbracket_h}{m}|-k^m\right)_{+}^2 \zeta^2 \chi.
\] 
Let us first concentrate on the diffusion term, and then on the right-hand side, whereas the time derivative will be dealt with at last. After passing to the limit $h\downarrow0$ we compute the spatial derivative 
\begin{align*}
    %D_\alpha \varphi
    %&=
    D_\alpha&\left[\power{u}{m}\big(|\power{u}{m}|-k^m\big)_{+}^2 \zeta^{2} \chi\right] \\
    &=
    D_\alpha \power{u}{m}\big(|\power{u}{m}|-k^m\big)_{+}^2 
    \zeta^2 \chi
    +
    2 \power{u}{m}\big(\left|\power{u}{m}\right|-k^m\big)_+
    \frac{D_\alpha\power{u}{m}\cdot\power{u}{m}}{|\power{u}{m}|} \zeta^2 \chi \\
    &\phantom{=\,}
    +2\power{u}{m}\big(|\power{u}{m}|-k^m\big)^2  \zeta D_\alpha \zeta \chi.
\end{align*}
Hence, the diffusion term becomes
$
    %\iint_{Q_{R,S}}
    %\boldsymbol{A}_\alpha(x,t,u,D\power{u}{m}) \cdot D_\alpha \varphi
    %\, \dx\dt
    %=: $
    \boldsymbol D_1+\boldsymbol D_2+\boldsymbol D_3
$,
where
\begin{align*}
    \boldsymbol D_1
    &:=
    \iint_{Q_{R,S}} 
    \boldsymbol{A}_\alpha (x,t,u,D\power{u}{m})
    \cdot D_\alpha \power{u}{m}
    \big(|\power{u}{m}|-k^m\big)_+^2 \zeta^2 \chi\,\dx\dt,\\
    \boldsymbol D_2
    &:=
    2 \iint_{Q_{R,S}} 
    \boldsymbol{A}_\alpha (x,t,u,D\power{u}{m})
    \cdot 
    \power{u}{m} \frac{D_\alpha \power{u}{m} \cdot \power{u}{m}}{|\power{u}{m}|}
    \big(\left|\power{u}{m}\right|-k^m\big)_+\zeta^2 \chi\,\dx\dt,\\
    \boldsymbol D_3
    &:=
    2\iint_{Q_{R,S}} \boldsymbol{A}_\alpha (x,t,u,D\power{u}{m})
    \cdot 
    \power{u}{m}\big(\left|\power{u}{m}\right|-k^m\big)_+^2\zeta D_\alpha\zeta \chi\,\dx\dt.
\end{align*}
We proceed and estimate the three terms. Using the coercivity 
\eqref{growth}$_1$ we have
\begin{align*}
    \boldsymbol D_1 
    &\ge 
    \nu \iint_{Q_{R,S}}\left|D \power{u}{m}\right|^2\left(\left|\power{u}{m}\right|-k^m\right)_{+}^2 \zeta^2\chi\,\dx\dt.
\end{align*}
On the other hand, \eqref{growth}$_3$ immediately leads to $\boldsymbol D_2\ge 0$, and no further estimate is needed for $\boldsymbol D_2$. Finally, using the growth assumption \eqref{growth}$_2$ and Young's inequality we obtain
\begin{align*}
    \boldsymbol D_3
    &\le
    2 L \iint_{Q_{R,S}}
    |D \power{u}{m}||\power{u}{m}|
    \big(|\power{u}{m}|-k^m\big)_{+}^2|\nabla \zeta| \zeta \chi\,\dx\dt\\
    &\le 
    \varepsilon \iint_{Q_{R,S}}
    |D \power{u}{m}|^2\big(|\power{u}{m}|-k^m\big)_+^2 \zeta^2 \chi\,\dx\dt\\
    &\phantom{\le\,}
    +\frac{L^2}{\varepsilon} 
    \iint_{Q_{R,S}}
    |\power{u}{m}|^2\big(|\power{u}{m}|-k^m\big)_{+}^2 
    |\nabla \zeta|^2 \chi\,\dx\dt.
\end{align*}

On the right-hand side of the mollified equation \eqref{Eq-Sys_Moll}, the last term vanishes when $h\downarrow0$ thanks to $\vp(0)=0$.
Let us now consider the contributions coming from the inhomogeneity $F$. Likewise, after passing to the limit $h\downarrow0$ we have
\begin{equation*}
    \iint_{Q_{R,S}} F_\alpha\cdot D_\alpha\varphi\,\dx\dt
    =:\boldsymbol R_1+\boldsymbol R_2+\boldsymbol R_3,
\end{equation*}
where
\begin{align*}
    \boldsymbol R_1
    &:=
     \iint_{Q_{R,S}}
     F_\alpha\cdot D_\alpha \power{u}{m}
     \big(|\power{u}{m}|-k^m\big)^2_+\zeta^2 \chi\,\dx\dt,\\
     \boldsymbol R_2
    &:=
     2\iint_{Q_{R,S}}
     F_\alpha \cdot \power{u}{m} \frac{D_\alpha \power{u}{m} \cdot \power{u}{m}}{|\power{u}{m}|}
    \big(|\power{u}{m}|-k^m\big)_+\zeta^2\chi\,\dx\dt,\\
    \boldsymbol R_3
    &:=
    2 \iint_{Q_{R,S}} 
    F_\alpha \cdot \power{u}{m}\big(|\power{u}{m}|-k^m\big)_+^2
    \zeta D_\alpha \zeta \chi\,\dx\dt.
\end{align*}
In  turn we estimate $\boldsymbol R_i$, $i=1,2,3$. By Young's inequality we have
\begin{align*}
    \boldsymbol R_1
    &\le
    \varepsilon 
    \iint_{Q_{R;S}}
    |D \power{u}{m}|^2
    \big(|\power{u}{m}|-k^m\big)_+^2 \zeta^2 \chi\,\dx\dt\\
    &\phantom{\le\,}
    +\tfrac1{4\varepsilon} \iint_{Q_{R,S}}|F|^2
    \big(|\power{u}{m}|-k^m\big)_+^2 \zeta^2 \chi\,\dx\dt,
\end{align*}
and similarly
\begin{align*}
    \boldsymbol R_2 
    & \le 2 \iint_{Q_{R,S}}|F||\power{u}{m}||D \power{u}{m}|
    \big(|\power{u}{m}|-k^m\big)_{+} \zeta^2 \chi\,\dx\dt \\
    &\le 
    \varepsilon \iint_{Q_{R,S}}\!
    |D \power{u}{m}|^2
    \big(|\power{u}{m}|-k^m\big)_{+}^2 \zeta^2 \chi\,\dx\dt
    +
    \tfrac{1}{\varepsilon} \iint_{Q_{R,S}}\!|F|^2|\power{u}{m}|^2 \mathbf 1_{\{ |u|>k\}}\zeta^2 \chi\,\dx\dt.
\end{align*}
Finally, using once again Young's inequality we have
\begin{align*}
    \boldsymbol R_3 
    & \le 
    2 \iint_{Q_{R,S}}
    |F||\power{u}{m}|
    \big(|\power{u}{m}|- k^m\big)_+^2 
    \zeta| \nabla\zeta|\chi
    \,\dx\dt \\
    &\le \iint_{Q_{R,S}}
    |F|^2|\power{u}{m}|^2 
    \mathbf 1_{\{|u|>k\}}
    \zeta^2 \chi\,\dx\dt
    +
    \iint_{Q_{R,S}}
    \big(|\power{u}{m}|-k^m\big)_{+}^4
    |\nabla \zeta|^2 \chi\,\dx\dt.
\end{align*}
When estimating $\boldsymbol R_2$
and $\boldsymbol R_3$, we also took advantage of the fact that the domain of integration is  $Q_{R,S}
\cap \{ |u|>k\}$.  We conclude with the estimate of the term containing the time derivative, namely
\begin{align*}
    \iint_{Q_{R,S}} &
    \partial_t\llbracket u\rrbracket_h\cdot \power{u}{m}\big(|\power{\llbracket u\rrbracket_h}{m}|-k^m\big)_+^2\zeta^2\chi\,\dx\dt\\
    &=\iint_{Q_{R,S}}\partial_t\llbracket u\rrbracket_h\cdot \big(\power{u}{m}-\power{\llbracket u\rrbracket_h}{m}\big)\big(|\power{\llbracket u\rrbracket_h}{m}|-k^m\big)_+^2\zeta^2\chi\,\dx\dt\\
    &\quad +\iint_{Q_{R,S}}\partial_t\llbracket u\rrbracket_h\cdot \power{\llbracket u\rrbracket_h}{m}\big(|\power{\llbracket u\rrbracket_h}{m}|-k^m\big)_+^2\zeta^2\chi\,\dx\dt\\
    &\ge\iint_{Q_{R,S}}\partial_t\llbracket u\rrbracket_h\cdot \power{\llbracket u\rrbracket_h}{m}\big(|\power{\llbracket u\rrbracket_h}{m}|-k^m\big)_+^2\zeta^2\chi\,\dx\dt,
\end{align*}
where we used~\eqref{moly-sign}. 
Moreover, we have
\begin{align*}
    \partial_t \llbracket u\rrbracket_h\cdot \power{\llbracket u\rrbracket_h}{m}\big(|\power{\llbracket u\rrbracket_h}{m}|-k^m\big)_+^2
    &=
    \tfrac{1}{m+1} \partial_t|\llbracket u\rrbracket_h|^{m+1}\big(|\llbracket u\rrbracket_h|^m-k^m\big)_{+}^2  \\
    &=
    \tfrac{1}{m+1} \partial_t|\llbracket u\rrbracket_h|^{m+1}
    \Big(\big(|\llbracket u\rrbracket_h|^{m+1}\big)^{\frac{m}{m+1}}-k^{m}\Big)_{+}^2 \\
    &=
    \tfrac{1}{m+1}
    \partial_t \bigg[
    \int_{k^{m+1}}^{|\llbracket u\rrbracket_h|^{m+1}} \Phi\big(s^{\frac{m}{m+1}}\big) 
    \ds\bigg].
\end{align*}
To obtain the last equality we introduced $\Phi(s):=(s-k^m)_{+}^2$. Hence,
\begin{align*}
    \iint_{Q_{R,S}}
    &\partial_t \llbracket u\rrbracket_h\cdot \power{\llbracket u\rrbracket_h}{m}\big(|\power{\llbracket u\rrbracket_h}{m}|-k^m\big)_+^2\zeta^2\chi\,\dx\dt\\
    &=
    \tfrac{1}{m+1} \iint_{Q_{R,S}}\bigg[\partial_t \int_{k^{m+1}}^{|\llbracket u\rrbracket_h|^{m+1}} \Phi\big(s^{\frac{m}{m+1}}\big) \ds\bigg] \zeta^2 \chi\,\dx\dt \\
    &=
    -\tfrac{1}{m+1} \iint_{Q_{R,S}}
    \bigg[\int_{k^{m+1}}^{|\llbracket u\rrbracket_h|^{m+1}} \Phi\big(s^{\frac{m}{m+1}}\big) \ds\bigg] \partial_t (\zeta^2 \chi)\,\dx\dt \\
    &\to
    -\tfrac{1}{m+1} \iint_{Q_{R,S}} v\, \zeta^2 \partial_t \chi\,\dx\dt
    -\tfrac{1}{m+1} \iint_{Q_{R,S}} v\, \chi\, \partial_t \zeta^2\,\dx\dt
\end{align*}
as $h\downarrow0$, where we abbreviated
\begin{equation*}
    v:=\int_{k^{m+1}}^{|u|^{m+1}} \Phi\big(s^{\frac{m}{m+1}}\big)\, \ds.
\end{equation*}
Collecting and re-arranging all the terms yields
\begin{align*}
    - & \tfrac{1}{m+1} \iint_{Q_{R,S}} v\, \zeta^2 \partial_t \chi\,\dx\dt
    +
    (\nu-3\varepsilon) 
    \iint_{Q_{R,S}}|D \power{u}{m}|^2
    \big(|\power{u}{m}|-k^m\big)_+^2 \zeta^2\chi \,\dx\dt \\
    &\le
    \tfrac{L^2}{\varepsilon} 
    \iint_{Q_{R,S}}\underbrace{|\power{u}{m}|^2
    \big(|\power{u}{m}|-k^m\big)_{+}^2}_{\le |u|^{4m}\mathbf{1}_{\{|u|>k\}}} |\nabla \zeta|^2\chi 
    \,\dx\dt +
    \iint_{Q_{R,S}}
    \underbrace{\big(|\power{u}{m}|-k^m\big)_{+}^4}_{\le |u|^{4m}\mathbf{1}_{\{|u|>k\}}}
    |\nabla \zeta|^2 \chi
    \,\dx\dt \\
    &\phantom{\le\,}
    +\tfrac{1}{4\varepsilon}
    \iint_{Q_{R,S}}|F|^2\underbrace{\big(|\power{u}{m}|-k^m\big)_+^2}_{\le |u|^{2m}\mathbf{1}_{\{|u|>k\}}}
    \zeta^2\chi \,\dx\dt\\
    &\phantom{\le\,}
    +
    \big( 1+\tfrac1{\varepsilon}\big) 
    \iint_{Q_{R,S}}
    |F|^2|\power{u}{m}|^2 \mathbf{1}_{\{|u|>k\}}
    \zeta^2 \chi\,\dx\dt \\
    &\phantom{\le\,}
    +
    \tfrac{1}{m+1} \iint_{Q_{R,S}} v|\partial_t \zeta^2|\chi\,\dx\dt\\
    &\le
    \big( 1+\tfrac{L^2}{\varepsilon} \big)
    \iint_{Q_{R,S}}|u|^{4m}\mathbf{1}_{\{|u|>k\}} |\nabla \zeta|^2\chi 
    \,\dx\dt\\
    &\phantom{\le\,}
    +
    \big(1+\tfrac{5}{4\varepsilon}\big)
    \iint_{Q_{R,S}}|F|^2|u|^{2m}\mathbf{1}_{\{|u|>k\}}
    \zeta^2\chi \,\dx\dt
    +
    \tfrac{1}{m+1} \iint_{Q_{R,S}} v|\partial_t \zeta^2|\chi\,\dx\dt.
\end{align*}
Choosing $\varepsilon =\frac16 \nu$ we end up with
\begin{align*}
    -\tfrac{1}{m+1} &\iint_{Q_{R,S}} v\, \zeta^2 \partial_t \chi\,\dx\dt
    +
    \tfrac12 \nu
    \iint_{Q_{R,S}}|D \power{u}{m}|^2
    \big(|\power{u}{m}|-k^m\big)_+^2 \zeta^2\chi 
    \,\dx\dt\\
    &\le 
    C\iint_{Q_{R,S}}|u|^{4m}\mathbf{1}_{\{|u|>k\}} |\nabla \zeta|^2\chi 
    \,\dx\dt +
    \tfrac{1}{m+1} \iint_{Q_{R,S}} v\chi|\partial_t \zeta^2|\,\dx\dt\\
    &\phantom{\le\,}
    +
    C\iint_{Q_{R,S}}|F|^2|u|^{2m}\mathbf{1}_{\{|u|>k\}}
    \zeta^2\chi \,\dx\dt,
\end{align*}
for a constant $C=C(\nu,L)$. Arrived at this point, we let $\delta\downarrow 0$ in the cut-off function $\chi\equiv \chi_\delta$. In the limit we have
$\chi\downarrow\mathbf 1_{(t_o-S,\tau]}$ and $-\pl_t\chi$ converges to the Dirac mass at $\tau$. Consequently, the first integral on the left-hand side
converges to
$$ 
    \tfrac1{m+1}\int_{B_{R}(x_o)
    \times\{\tau\}}v\zeta^2\,\dx, 
$$
while the second integral  on the left-hand side can be estimated  from below by 
$$
\tfrac12\nu \iint_{B_R(x_o)\times (t_o-S,\tau]} |D\power{u}{m}|^2 \big(|\power{u}{m}|-k^m\big)_+^2\zeta^{2}\,\dx\dt.
$$ 
In the right-hand side integrals we bound the cut-off function $\chi$ from above by 1. After that, we take the supremum over $\tau\in (t_o-S,t_o)$ in the first integral  and let $\tau\uparrow t_o$ in the second one. In this way we obtain
\begin{align*}
    \tfrac{1}{m+1}&  \sup_{\tau\in (t_o-S,t_o)} 
    \int_{B_R \times\{\tau\}} v\, \zeta^2\,\dx
    +
    \tfrac12 \nu
    \iint_{Q_{R,S}}|D \power{u}{m}|^2
    \big(|\power{u}{m}|-k^m\big)_+^2 \zeta^2
    \,\dx\dt\\
    &\le 
    C\iint_{Q_{R,S}}|u|^{4m}\mathbf{1}_{\{|u|>k\}} |\nabla \zeta|^2
    \,\dx\dt +
    \tfrac{2}{m+1} \iint_{Q_{R,S}} v|\partial_t \zeta^2|\,\dx\dt\\
    &\phantom{\le\,}
    +
    C\iint_{Q_{R,S}}|F|^2|u|^{2m}\mathbf{1}_{\{|u|>k\}}
    \zeta^2 \,\dx\dt ,
\end{align*}
where $C=C(\nu,L)$. To treat the second term on the left and deduce the final form of the energy inequality we use Kato's inequality, that is $|Dw|\ge |\nabla |w||$. This  allows us to  estimate $|D\power{u}{m}|\ge |\nabla |u|^m|$, and with that 
\begin{align*}
    |\nabla |u|^m|^2&\big(|\power{u}{m}|-k^m\big)_{+}^2
    =
    |\nabla|u|^m|^2\big(|u|^m-k^m\big)_{+}^2 \\
    &=
    \big[|\nabla|u|^m|\big(|u|^m-k^m\big)_{+}\big]^2
    =
    \tfrac{1}{4}\big|\nabla\big(|u|^m-k^m\big)_{+}^2\big|^2.
\end{align*}
Hence
\begin{equation*}
     \iint_{Q_{R, S}}|D \power{u}{m}|^2
     \big(|u|^m-k^m\big)_{+}^2 \zeta^2\,\dx\dt
     \ge 
     \tfrac{1}{4} 
     \iint_{Q_{R, S}}\big|\nabla\big(|u|^m-k^m\big)_{+}^2\big|^2 \zeta^2\,\dx\dt,
\end{equation*}
which immediately leads to 
\begin{align*}
    \tfrac{1}{m+1}&  \sup_{\tau\in (t_o-S,t_o)} 
    \int_{B_R \times\{\tau\}} v\, \zeta^2\,\dx
    +
    \tfrac18 \nu
    \iint_{Q_{R, S}}\big|\nabla\big(|u|^m-k^m\big)_{+}^2\big|^2 \zeta^2\,\dx\dt\\
    &\le 
    C\iint_{Q_{R,S}}
    |u|^{4m}\mathbf{1}_{\{|u|>k\}}|\nabla \zeta|^2
    \,\dx\dt +
    \tfrac{2}{m+1} \iint_{Q_{R,S}} v|\partial_t \zeta^2|\,\dx\dt\\
    &\phantom{\le\,}
    +
    C\iint_{Q_{R,S}}|F|^2|u|^{2m}\mathbf{1}_{\{|u|>k\}}
    \zeta^2 \,\dx\dt.
\end{align*}
It remains to estimate $v$ from below and above. 
If we perform the change of variable $s^{\frac m{m+1}}=y$, it is apparent that 
\[
v=\frac{m+1}m\int_{k^m}^{|u|^m}(y-k^m)_+^2\, y^{\frac1m}\,\dy.
\]
On one hand,
\begin{equation*}
\int_{k^m}^{|u|^m}(y-k^m)_+^2\,y^{\frac1m}\,\dy\ge\int_{k^m}^{|u|^m}(y-k^m)_+^{2+\frac1m}\,\dy=\frac m{3m+1}(|u|^m-k^m)_+^{3+\frac1m},
\end{equation*}
which yields
\[
v\ge\frac{m+1}{3m+1}(|u|^m-k^m)_+^{3+\frac1m}.
\]
On the other hand,
\begin{equation*}
\int_{k^m}^{|u|^m}(y-k^m)_+^2\,y^{\frac1m}\,\dy\le |u|\int_{k^m}^{|u|^m}(y-k^m)_+^2\,\dy=\frac{1}3|u|(|u|^m-k^m)_+^{3},
\end{equation*}
so that
\[
v\le\frac{m+1}{3m}|u|(|u|^m-k^m)_+^{3}.
\]
Altogether this yields
\begin{align*}
    &  \sup_{\tau\in (t_o-S,t_o)} 
    \int_{B_R \times\{\tau\}} 
    \big(|u|^m-
    k^m\big)_+^{3+\frac{1}{m}}
    \, \zeta^2\,\dx
    +
    \iint_{Q_{R, S}}\big|\nabla\big(|u|^m-k^m\big)_{+}^2\big|^2 \zeta^2\,\dx\dt\\
    &\qquad\le 
    C\iint_{Q_{R,S}}|u|^{4m}\mathbf{1}_{\{|u|>k\}}|\nabla \zeta|^2
    \,\dx\dt
    +
    C \iint_{Q_{R,S}} 
    |u|
    \big(|u|^{m}-k^m\big)_+^3
    |\partial_t \zeta^2|\,\dx\dt\\
    &\qquad\phantom{\le\,} +
    C\iint_{Q_{R,S}}|F|^2|u|^{2m}\mathbf{1}_{\{|u|>k\}}
    \zeta^2 \,\dx\dt,
\end{align*}
for a constant $C=C(m,\nu,L)$.
\end{proof}

%\noindent{\bf Data Availability.} Data sharing not applicable to this article as no datasets were generated or analyzed during the current study.

%\

%\noindent{\bf Conflict of Interest.} The authors declare no conflict of interest.

\end{document}